\documentclass{article}
\usepackage{amsmath}
\usepackage{amssymb}
\begin{document}
\def\e#1\e{\begin{equation}#1\end{equation}}
\def\ea#1\ea{\begin{align}#1\end{align}}
\def\eq#1{{\rm(\ref{#1})}}
\newtheorem{thm}{Theorem}[section]
\newtheorem{lem}[thm]{Lemma}
\newtheorem{prop}[thm]{Proposition}
\newtheorem{conj}[thm]{Conjecture}
\newtheorem{cor}[thm]{Corollary}
\newenvironment{dfn}{\medskip\refstepcounter{thm}
\noindent{\bf Definition \thesection.\arabic{thm}\ }}{\medskip}
\newenvironment{ex}{\medskip\refstepcounter{thm}
\noindent{\bf Example \thesection.\arabic{thm}\ }}{\medskip}
\newenvironment{proof}[1][,]{\medskip\ifcat,#1
\noindent{\it Proof.\ }\else\noindent{\it Proof of #1.\ }\fi}
{\relax\unskip\nobreak ~\hfill$\square$\medskip}
\def\narrow{\par\addtolength{\leftskip}{20pt}\addtolength{\rightskip}{5pt}}
\def\oldmargins{\par\setlength{\leftskip}{0pt}\setlength{\rightskip}{0pt}}
\newcounter{quest}[section]
\newcounter{qa}[quest]
\newcounter{qi}[quest]
\newenvironment{question}{\narrow\setlength{\parindent}{0pt}\medskip
\refstepcounter{quest}\noindent\hbox to 0pt{\hss\bf
\thesection.\arabic{quest}\hskip .3em}}{\oldmargins\medskip}
\def\inext{\par\ifnum\value{qi}=0 \narrow\fi \medskip\addtocounter{qi}{1}
\noindent\hbox to 0pt{\hss\bf(\roman{qi})\hskip .3em}}
\def\anext{\par\ifnum\value{qa}=0 \narrow\fi \medskip\addtocounter{qa}{1}
\noindent\hbox to 0pt{\hss\bf(\alph{qa})\hskip .3em}}
\def\dim{\mathop{\rm dim}}
\def\Re{\mathop{\rm Re}}
\def\Im{\mathop{\rm Im}}
\def\Ker{\mathop{\rm Ker}}
\def\Coker{\mathop{\rm Coker}}
\def\ind{\mathop{\rm ind}}
\def\lind{\mathop{\text{\rm l-ind}}}
\def\sign{\mathop{\rm sign}}
\def\End{\mathop{\rm End}}
\def\Vect{\mathop{\rm Vect}}
\def\vol{\mathop{\rm vol}}
\def\Hol{{\textstyle\mathop{\rm Hol}}}
\def\hol{{\textstyle\mathop{\mathfrak{hol}}}}
\def\Ric{\mathop{\rm Ric}}
\def\Vol{\mathop{\rm Vol}}
\def\Hess{\mathop{\rm Hess}}
\def\Image{\mathop{\rm Image}}
\def\Tr{\mathop{\rm Tr}}
\def\Spin{\mathop{\rm Spin}}
\def\Sp{\mathop{\rm Sp}}
\def\GL{\mathop{\rm GL}}
\def\SO{\mathop{\rm SO}}
\def\O{\mathop{\rm O}}
\def\U{\mathbin{\rm U}}
\def\SL{\mathop{\rm SL}}
\def\SU{\mathop{\rm SU}}
\def\ge{\geqslant} 
\def\le{\leqslant} 
\def\cal{\mathcal}
\def\H{\mathbin{\mathbb H}}
\def\R{\mathbin{\mathbb R}}
\def\Z{\mathbin{\mathbb Z}}
\def\N{\mathbin{\mathbb N}}
\def\C{\mathbin{\mathbb C}}
\def\g{\mathbin{\mathfrak g}}
\def\h{\mathbin{\mathfrak h}}
\def\CP{\mathbb{CP}}
\def\u{\mathfrak{u}} 
\def\su{\mathfrak{su}} 
\def\so{\mathfrak{so}} 
\def\al{\alpha}
\def\be{\beta}
\def\ga{\gamma}
\def\de{\delta}
\def\ep{\epsilon}
\def\th{\theta}
\def\la{\lambda}
\def\ka{\kappa}
\def\vp{\varphi}
\def\si{\sigma}
\def\De{\Delta}
\def\La{\Lambda}
\def\Th{\Theta}
\def\Om{\Omega}
\def\Ga{\Gamma}
\def\Si{\Sigma}
\def\om{\omega}
\def\d{{\rm d}}
\def\pd{\partial}
\def\db{{\bar\partial}}
\def\ts{\textstyle}
\def\sst{\scriptscriptstyle}
\def\pha{\phantom}
\def\w{\wedge}
\def\lt{\ltimes}
\def\ti{\tilde}
\def\sm{\setminus}
\def\ov{\overline}
\def\ot{\otimes}
\def\bigot{\bigotimes}
\def\iy{\infty}
\def\ra{\rightarrow}
\def\t{\times}
\def\ha{{\textstyle\frac{1}{2}}}
\def\op{\oplus}
\def\nae{\nabla^{\sst E}}
\def\bs{\boldsymbol}
\def\ms#1{\vert#1\vert^2}
\def\bms#1{\bigl\vert#1\bigr\vert^2}
\def\md#1{\vert #1 \vert}
\def\bmd#1{\bigl\vert #1 \bigr\vert}
\def\an#1{\langle#1\rangle}
\title{Lectures on Calabi--Yau and \\ special Lagrangian geometry}
\author{Dominic Joyce \\ Lincoln College, Oxford, OX1 3DR}
\date{}
\maketitle

\section{Introduction}
\label{l1}

{\it Calabi--Yau $m$-folds} are compact Ricci-flat K\"ahler manifolds
$(M,J,g)$ of complex dimension $m$, with trivial canonical bundle 
$K_M$. Taken together, the complex structure $J$, K\"ahler metric $g$, 
and a holomorphic section $\Om$ of $K_M$ make up a rich, fairly 
rigid geometrical structure with very interesting properties --- for 
instance, Calabi--Yau $m$-folds occur in smooth, finite-dimensional 
moduli spaces with known dimension. 

Using Algebraic Geometry and Yau's solution of the Calabi Conjecture,
one can construct Calabi--Yau $m$-folds in huge numbers. String
Theorists (a species of theoretical physicist) are very interested
in Calabi--Yau 3-folds, and have made some extraordinary conjectures
about them, in the subject known as Mirror Symmetry.

{\it Special Lagrangian $m$-folds (SL\/ $m$-folds)} are a 
distinguished class of real $m$-dimensional minimal submanifolds 
$N$ in Calabi--Yau $m$-folds $M$. They are fairly rigid and 
well-behaved, so that compact SL $m$-folds occur in smooth 
moduli spaces of dimension $b^1(N)$, for instance. 

The existence of compact special Lagrangian $m$-folds in general 
Calabi--Yau $m$-folds is not well understood, but there are 
reasons to believe they are very abundant. They are important
in String Theory, and are expected to play a r\^ole in the
eventual explanation of Mirror Symmetry.

In this paper we aim to do three things. Firstly, we introduce
Calabi--Yau manifolds, going via Riemannian holonomy groups
and K\"ahler geometry. Secondly, we introduce special Lagrangian
submanifolds, going via calibrated geometry. Finally, we
survey an area of current research on the singularities
of SL $m$-folds, and its application to Mirror Symmetry 
and the SYZ Conjecture. Exercises are given at the end
of each section.

I hope the paper will be useful to graduate students in Geometry, 
String Theorists, and others who wish to learn the subject. Apart from 
the last two sections, the paper is intended as a straightforward 
exposition of standard material, to take the reader from a starting 
point of a good background in Differential Geometry as far as the
boundaries of current research in an exciting area.

Our approach to Calabi--Yau $m$-folds is resolutely Differential 
Geometric. That is, we regard them as smooth real manifolds 
equipped with a geometric structure. The alternative is to use 
Algebraic Geometry, regard them as complex algebraic varieties, 
and mostly forget about the metric $g$. Though very important, 
this side of the story will barely enter these notes, mainly 
because SL $m$-folds are not algebraic objects, and Algebraic 
Geometry has (so far) little to say about them.

I have chosen to introduce Calabi--Yau $m$-folds in the context 
of Riemannian holonomy groups, and special Lagrangian $m$-folds
in the context of calibrated geometry. Though these are perhaps
unnecessary diversions, I hope the reader will gain something
through understanding the wider horizon into which Calabi--Yau 
and special Lagrangian geometry fits. 

Also, it is my strong conviction that holonomy groups and 
calibrated geometry belong together as partner subjects,
and I want to take the opportunity to teach them together.
Though the field of Riemannian holonomy is now mature, the 
subject of calibrated submanifolds of Riemannian manifolds 
with special holonomy is really only beginning to be explored.

We begin in \S\ref{l2} with some background from Differential
Geometry, and define holonomy groups of connections and of
Riemannian metrics. Section \ref{l3} explains Berger's
classification of holonomy groups of Riemannian manifolds.
Section \ref{l4} discusses K\"ahler geometry and Ricci curvature
of K\"ahler manifolds and defines Calabi--Yau manifolds, and 
\S\ref{l5} sketches the proof of the Calabi Conjecture,
and how it is used to construct examples of Calabi--Yau
manifolds via Algebraic Geometry.

The second part of the paper begins in \S\ref{l6} with
an introduction to calibrated geometry. Section \ref{l7}
covers general properties of special Lagrangian $m$-folds
in $\C^m$, and \S\ref{l8} construction of examples. Section
\ref{l9} discusses compact SL $m$-folds in Calabi--Yau
$m$-folds, and \S\ref{l10} the singularities of SL $m$-folds.
Finally, \S\ref{l11} briefly introduces String Theory and
Mirror Symmetry, explains the {\it SYZ Conjecture}, and
summarizes some research on the singularities of special
Lagrangian fibrations.

Readers are warned that sections \ref{l8}, \ref{l10} and \ref{l11}
are unashamedly biased in favour of presenting the author's ideas
and opinions; this is not intended as an even-handed survey of the
whole field. Further, sections \ref{l10} and \ref{l11} are fairly
speculative, as I am setting out what I think the interesting 
problems in the field are, and where I would like it to go in the 
next few years.

The author's paper \cite{Joyc8} is a much shortened version of
this paper, containing only the special Lagrangian material,
roughly \S\ref{l7}--\S\ref{l11} below, slightly rewritten.
People already well-informed about Calabi--Yau geometry may
prefer to read~\cite{Joyc8}.
\medskip

\noindent{\it Acknowledgements.} These notes are based on lecture
courses given to the Summer School on Symplectic Geometry in
Nordfjordeid, Norway, in June 2001, and to the Clay Institute's
Summer School on the Global Theory of Minimal Surfaces, MSRI,
California, in July 2001. I would like to thank the organizers
of these for inviting me to speak.

Many people have helped me develop my ideas on special Lagrangian
geometry; amongst them I would particularly like to thank Robert 
Bryant, Mark Gross, Mark Haskins, Nigel Hitchin, Ian McIntosh,
Richard Thomas, and Karen Uhlenbeck.

\section{Introduction to holonomy groups}
\label{l2}

We begin by giving some background from differential and Riemannian
geometry, principally to establish notation, and move on to discuss 
connections on vector bundles, parallel transport, and the definition 
of holonomy groups. Some suitable reading for this section is my book 
\cite[\S 2 \& \S 3]{Joyc2}, and Kobayashi and
Nomizu~\cite[\S I--\S IV]{KoNo1}.

\subsection{Tensors and forms}
\label{l21}

Let $M$ be a smooth $n$-dimensional manifold, with tangent bundle 
$TM$ and cotangent bundle $T^*M$. Then $TM$ and $T^*M$ are 
{\it vector bundles} over $M$. If $E$ is a vector bundle over
$M$, we use the notation $C^\iy(E)$ for the vector space of
smooth sections of $E$. Elements of $C^\iy(TM)$ are called
{\it vector fields}, and elements of $C^\iy(T^*M)$ are called
1-{\it forms}. By taking tensor products of the vector bundles 
$TM$ and $T^*M$ we obtain the bundles of {\it tensors} on $M$. 
A {\it tensor} $T$ on $M$ is a smooth section of a bundle 
$\bigot^kTM\ot\bigot^lT^*M$ for some~$k,l\in\N$. 

It is convenient to write tensors using the {\it index notation}.
Let $U$ be an open set in $M$, and $(x^1,\ldots,x^n)$ coordinates
on $U$. Then at each point $x\in U$, $\frac{\pd}{\pd x^1},\ldots,
\frac{\pd}{\pd x^n}$ are a basis for $T_xU$. Hence, any vector field $v$ 
on $U$ may be uniquely written $v=\sum_{a=1}^nv^a\frac{\pd}{\pd x^a}$
for some smooth functions $v^1,\ldots,v^n:U\ra\R$\,. We denote $v$ by
$v^a$, which is understood to mean the collection of $n$ functions
$v^1,\ldots,v^n$, so that $a$ runs from 1 to $n$.

Similarly, at each $x\in U$, $\d x^1,\ldots,\d x^n$ are a basis for 
$T^*_xU$. Hence, any 1-form $\al$ on $U$ may be uniquely written 
$\al=\sum_{b=1}^n\al_b\d x^b$ for some smooth functions $\al_1,\ldots,
\al_n:U\ra\R$\,. We denote $\al$ by $\al_b$, where $b$ runs from 1 to $n$.
In the same way, a general tensor $T$ in $C^\iy(\bigot^kTM\ot\bigot^lT^*M)$ 
is written $T^{a_1\ldots a_k}_{b_1\ldots b_l}$, where
\begin{equation*}
T=\sum_{\substack{1\le a_i\le n,\ 1\le i\le k \\
1\le b_j\le n,\ 1\le j\le l}}T^{a_1\ldots a_k}_{b_1\ldots b_l}
{\ts\frac{\pd}{\pd x^{a_1}}\ot\cdots\frac{\pd}{\pd x^{a_k}}}\ot
\d x^{b_1}\ot\cdots\ot\d x^{b_l}.
\end{equation*}

The $k^{\rm th}$ exterior power of the cotangent bundle $T^*M$ is written 
$\La^kT^*M$. Smooth sections of $\La^kT^*M$ are called {\it $k$-forms}, 
and the vector space of $k$-forms is written $C^\iy(\La^kT^*M)$. They
are examples of tensors. In the index notation they are written 
$T_{b_1\ldots b_k}$, and are antisymmetric in the indices $b_1,\ldots,b_k$.
The {\it exterior product}\/ $\w$ and the {\it exterior derivative} $\d$ 
are important natural operations on forms. If $\al$ is a $k$-form and 
$\be$ an $l$-form then $\al\w\be$ is a $(k\!+\!l)$-form and $\d\al$ a 
$(k\!+\!1)$-form, which are given in index notation by
\begin{equation*}
(\al\w\be)_{a_1\ldots a_{k+l}}=\al_{[a_1\ldots a_k}\be_{a_{k+1}
\ldots a_{k+l}]}\quad\text{and}\quad (\d\al)_{a_1\ldots a_{k+1}}=
\frac{\pd}{\pd x^{[a_1}}\al_{a_2\ldots a_{k+1}]},
\end{equation*}
where $[\cdots]$ denotes antisymmetrization over the enclosed group 
of indices.

Let $v,w$ be vector fields on $M$. The {\it Lie bracket} $[v,w]$ of
$v$ and $w$ is another vector field on $M$, given in index notation by
\e
[v,w]^a=v^b\frac{\pd w^a}{\pd x^b}-w^b\frac{\pd v^a}{\pd x^b}.
\label{liebracketeq}
\e
Here we have used the {\it Einstein summation convention}, that is,
the repeated index $b$ on the right hand side is summed from 1 to $n$.
The important thing about this definition is that it is independent
of choice of coordinates $(x^1,\ldots,x^n)$.

\subsection{Connections on vector bundles and curvature}
\label{l22}

Let $M$ be a manifold, and $E\ra M$ a vector bundle. A {\it connection} 
$\nae$ on $E$ is a linear map $\nae:C^\iy(E)\ra C^\iy(E\ot T^*M)$ 
satisfying the condition
\begin{equation*}
\nae(\al\,e)=\al\nae e+e\ot\d\al,
\end{equation*}
whenever $e\in C^\iy(E)$ is a smooth section of $E$ and $\al$ is a 
smooth function on~$M$. 

If $\nae$ is such a connection, $e\in 
C^\iy(E)$, and $v\in C^\iy(TM)$ is a vector field, then we write 
$\nae_ve=v\cdot\nae e\in C^\iy(E)$, where `$\cdot$' contracts
together the $TM$ and $T^*M$ factors in $v$ and $\nae e$. Then
if $v\in C^\iy(TM)$ and $e\in C^\iy(E)$ and $\al,\be$ are smooth
functions on $M$, we have
\begin{equation*}
\nae_{\al v}(\be e)=\al\be\nae_ve+\al(v\cdot\be)e.
\end{equation*}
Here $v\cdot\be$ is the Lie derivative of $\be$ by $v$. It is a 
smooth function on $M$, and could also be written~$v\cdot\d\be$.

There exists a unique, smooth section $R(\nae)\in C^\iy\bigl(\End(E)
\ot\La^2T^*M\bigr)$ called the {\it curvature} of $\nae$, that satisfies
the equation
\e
R(\nae)\cdot(e\ot v\w w)=\nae_v\nae_we-\nae_w\nae_ve-\nae_{[v,w]}e
\label{redefeq}
\e
for all $v,w\in C^\iy(TM)$ and $e\in C^\iy(E)$, where $[v,w]$ is the
Lie bracket of~$v,w$.

Here is one way to understand the curvature of $\nae$. Define 
$v_i=\pd/\pd x^i$ for $i=1,\ldots,n$. Then $v_i$ is a vector field 
on $U$, and $[v_i,v_j]=0$. Let $e$ be a smooth section of $E$. Then 
we may interpret $\nae_{v_i}e$ as a kind of {\it partial derivative} 
$\pd e/\pd x^i$ of $e$. Equation \eq{redefeq} then implies that
\e
R(\nae)\cdot(e\ot v_i\w v_j)=
\frac{\pd^2 e}{\pd x^i\pd x^j}-\frac{\pd^2 e}{\pd x^j\pd x^i}.
\label{reijeq}
\e
Thus, {\it the curvature $R(\nae)$ measures how much partial 
derivatives in $E$ fail to commute}.

Now let $\nabla$ be a connection on the tangent bundle $TM$ of
$M$, rather than a general vector bundle $E$. Then there is a 
unique tensor $T=T^a_{bc}$ in $C^\iy\bigl(TM\ot\La^2T^*M\bigr)$ 
called the {\it torsion} of $\nabla$, satisfying 
\begin{equation*}
T\cdot(v\w w)=\nabla_vw-\nabla_wv-[v,w]
\quad\text{for all $v,w\in C^\iy(TM)$.}
\end{equation*}
A connection $\nabla$ with zero torsion is called {\it torsion-free}.
Torsion-free connections have various useful properties, so we
usually restrict attention to torsion-free connections on~$TM$.

A connection $\nabla$ on $TM$ extends naturally to connections
on all the bundles of tensors $\bigot^kTM\ot\bigot^lT^*M$ for 
$k,l\in\N$, which we will also write $\nabla$. That is, we
can use $\nabla$ to differentiate not just vector fields,
but any tensor on~$M$.

\subsection{Parallel transport and holonomy groups}
\label{l23}

Let $M$ be a manifold, $E\ra M$ a vector bundle over $M$, and 
$\nae$ a connection on $E$. Let $\ga:[0,1]\ra M$ be a smooth 
curve in $M$. Then the pull-back $\ga^*(E)$ of $E$ to $[0,1]$ 
is a vector bundle over $[0,1]$ with fibre $E_{\ga(t)}$ over 
$t\in[0,1]$, where $E_x$ is the fibre of $E$ over~$x\in M$. 
The connection $\nae$ pulls back under $\ga$ to give a 
connection on $\ga^*(E)$ over~$[0,1]$. 

\begin{dfn} Let $M$ be a manifold, $E$ a vector bundle over $M$, 
and $\nae$ a connection on $E$. Suppose $\ga:[0,1]\ra M$ is 
(piecewise) smooth, with $\ga(0)=x$ and $\ga(1)=y$, where 
$x,y\in M$. Then for each $e\in E_x$, there exists a unique 
smooth section $s$ of $\ga^*(E)$ satisfying $\nae_{\dot\ga(t)}s(t)=0$ 
for $t\in[0,1]$, with $s(0)=e$. Define $P_\ga(e)=s(1)$. Then 
$P_\ga:E_x\ra E_y$ is a well-defined linear map, called the 
{\it parallel transport map}. 
\end{dfn}

We use parallel transport to define the {\it holonomy group} of~$\nae$.

\begin{dfn} Let $M$ be a manifold, $E$ a vector bundle over $M$, 
and $\nae$ a connection on $E$. Fix a point $x\in M$. We say that 
$\ga$ is a {\it loop based at\/ $x$} if $\ga:[0,1]\ra M$ is a 
piecewise-smooth path with $\ga(0)=\ga(1)=x$. The parallel 
transport map $P_\ga:E_x\ra E_x$ is an invertible linear map, 
so that $P_\ga$ lies in $\GL(E_x)$, 
the group of invertible linear transformations of $E_x$. Define 
the {\it holonomy group} $\Hol_x(\nae)$ of $\nae$ based at $x$ to be
\e
\Hol_x(\nae)=\bigl\{P_\ga:
\hbox{$\ga$ is a loop based at $x$}\bigr\}\subset\GL(E_x).
\e
\end{dfn}

The holonomy group has the following important properties.

\begin{itemize}
\item It is a {\it Lie subgroup} of $\GL(E_x)$. To show that 
$\Hol_x(\nae)$ is a subgroup of $\GL(E_x)$, let $\ga,\de$ be 
loops based at $x$, and define loops $\ga\de$ and $\ga^{-1}$ by
\begin{equation*}
\ga\de(t)=\begin{cases} \de(2t) & t\in[0,\ha] \\
\ga(2t-1) & t\in[\ha,1] \end{cases}\quad\text{and}\quad
\ga^{-1}(t)=\ga(1-t)\quad\text{for $t\in[0,1]$.}
\end{equation*}
Then $P_{\ga\de}=P_\ga\circ P_\de$ and $P_{\ga^{-1}}=P_\ga^{-1}$,
so $\Hol_x(\nae)$ is closed under products and inverses.
\item It is {\it independent of basepoint} $x\in M$, in the following
sense. Let $x,y\in M$, and let $\ga:[0,1]\ra M$ be a smooth path from
$x$ to $y$. Then $P_\ga:E_x\ra E_y$, and $\Hol_x(\nae)$ and 
$\Hol_y(\nae)$ satisfy $\Hol_y(\nae)=P_\ga\Hol_x(\nae)P_\ga^{-1}$.

Suppose $E$ has fibre $\R^k$, so that $\GL(E_x)\cong\GL(k,\R)$. Then 
we may regard $\Hol_x(\nae)$ as a subgroup of $\GL(k,\R)$ defined up 
to conjugation, and it is then independent of basepoint~$x$.

\item If $M$ is simply-connected, then $\Hol_x(\nae)$ is connected.
To see this, note that any loop $\ga$ based at $x$ can be continuously 
shrunk to the constant loop at $x$. The corresponding family of parallel
transports is a continuous path in $\Hol_x(\nae)$ joining $P_\ga$ to
the identity.
\end{itemize}

The holonomy group of a connection is closely related to its 
curvature. Here is one such relationship. As $\Hol_x(\nae)$ 
is a Lie subgroup of $\GL(E_x)$, it has a {\it Lie algebra} 
$\hol_x(\nae)$, which is a Lie subalgebra of $\End(E_x)$. It
can be shown that the curvature $R(\nae)_x$ at $x$ lies in
the linear subspace $\hol_x(\nae)\ot\La^2T_x^*M$ of $\End(E_x)
\ot\La^2T_x^*M$. Thus, {\it the holonomy group of a connection 
places a linear restriction upon its curvature}.

Now let $\nabla$ be a connection on $TM$. Then from \S\ref{l22},
$\nabla$ extends to connections on all the tensor bundles $\bigot^k
TM\ot\bigot^lT^*M$. We call a tensor $S$ on $M$ {\it constant} if
$\nabla S=0$. The constant tensors on $M$ are determined by the 
holonomy group~$\Hol(\nabla)$.

\begin{thm} Let\/ $M$ be a manifold, and\/ $\nabla$ a connection 
on $TM$. Fix $x\in M$, and let\/ $H=\Hol_x(\nabla)$. Then $H$ acts
naturally on the tensor powers $\bigot^kT_xM\ot\bigot^lT_x^*M$. 
Suppose $S\in C^\iy\bigl(\bigot^kTM\ot\bigot^lT^*M\bigr)$ is
a constant tensor. Then $S\vert_x$ is fixed by the action of\/ 
$H$ on $\bigot^kT_xM\ot\bigot^lT_x^*M$. Conversely, if\/ 
$S\vert_x\in\bigot^kT_xM\ot\bigot^lT_x^*M$ is fixed by $H$, 
it extends to a unique constant tensor $S\in C^\iy\bigl(
\bigot^kTM\ot\bigot^lT^*M\bigr)$.
\label{l2thm1}
\end{thm}

The main idea in the proof is that if $S$ is a constant tensor and
$\ga:[0,1]\ra M$ is a path from $x$ to $y$, then $P_\ga(S\vert_x)=
S\vert_y$. Thus, constant tensors are invariant under parallel transport.

\subsection{Riemannian metrics and the Levi-Civita connection}
\label{l24}

Let $g$ be a {\it Riemannian metric} on $M$. We refer to the pair
$(M,g)$ as a {\it Riemannian manifold}. Here $g$ is a tensor in 
$C^\iy(S^2T^*M)$, so that $g=g_{ab}$ in index notation with 
$g_{ab}=g_{ba}$. There exists a unique, torsion-free connection 
$\nabla$ on $TM$ with $\nabla g=0$, called the 
{\it Levi-Civita connection}, which satisfies
\begin{align*}
2g(\nabla_uv,w)=\,&u\cdot g(v,w)+v\cdot g(u,w)-w\cdot g(u,v)\\
&+g([u,v],w)-g([v,w],u)-g([u,w],v)
\end{align*}
for all $u,v,w\in C^\iy(TM)$. This result is known as the 
{\it fundamental theorem of Riemannian geometry}.

The curvature $R(\nabla)$ of the Levi-Civita connection is a tensor 
$R^a_{\pha{a}bcd}$ on $M$. Define $R_{abcd}=g_{ae}R^e_{\pha{e}bcd}$. 
We shall refer to both $R^a_{\pha{a}bcd}$ and $R_{abcd}$ as the 
{\it Riemann curvature} of $g$. The following theorem gives a number 
of symmetries of $R_{abcd}$. Equations \eq{r1bianchi} and 
\eq{r2bianchi} are known as the {\it first\/} and {\it second 
Bianchi identities}, respectively.

\begin{thm} Let\/ $(M,g)$ be a Riemannian manifold, $\nabla$ the 
Levi-Civita connection of\/ $g$, and\/ $R_{abcd}$ the Riemann 
curvature of\/ $g$. Then $R_{abcd}$ and\/ $\nabla_eR_{abcd}$ 
satisfy the equations
\ea
&R_{abcd}=-R_{abdc}=-R_{bacd}=R_{cdab},
\label{rsymm}\\
&R_{abcd}+R_{adbc}+R_{acdb}=0,
\label{r1bianchi}\\
\hbox{and}\qquad
&\nabla_eR_{abcd}+\nabla_cR_{abde}+\nabla_dR_{abec}=0.
\label{r2bianchi}
\ea
\label{l2thm2}
\end{thm}

Let $(M,g)$ be a Riemannian manifold, with Riemann curvature 
$R^a_{\pha{a}bcd}$. The {\it Ricci curvature} of $g$ is 
$R_{ab}=R^c_{\pha{c}acb}$. It is a component of the full Riemann 
curvature, and satisfies $R_{ab}=R_{ba}$. We say that $g$ is 
{\it Einstein} if $R_{ab}=\la g_{ab}$ for some constant $\la\in\R$, 
and {\it Ricci-flat} if $R_{ab}=0$. Einstein and Ricci-flat metrics 
are of great importance in mathematics and physics.

\subsection{Riemannian holonomy groups}
\label{l25}

Let $(M,g)$ be a Riemannian manifold. We define the {\it holonomy 
group} $\Hol_x(g)$ of $g$ to be the holonomy group $\Hol_x(\nabla)$
of the Levi-Civita connection $\nabla$ of $g$, as in \S\ref{l23}. 
Holonomy groups of Riemannian metrics, or {\it Riemannian holonomy 
groups}, have stronger properties than holonomy groups of connections 
on arbitrary vector bundles. We shall explore some of these.

Firstly, note that $g$ is a {\it constant tensor} as $\nabla g=0$,
so $g$ is invariant under $\Hol(g)$ by Theorem \ref{l2thm1}. That
is, $\Hol_x(g)$ lies in the subgroup of $\GL(T_xM)$ which preserves
$g\vert_x$. This subgroup is isomorphic to $\O(n)$. Thus, $\Hol_x(g)$
may be regarded as a {\it subgroup of\/ $\O(n)$ defined up to 
conjugation}, and it is then independent of $x\in M$, so we will often 
write it as $\Hol(g)$, dropping the basepoint~$x$.

Secondly, the holonomy group $\Hol(g)$ constrains the Riemann 
curvature of $g$, in the following way. The Lie algebra 
$\hol_x(\nabla)$ of $\Hol_x(\nabla)$ is a vector subspace of 
$T_xM\ot T^*_xM$. From \S\ref{l23}, we have $R^a_{\pha{a}bcd}
\vert_x\in\hol_x(\nabla)\ot\La^2T_x^*M$.

Use the metric $g$ to identify $T_xM\ot T^*_xM$ and $\ot^2T^*_xM$, 
by equating $T^a_{\pha{a}b}$ with $T_{ab}=g_{ac}T^c_{\pha{c}b}$. This 
identifies $\hol_x(\nabla)$ with a vector subspace of $\ot^2T^*_xM$ 
that we will write as $\hol_x(g)$. Then $\hol_x(g)$ lies in $\La^2T^*_xM$,
and $R_{abcd}\vert_x\in\hol_x(g)\ot\La^2T_x^*M$. Applying the symmetries 
\eq{rsymm} of $R_{abcd}$, we have:

\begin{thm} Let\/ $(M,g)$ be a Riemannian manifold with Riemann curvature 
$R_{abcd}$. Then $R_{abcd}$ lies in the vector subspace 
$S^2\hol_x(g)$ in $\La^2T_x^*M\ot\La^2T_x^*M$ at each\/~$x\in M$.
\label{l2thm3}
\end{thm}

Combining this theorem with the Bianchi identities, \eq{r1bianchi} 
and \eq{r2bianchi}, gives strong restrictions on the curvature tensor 
$R_{abcd}$ of a Riemannian metric $g$ with a prescribed holonomy group
$\Hol(g)$. These restrictions are the basis of the classification 
of Riemannian holonomy groups, which will be explained in~\S\ref{l3}.

\subsection{Exercises}
\label{l26}

\begin{question}Let $M$ be a manifold and $u,v,w$ be vector fields on 
$M$. The {\it Jacobi identity} for the Lie bracket of vector fields is
\begin{equation*}
[u,[v,w]]+[v,[w,u]]+[w,[u,v]]=0.
\end{equation*}
Prove the Jacobi identity in coordinates $(x^1,\dots,x^n)$ on a 
coordinate patch $U$. Use the coordinate expression \eq{liebracketeq}
for the Lie bracket of vector fields.
\end{question}

\begin{question}In \S\ref{l23} we explained that if $M$ is a manifold, 
$E\ra M$ a vector bundle and $\nae$ a connection, then $\Hol(\nae)$ 
is connected when $M$ is simply-connected. If $M$ is not 
simply-connected, what is the relationship between the 
fundamental group $\pi_1(M)$ and $\Hol(\nae)$?
\end{question}

\begin{question}Work out your own proof of Theorem~\ref{l2thm1}.
\end{question}

\section{Berger's classification of holonomy groups}
\label{l3}

Next we describe Berger's classification of Riemannian holonomy 
groups, and briefly discuss the possibilities in the classification. 
Some references for the material of this section are my book 
\cite[\S 3]{Joyc2} and Kobayashi and Nomizu \cite[\S XI]{KoNo2}.
Berger's original paper is \cite{Berg}, but owing to language and
notation most will now find it difficult to read.

\subsection{Reducible Riemannian manifolds}
\label{l31}

Let $(P,g)$ and $(Q,h)$ be Riemannian manifolds with positive 
dimension, and $P\t Q$ the product manifold. Then at each $(p,q)$ 
in $P\t Q$ we have $T_{(p,q)}(P\t Q)\cong T_pP\op T_qQ$. Define 
the {\it product metric} $g\t h$ on $P\t Q$ by $g\t h\vert_{(p,q)}
=g\vert_p+h\vert_q$ for all $p\in P$ and $q\in Q$. We call 
$(P\t Q,g\t h)$ a {\it Riemannian product}.

A Riemannian manifold $(M,g')$ is said to be {\it (locally) reducible} 
if every point has an open neighbourhood isometric to a Riemannian
product $(P\t Q,g\t h)$, and {\it irreducible} if it is not locally 
reducible. It is easy to show that the holonomy of a product metric 
$g\t h$ is the product of the holonomies of $g$ and~$h$.

\begin{prop} If\/ $(P,g)$ and\/ $(Q,h)$ are Riemannian 
manifolds, then $\Hol(g\t h)=\Hol(g)\t\Hol(h)$.
\label{l3prop1}
\end{prop}

Here is a kind of converse to this.

\begin{thm} Let\/ $M$ be an $n$-manifold, and\/ $g$ an 
irreducible Riemannian metric on $M$. Then the representation 
of\/ $\Hol(g)$ on $\R^n$ is irreducible.
\label{l3thm1}
\end{thm}

To prove the theorem, suppose $\Hol(g)$ acts reducibly on $\R^n$, 
so that $\R^n$ is the direct sum of representations $\R^k$, $\R^l$
of $\Hol(g)$ with $k,l>0$. Using parallel transport, one can define
a splitting $TM=E\op F$, where $E,F$ are vector subbundles with
fibres $\R^k,\R^l$. These vector subbundles are {\it integrable},
so locally $M\cong P\t Q$ with $E=TP$ and $F=TQ$. One can then
show that the metric on $M$ is the product of metrics on $P$ and 
$Q$, so that $g$ is locally reducible.

\subsection{Symmetric spaces}
\label{l32}

Next we discuss Riemannian symmetric spaces. 

\begin{dfn} A Riemannian manifold $(M,g)$ is said to be a
{\it symmetric space} if for every point $p\in M$ there 
exists an isometry $s_p:M\ra M$ that is an involution
(that is, $s_p^2$ is the identity), such that $p$ is an
isolated fixed point of~$s_p$.
\end{dfn}

Examples include $\R^n$, spheres ${\cal S}^n$, projective spaces
$\CP^m$ with the Fubini--Study metric, and so on. Symmetric spaces
have a transitive group of isometries.

\begin{prop} Let $(M,g)$ be a connected, simply-connected 
symmetric space. Then $g$ is complete. Let $G$ be the group 
of isometries of\/ $(M,g)$ generated by elements of the form 
$s_q\circ s_r$ for $q,r\in M$. Then $G$ is a connected Lie 
group acting transitively on $M$. Choose $p\in M$, and let\/ 
$H$ be the subgroup of\/ $G$ fixing $p$. Then $H$ is a closed, 
connected Lie subgroup of\/ $G$, and\/ $M$ is the homogeneous 
space~$G/H$. 
\label{l3prop2}
\end{prop}

Because of this, symmetric spaces can be classified completely using
the theory of Lie groups. This was done in 1925 by \'Elie Cartan.
From Cartan's classification one can quickly deduce the list of 
holonomy groups of symmetric spaces.

A Riemannian manifold $(M,g)$ is called {\it locally symmetric} if
every point has an open neighbourhood isometric to an open set in a
symmetric space, and {\it nonsymmetric} if it is not locally 
symmetric. It is a surprising fact that Riemannian manifolds are locally
symmetric if and only if they have {\it constant curvature}.

\begin{thm} Let\/ $(M,g)$ be a Riemannian manifold, with
Levi-Civita connection $\nabla$ and Riemann curvature $R$. Then 
$(M,g)$ is locally symmetric if and only if\/~$\nabla R=0$.
\label{l3thm2}
\end{thm}

\subsection{Berger's classification}
\label{l33}

In 1955, Berger proved the following result.

\begin{thm}[Berger] Suppose $M$ is a simply-connected 
manifold of dimension $n$, and that\/ $g$ is a Riemannian 
metric on $M$, that is irreducible and nonsymmetric. 
Then exactly one of the following seven cases holds.
\begin{itemize}
\setlength{\parsep}{0pt}
\setlength{\itemsep}{0pt}
\item[{\rm(i)}] $\Hol(g)=\SO(n)$,
\item[{\rm(ii)}] $n=2m$ with\/ $m\ge 2$, and\/ 
$\Hol(g)=\U(m)$ in $\SO(2m)$,
\item[{\rm(iii)}] $n=2m$ with\/ $m\ge 2$, and\/ 
$\Hol(g)=\SU(m)$ in $\SO(2m)$,
\item[{\rm(iv)}] $n=4m$ with\/ $m\ge 2$, and\/ 
$\Hol(g)=\Sp(m)$ in $\SO(4m)$,
\item[{\rm(v)}] $n=4m$ with\/ $m\ge 2$, and\/ 
$\Hol(g)=\Sp(m)\Sp(1)$ in $\SO(4m)$,
\item[{\rm(vi)}] $n=7$ and\/ $\Hol(g)=G_2$ in $\SO(7)$, or
\item[{\rm(vii)}] $n=8$ and\/ $\Hol(g)=\Spin(7)$ in $\SO(8)$.
\end{itemize}
\label{l3thm3}
\end{thm}

Notice the three simplifying assumptions on $M$ and $g$: that $M$ is 
simply-connected, and $g$ is irreducible and nonsymmetric. Each 
condition has consequences for the holonomy group~$\Hol(g)$.
\begin{itemize}
\item As $M$ is simply-connected, $\Hol(g)$ is connected,
from \S\ref{l23}.
\item As $g$ is irreducible, $\Hol(g)$ acts irreducibly on
$\R^n$ by Theorem~\ref{l3thm1}.
\item As $g$ is nonsymmetric, $\nabla R\not\equiv 0$ by 
Theorem~\ref{l3thm2}.
\end{itemize}
The point of the third condition is that there are some
holonomy groups $H$ which can {\it only} occur for metrics
$g$ with $\nabla R=0$, and these holonomy groups are excluded
from the theorem.

One can remove the three assumptions, at the cost of making
the list of holonomy groups much longer. To allow $g$ to be 
symmetric, we must include the holonomy groups of Riemannian 
symmetric spaces, which are known from Cartan's classification.
To allow $g$ to be reducible, we must include all products 
of holonomy groups already on the list. To allow $M$ not 
simply-connected, we must include non-connected Lie groups
whose identity components are already on the list.

Berger proved that the groups on his list were the only
possibilities, but he did not show whether the groups actually
do occur as holonomy groups. It is now known (but this took
another thirty years to find out) that all of the groups on 
Berger's list do occur as the holonomy groups of irreducible, 
nonsymmetric metrics.

\subsection{A sketch of the proof of Berger's Theorem}
\label{l34}

Let $(M,g)$ be a Riemannian $n$-manifold with $M$ simply-connected
and $g$ irreducible and nonsymmetric, and let $H=\Hol(g)$. Then it
is known that $H$ is a closed, connected Lie subgroup of $\SO(n)$.
The classification of such subgroups follows from the classification
of Lie groups. Berger's method was to take the list of all closed, 
connected Lie subgroups $H$ of $\SO(n)$, and apply two tests to each 
possibility to find out if it could be a holonomy group. The only 
groups $H$ which passed both tests are those in the Theorem~\ref{l3thm3}. 

Berger's tests are algebraic and involve the curvature tensor. 
Suppose $R_{abcd}$ is the Riemann curvature of a metric $g$ with 
$\Hol(g)=H$, and let $\h$ be the Lie algebra of $H$. Then Theorem 
\ref{l2thm2} shows that $R_{abcd}\in S^2\h$, and the first Bianchi 
identity \eq{r1bianchi} applies. 

If $\h$ has large codimension in $\so(n)$, then the vector space 
$\mathfrak{R}^H$ of elements of $S^2\h$ satisfying \eq{r1bianchi} will 
be small, or even zero. But the {\it Ambrose--Singer Holonomy Theorem} 
shows that $\mathfrak{R}^H$ must be big enough to generate $\h$, in a 
certain sense. For many of the candidate groups $H$ this does not hold, 
and so $H$ cannot be a holonomy group. This is the first test. 

Now $\nabla_eR_{abcd}$ lies in $(\R^n)^*\ot\mathfrak{R}^H$, and also 
satisfies the second Bianchi identity \eq{r2bianchi}. Frequently these 
requirements imply that $\nabla R=0$, so that $g$ is locally symmetric. 
Therefore we may exclude such $H$, and this is Berger's second test. 

\subsection{The groups on Berger's list}
\label{l35}

Here are some brief remarks about each group on Berger's list.

\begin{list}{}{
\setlength{\leftmargin}{25pt}
\setlength{\labelwidth}{20pt}
\setlength{\parsep}{1pt}
\setlength{\itemsep}{1pt}}
\item[(i)] $\SO(n)$ is the holonomy group of generic Riemannian metrics. 
\item[(ii)] Riemannian metrics $g$ with $\Hol(g)\subseteq\U(m)$ are 
called {\it K\"ahler metrics}. K\"ahler metrics are a natural class of 
metrics on complex manifolds, and generic K\"ahler metrics on a 
given complex manifold have holonomy $\U(m)$. 
\item[(iii)] Metrics $g$ with $\Hol(g)\subseteq\SU(m)$ are called
{\it Calabi--Yau metrics}. Since $\SU(m)$ is a subgroup of $\U(m)$, 
all Calabi--Yau metrics are K\"ahler. If $g$ is K\"ahler and $M$
is simply-connected, then $\Hol(g)\subseteq\SU(m)$ if and only 
if $g$ is Ricci-flat. Thus Calabi--Yau metrics are locally the 
same as Ricci-flat K\"ahler metrics.
\item[(iv)] Metrics $g$ with $\Hol(g)\subseteq\Sp(m)$ are called
{\it hyperk\"ahler}. As $\Sp(m)\subseteq\SU(2m)\subset\U(2m)$,
hyperk\"ahler metrics are Ricci-flat and K\"ahler. 
\item[(v)] Metrics $g$ with holonomy group $\Sp(m)\Sp(1)$ for $m\ge 2$
are called {\it quaternionic K\"ahler}. (Note that quaternionic K\"ahler 
metrics are not in fact K\"ahler.) They are Einstein, but not Ricci-flat. 
\item[(vi) and (vii)] The holonomy groups $G_2$ and $\Spin(7)$ 
are called the {\it exceptional holonomy groups}. Metrics with these
holonomy groups are Ricci-flat.
\end{list}

The groups can be understood in terms of the four {\it division 
algebras}: the {\it real numbers} $\R$, the {\it complex numbers} 
$\C$, the {\it quaternions} $\H$, and the {\it octonions} or 
{\it Cayley numbers} $\mathbb O$. 
\begin{itemize}
\item $\SO(n)$ is a group of automorphisms of $\R^n$. 
\item $\U(m)$ and $\SU(m)$ are groups of automorphisms of $\C^m$ 
\item $\Sp(m)$ and $\Sp(m)\Sp(1)$ are automorphism groups of $\H^m$. 
\item $G_2$ is the automorphism group of $\Im\mathbb{O}\cong\R^7$. 
$\Spin(7)$ is a group of automorphisms of ${\mathbb O}\cong\R^8$.
\end{itemize}

Here are three ways in which we can gather together the holonomy 
groups on Berger's list into subsets with common features.

\begin{itemize}
\item The {\it K\"ahler holonomy groups} are $\U(m)$, $\SU(m)$ and
$\Sp(m)$. Any Riemannian manifold with one of these holonomy groups 
is a K\"ahler manifold, and thus a complex manifold. 

\item The {\it Ricci-flat holonomy groups} are $\SU(m)$, 
$\Sp(m)$, $G_2$ and $\Spin(7)$. Any metric with one of these 
holonomy groups is Ricci-flat. This follows from the effect of
holonomy on curvature discussed in \S\ref{l25} and \S\ref{l34}:
if $H$ is one of these holonomy groups and $R_{abcd}$ any curvature
tensor lying in $S^2\h$ and satisfying \eq{r1bianchi}, then
$R_{abcd}$ has zero Ricci component.

\item The {\it exceptional holonomy groups} are $G_2$ and $\Spin(7)$. 
They are the exceptional cases in Berger's classification, and they 
are rather different from the other holonomy groups. 
\end{itemize}

\subsection{Exercises}
\label{l36}

\begin{question}Work out your own proofs of Proposition \ref{l3prop1} 
and (harder) Theorem~\ref{l3thm1}.
\end{question}

\begin{question}Suppose that $(M,g)$ is a simply-connected Ricci-flat 
K\"ahler manifold of complex dimension 4. What are the possibilities 
for~$\Hol(g)$?

[You may use the fact that the only simply-connected Ricci-flat 
symmetric spaces are $\R^n$, $n\in\N$.]
\end{question}

\section{K\"ahler geometry and holonomy}
\label{l4}

We now focus our attention on K\"ahler geometry, and the Ricci curvature 
of K\"ahler manifolds. This leads to the definition of {\it Calabi--Yau 
manifolds}, compact Ricci-flat K\"ahler manifolds with holonomy $\SU(m)$. 
A reference for this section is my book \cite[\S 4, \S 6]{Joyc2}.

\subsection{Complex manifolds}
\label{l41}

We begin by defining {\it complex manifolds} $M$. The usual definition
of complex manifolds involves an atlas of complex coordinate patches
covering $M$, whose transition functions are holomorphic. However,
for our purposes we need a more differential geometric definition,
involving a tensor $J$ on $M$ called a {\it complex structure}.

Let $M$ be a real manifold of dimension $2m$. An {\it almost complex 
structure} $J$ on $M$ is a tensor $J_a^b$ on $M$ satisfying 
$J_a^bJ_b^c=-\de_a^c$. For each vector field $v$ on $M$ define $Jv$ 
by $(Jv)^b=J_a^bv^a$. Then $J^2=-1$, so $J$ gives each tangent space 
$T_pM$ the structure of a {\it complex vector space}.

We can associate a tensor $N=N^a_{bc}$ to $J$, called the
{\it Nijenhuis tensor}, which satisfies
\begin{equation*}
N^a_{bc}v^bw^c=\bigl([v,w]+J\bigl([Jv,w]+[v,Jw]\bigr)-[Jv,Jw]\bigr)^a
\end{equation*}
for all vector fields $v,w$ on $M$, where $[\,,\,]$ is the Lie 
bracket of vector fields. The almost complex structure $J$ is
called a {\it complex structure} if $N\equiv 0$. A {\it complex
manifold} $(M,J)$ is a manifold $M$ with a complex structure~$J$.

Here is why this is equivalent to the usual definition. A smooth 
function $f:M\ra\C$ is called {\it holomorphic} if $J_a^b(\d f)_b
\equiv i(\d f)_a$ on $M$. These are called the {\it Cauchy--Riemann 
equations}. It turns out that the Nijenhuis tensor $N$ is the 
obstruction to the existence of holomorphic functions. If $N\equiv 0$ 
there are many holomorphic functions locally, enough to form a set
of holomorphic coordinates around every point. 

\subsection{K\"ahler manifolds}
\label{l42}

Let $(M,J)$ be a complex manifold, and let $g$ be a Riemannian metric on 
$M$. We call $g$ a {\it Hermitian metric} if $g(v,w)=g(Jv,Jw)$ for all 
vector fields $v,w$ on $M$, or $g_{ab}=J_a^cJ_b^dg_{cd}$ in index 
notation. When $g$ is Hermitian, define the {\it Hermitian form} $\om$ 
of $g$ by $\om(v,w)=g(Jv,w)$ for all vector fields $v,w$ on $M$, or 
$\om_{ac}=J_a^bg_{bc}$ in index notation. Then $\om$ is a (1,1)-form, 
and we may reconstruct $g$ from $\om$ by~$g(v,w)=\om(v,Jw)$.

A Hermitian metric $g$ on a complex manifold $(M,J)$ is called {\it 
K\"ahler} if one of the following three equivalent conditions holds:
\begin{itemize}
\setlength{\parsep}{0pt}
\setlength{\itemsep}{0pt}
\item[{\rm(i)}] $\d\om=0$,
\item[{\rm(ii)}] $\nabla J=0$, or
\item[{\rm(iii)}] $\nabla\om=0$,
\end{itemize}
where $\nabla$ is the Levi-Civita connection of $g$. We then call
$(M,J,g)$ a {\it K\"ahler manifold}. K\"ahler metrics are a natural 
and important class of metrics on complex manifolds.

By parts (ii) and (iii), if $g$ is K\"ahler then $J$ and $\om$ are
{\it constant tensors} on $M$. Thus by Theorem \ref{l2thm1}, the 
holonomy group $\Hol(g)$ must preserve a complex structure $J_0$ 
and 2-form $\om_0$ on $\R^{2m}$. The subgroup of $\O(2m)$ preserving 
$J_0$ and $\om_0$ is $\U(m)$, so $\Hol(g)\subseteq\U(m)$. So we prove:

\begin{prop} A metric $g$ on a $2m$-manifold\/ $M$ is K\"ahler with
respect to some complex structure $J$ on $M$ if and only 
if\/~$\Hol(g)\subseteq\U(m)\subset\O(2m)$.
\label{l4prop1}
\end{prop}

\subsection{K\"ahler potentials}
\label{l43}

Let $(M,J)$ be a complex manifold. We have seen that to each
K\"ahler metric $g$ on $M$ there is associated a closed real (1,1)-form
$\om$, called the K\"ahler form. Conversely, if $\om$ is a closed real
(1,1)-form on $M$, then $\om$ is the K\"ahler form of a K\"ahler metric 
if and only if $\om$ is {\it positive}, that is, $\om(v,Jv)>0$ for 
all nonzero vectors~$v$. 

Now there is an easy way to manufacture closed real (1,1)-forms,
using the $\pd$ and $\db$ operators on $M$. If $\phi:M\ra\R$ is
smooth, then $i\pd\db\phi$ is a closed real (1,1)-form, and
every closed real (1,1)-form may be locally written in this way.
Therefore, every K\"ahler metric $g$ on $M$ may be described locally
by a function $\phi:M\ra\R$ called a {\it K\"ahler potential},
such that the K\"ahler form $\om$ satisfies~$\om=i\pd\db\phi$.

However, in general one cannot write $\om=i\pd\db\phi$ globally 
on $M$, because $i\pd\db\phi$ is {\it exact}, but $\om$ is usually
not exact (never, if $M$ is compact). Thus we are led to consider 
the {\it de Rham cohomology class} $[\om]$ of $\om$ in $H^2(M,\R)$. 
We call $[\om]$ the {\it K\"ahler class} of $g$. If two K\"ahler 
metrics $g,g'$ on $M$ lie in the same K\"ahler class, then they 
differ by a K\"ahler potential.

\begin{prop} Let\/ $(M,J)$ be a compact complex manifold, and 
let\/ $g,g'$ be K\"ahler metrics on $M$ with K\"ahler forms $\om,\om'$.
Suppose that\/ $[\om]=[\om']\in H^2(M,\R)$. Then there exists a 
smooth, real function $\phi$ on $M$ such that\/ $\om'=\om+i\pd\db\phi$.
This function $\phi$ is unique up to the addition of a constant.
\label{l4prop2}
\end{prop}

Note also that if $\om$ is the K\"ahler form of a fixed K\"ahler metric 
$g$ and $\phi$ is sufficiently small in $C^2$, then $\om'=\om+i\pd\db\phi$
is the K\"ahler form of another K\"ahler metric $g'$ on $M$, in the same 
K\"ahler class as $g$. This implies that if there exists one K\"ahler 
metric on $M$, then there exists an infinite-dimensional family --- 
K\"ahler metrics are very abundant.

\subsection{Ricci curvature and the Ricci form}
\label{l44}

Let $(M,J,g)$ be a K\"ahler manifold, with Ricci curvature $R_{ab}$.
Define the {\it Ricci form} $\rho$ by $\rho_{ac}=J_a^bR_{bc}$.
Then it turns out that $\rho_{ac}=-\rho_{ca}$, so that $\rho$ is 
a 2-form. Furthermore, it is a remarkable fact that $\rho$ is a 
{\it closed, real\/ $(1,1)$-form}. Note also that the Ricci
curvature can be recovered from $\rho$ by the 
formula~$R_{ab}=\rho_{ac}J_b^c$.

To explain this, we will give an explicit expression for the Ricci 
form. Let $(z_1,\ldots,z_m)$ be holomorphic coordinates on an open 
set $U$ in $M$. Define a smooth function $f:U\ra(0,\iy)$ by 
\e
\om^m=f\cdot\frac{(-1)^{m(m-1)/2}i^mm!}{2^m}\cdot
\d z_1\w\cdots\w\d z_m\w\d\bar z_1\w\cdots\w\d\bar z_m.
\label{ommfeq}
\e
Here the constant factor ensures that $f$ is positive, and gives 
$f\equiv 1$ when $\om$ is the standard Hermitian form on $\C^m$.
Then it can be shown that 
\e
\rho=-i\pd\db(\log f)\quad\text{on $U$,} 
\label{rhodefeq}
\e
so that $\rho$ is indeed a closed real (1,1)-form.

Using some algebraic geometry, we can interpret this. The 
{\it canonical bundle} $K_M=\La^{(m,0)}T^*M$ is a holomorphic
line bundle over $M$. The K\"ahler metric $g$ on $M$ induces a metric
on $K_M$, and the combination of metric and holomorphic structure
induces a connection $\nabla^{\sst K}$ on $K_M$. The curvature of 
this connection is a closed 2-form with values in the Lie algebra 
$\u(1)$, and identifying $\u(1)\cong\R$ we get a closed 2-form, 
which is the Ricci form.

Thus the Ricci form $\rho$ may be understood as the curvature 2-form 
of a connection $\nabla^{\sst K}$ on the canonical bundle $K_M$. So 
by characteristic class theory we may identify the de Rham cohomology 
class $[\rho]$ of $\rho$ in $H^2(M,\R)$: it satisfies
\e
[\rho]=2\pi\,c_1(K_M)=2\pi\,c_1(M), 
\label{rhoc1meq}
\e
where $c_1(M)$ is the first Chern class of $M$ in $H^2(M,\Z)$. It
is a topological invariant depending on the homotopy class of the
(almost) complex structure~$J$.

\subsection{Calabi--Yau manifolds}
\label{l45}

Here is our definition of Calabi--Yau manifold.

\begin{dfn} Let $m\ge 2$. A {\it Calabi--Yau $m$-fold} is 
a quadruple $(M,J,g,\Om)$ such that $(M,J)$ is a compact
$m$-dimensional complex manifold, $g$ a K\"ahler metric on $(M,J)$ 
with holonomy group $\Hol(g)=\SU(m)$, and $\Om$ a nonzero
constant $(m,0)$-form on $M$ called the {\it holomorphic 
volume form}, which satisfies
\e
\om^m/m!=(-1)^{m(m-1)/2}(i/2)^m\Om\w\bar\Om,
\label{omOmeq}
\e
where $\om$ is the K\"ahler form of $g$. The constant factor in 
\eq{omOmeq} is chosen to make $\Re\Om$ a {\it calibration}.
\label{l4def}
\end{dfn}

Readers are warned that there are several {\it different\/} 
definitions of Calabi--Yau manifolds in use in the literature. 
Ours is unusual in regarding $\Om$ as part of the given structure.
Some authors define a Calabi--Yau $m$-fold to be a compact K\"ahler
manifold $(M,J,g)$ with holonomy $\SU(m)$. We shall show that one
can associate a holomorphic volume form $\Om$ to such $(M,J,g)$ to
make it Calabi--Yau in our sense, and $\Om$ is unique up to phase.

\begin{lem} Let\/ $(M,J,g)$ be a compact K\"ahler manifold
with\/ $\Hol(g)=\SU(m)$. Then $M$ admits a holomorphic
volume form $\Om$, unique up to change of phase $\Om\mapsto
{\rm e}^{i\th}\Om$, such that\/ $(M,J,g,\Om)$ is a Calabi--Yau
manifold.
\label{l4lem}
\end{lem}

\begin{proof} Let $(M,J,g)$ be compact and K\"ahler with 
$\Hol(g)=\SU(m)$. Now the holonomy group $\SU(m)$ preserves 
the standard metric $g_0$ and K\"ahler form $\om_0$ on 
$\C^m$, and an $(m,0)$-form $\Om_0$ given by
\begin{align*}
g_0=\ms{\d z_1}+\cdots+\ms{\d z_m},\quad
\om_0&=\frac{i}{2}(\d z_1\w\d\bar z_1+\cdots+\d z_m\w\d\bar z_m),\\
\text{and}\quad\Om_0&=\d z_1\w\cdots\w\d z_m.
\end{align*}

Thus, by Theorem \ref{l2thm1} there exist corresponding constant 
tensors $g$, $\om$ (the K\"ahler form), and $\Om$ on $(M,J,g)$. 
Since $\om_0$ and $\Om_0$ satisfy
\begin{equation*}
\om_0^m/m!=(-1)^{m(m-1)/2}(i/2)^m\Om_0\w\bar\Om_0
\end{equation*}
on $\C^m$, it follows that $\om$ and $\Om$ satisfy \eq{omOmeq}
at each point, so $(M,J,g,\Om)$ is Calabi--Yau. It is easy to 
see that $\Om$ is unique up to change of phase.
\end{proof}

Suppose $(M,J,g,\Om)$ is a Calabi--Yau $m$-fold. Then $\Om$
is a constant section of the canonical bundle $K_M$. As $\Om$
is constant, it is holomorphic. Thus the canonical bundle 
$K_M$ admits a nonvanishing holomorphic section, so $(M,J)$ has 
{\it trivial canonical bundle}, which implies that~$c_1(M)=0$.

Further, the connection $\nabla^{\sst K}$ on $K_M$ must be {\it flat}. 
However, from \S\ref{l44} the curvature of $\nabla^{\sst K}$ is the 
Ricci form $\rho$. Therefore $\rho\equiv 0$, and $g$ is Ricci-flat. 
That is, Calabi--Yau $m$-folds are automatically {\it Ricci-flat}.
More generally, the following proposition explains the relationship
between the Ricci curvature and holonomy group of a K\"ahler metric.

\begin{prop} Let\/ $(M,J,g)$ be a K\"ahler $m$-fold with\/ $\Hol(g)
\subseteq\SU(m)$. Then $g$ is Ricci-flat. Conversely, let\/ $(M,J,g)$ 
be a Ricci-flat K\"ahler $m$-fold. If\/ $M$ is simply-connected or 
$K_M$ is trivial, then~$\Hol(g)\subseteq\SU(m)$.
\label{l4prop3}
\end{prop}

In the last part, $M$ simply-connected implies that $K_M$ is
trivial for Ricci-flat K\"ahler manifolds, but not vice versa.

\subsection{Exercises}
\label{l47}

\begin{question}Let $U$ be a simply-connected subset of $\C^m$ with 
coordinates $(z_1,\ldots,z_m)$, and $g$ a Ricci-flat K\"ahler metric on $U$ 
with K\"ahler form $\om$. Use equations \eq{ommfeq} and \eq{rhodefeq} to 
show that there exists a holomorphic $(m,0)$-form $\Om$ on $U$ satisfying
\begin{equation*}
\om^m/m!=(-1)^{m(m-1)/2}(i/2)^m\Om\w\bar\Om.
\end{equation*}

{\bf Hint:} Write $\Om=F\,\d z_1\w\cdots\w\d z_m$ for some holomorphic
function $F$. Use the fact that if $f$ is a real function on a 
simply-connected subset $U$ of $\C^m$ and $\pd\db f\equiv 0$, then 
$f$ is the real part of a holomorphic function on~$U$.
\end{question}

\begin{question}Let $\C^2$ have complex coordinates $(z_1,z_2)$, 
and define $u=\ms{z_1}+\ms{z_2}$. Let $f:[0,\iy)\ra\R$ be a smooth function,
and define a closed real (1,1)-form $\om$ on $\C^2$ by~$\om=i\pd\db f(u)$.
\anext Calculate the conditions on $f$ for $\om$ to be the
K\"ahler form of a K\"ahler metric $g$ on $\C^2$.

(You can define $g$ by $g(v,w)=\om(v,Jw)$, and need to ensure that
$g$ is positive definite).

\anext Supposing $g$ is a metric, calculate the conditions on
$f$ for $g$ to be Ricci-flat. You should get an o.d.e.\ on $f$.
If you can, solve this o.d.e., and write down the corresponding
K\"ahler metrics in coordinates.
\end{question}

\section{The Calabi Conjecture and \\ 
constructions of Calabi--Yau $m$-folds}
\label{l5}

The {\it Calabi Conjecture} specifies which closed (1,1)-forms on
a compact complex manifold can be the Ricci form of a K\"ahler metric.
It was posed by Calabi in 1954, and proved by Yau in 1976. We
shall explain the conjecture, and sketch its proof. An important 
application of the Calabi Conjecture is the construction of 
large numbers of {\it Calabi--Yau manifolds}. We explain some 
ways to do this, using algebraic geometry. 

A good general reference for this section is my book 
\cite[\S\S 5, 6 \& 7.3]{Joyc2}, which includes a proof of
the Calabi Conjecture. Other references on the Calabi Conjecture
are Aubin's book \cite{Aubi} and Yau \cite{Yau}, which is the
original proof of the conjecture, but fairly hard going unless
you know a lot of analysis.

\subsection{The Calabi Conjecture}
\label{l51}

Let $(M,J)$ be a compact, complex manifold, and $g$ a K\"ahler metric 
on $M$, with Ricci form $\rho$. From \S\ref{l44}, $\rho$ is a closed 
real (1,1)-form and $[\rho]=2\pi\,c_1(M)\in H^2(M,\R)$. The Calabi
Conjecture specifies which closed (1,1)-forms can be the Ricci forms 
of a K\"ahler metric on~$M$. 
\medskip

\noindent{\bf The Calabi Conjecture} {\it Let\/ $(M,J)$ be a compact, 
complex manifold, and\/ $g$ a K\"ahler metric on $M$, with K\"ahler form 
$\om$. Suppose that\/ $\rho'$ is a real, closed\/ $(1,1)$-form on $M$ with\/ 
$[\rho']=2\pi\,c_1(M)$. Then there exists a unique K\"ahler metric $g'$ on 
$M$ with K\"ahler form $\om'$, such that\/ $[\om']=[\om]\in H^2(M,\R)$, and 
the Ricci form of\/ $g'$ is~$\rho'$.}
\medskip

Note that $[\om']=[\om]$ says that $g$ and $g'$ are in the {\it same 
K\"ahler class}. The conjecture was posed by Calabi in 1954, and was 
eventually proved by Yau in 1976. Its importance to us is that when 
$c_1(M)=0$ we can take $\rho'\equiv 0$, and then $g'$ is Ricci-flat. 
Thus, assuming the Calabi Conjecture we prove:

\begin{cor} Let\/ $(M,J)$ be a compact complex manifold with\/ 
$c_1(M)=0$ in $H^2(M,\R)$. Then every K\"ahler class on $M$ contains
a unique Ricci-flat K\"ahler metric~$g$.
\label{l5cor1}
\end{cor}

If in addition $M$ is simply-connected or $K_M$ is trivial, then 
Proposition \ref{l4prop3} shows that $\Hol(g)\subseteq\SU(m)$.
When $\Hol(g)=\SU(m)$, which will happen under certain fairly
simple topological conditions on $M$, then by Lemma \ref{l4lem} we 
can make $(M,J,g)$ into a {\it Calabi--Yau manifold} $(M,J,g,\Om)$. 
So Yau's proof of the Calabi Conjecture gives a way to find examples 
of Calabi--Yau manifolds, which is how Calabi--Yau manifolds got 
their name.

Note that we know {\it almost nothing} about the Ricci-flat K\"ahler 
metric $g$ except that it exists; we cannot write it down explicitly in 
coordinates, for instance. In fact, {\it no} explicit non-flat examples 
of Calabi--Yau metrics on compact manifolds are known at all.

\subsection{Sketch of the proof of the Calabi Conjecture}
\label{l52}

The Calabi Conjecture is proved by rewriting it as a second-order
nonlinear elliptic p.d.e.\ upon a real function $\phi$ on $M$,
and then showing that this p.d.e.\ has a unique solution. We first
explain how to rewrite the Calabi Conjecture as a~p.d.e.

Let $(M,J)$ be a compact, complex manifold, and let $g,g'$ be two
K\"ahler metrics on $M$ with K\"ahler forms $\om,\om'$ and Ricci forms
$\rho,\rho'$. Suppose $g,g'$ are in the same K\"ahler class, so 
that $[\om']=[\om]\in H^2(M,\R)$. Define a smooth function 
$f:M\ra\R$ by $(\om')^m={\rm e}^f\om^m$. Then from equations 
\eq{ommfeq} and \eq{rhodefeq} of \S\ref{l44}, we find that 
$\rho'=\rho-i\pd\db f$. Furthermore, as $[\om']=[\om]$ in 
$H^2(M,\R)$, we have $[\om']^m=[\om]^m$ in $H^{2m}(M,\R)$, 
and thus~$\int_M{\rm e}^f\om^m=\int_M\om^m$.

Now suppose that we are given the real, closed $(1,1)$-form 
$\rho'$ with $[\rho']=2\pi\,c_1(M)$, and want to construct
a metric $g'$ with $\rho'$ as its Ricci form. Since 
$[\rho]=[\rho']=2\pi\,c_1(M)$, $\rho-\rho'$ is an {\it exact}\/ 
real (1,1)-form, and so by the $\pd\db$-Lemma there exists a 
smooth function $f:M\ra\R$ with $\rho-\rho'=i\pd\db f$. This
$f$ is unique up to addition of a constant, but the constant
is fixed by requiring that $\int_M{\rm e}^f\om^m=\int_M\om^m$.
Thus we have proved:

\begin{prop} Let\/ $(M,J)$ be a compact complex manifold, $g$ a 
K\"ahler metric on $M$ with K\"ahler form $\om$ and Ricci form 
$\rho$, and $\rho'$ a real, closed\/ $(1,1)$-form on $M$ 
with\/ $[\rho']=2\pi\,c_1(M)$. Then there is a unique 
smooth function $f:M\ra\R$ such that
\e
\rho'=\rho-i\pd\db f
\quad\text{and}\quad
\int_M{\rm e}^f\om^m=\int_M\om^m,
\e
and a K\"ahler metric $g$ on $M$ with K\"ahler form $\om'$ satisfying
$[\om']=[\om]$ in $H^2(M,\R)$ has Ricci form $\rho'$
if and only if\/~$(\om')^m={\rm e}^f\om^m$. 
\label{l5prop1}
\end{prop}

Thus we have transformed the Calabi Conjecture from seeking a 
metric $g'$ with {\it prescribed Ricci curvature} $\rho'$ to 
seeking a metric $g'$ with {\it prescribed volume form} $(\om')^m$.
This is an important simplification, because the Ricci curvature
depends on the second derivatives of $g'$, but the volume form
depends only on $g'$ and not on its derivatives.

Now by Proposition \ref{l4prop2}, as $[\om']=[\om]$ we may write
$\om'=\om+i\pd\db\phi$ for $\phi$ a smooth real function on $M$,
unique up to addition of a constant. We can fix the constant by
requiring that $\int_M\phi\,\d V_g=0$. So, from Proposition
\ref{l5prop1} we deduce that the Calabi Conjecture is equivalent to: 
\medskip

\noindent{\bf The Calabi Conjecture (second version)} {\it Let\/ 
$(M,J)$ be a compact, complex manifold, and\/ $g$ a K\"ahler metric on 
$M$, with K\"ahler form $\om$. Let\/ $f$ be a smooth real function on $M$
satisfying $\int_M{\rm e}^f\om^m=\int_M\om^m$. Then there exists a 
unique smooth real function $\phi$ such that
\begin{itemize}
\item[{\rm(i)}] $\om+i\pd\db\phi$ is a positive $(1,1)$-form,
that is, it is the K\"ahler form of some K\"ahler metric~$g'$,
\item[{\rm(ii)}] $\int_M\phi\,\d V_g=0$, and 
\item[{\rm(iii)}] $(\om+i\pd\db\phi)^m={\rm e}^f\om^m$ on~$M$.
\end{itemize}}
\medskip

This reduces the Calabi Conjecture to a problem in analysis, that of 
showing that the nonlinear p.d.e.\ $(\om+i\pd\db\phi)^m={\rm e}^f\om^m$ 
has a solution $\phi$ for every suitable function $f$. To prove this 
second version of the Calabi Conjecture, Yau used the {\it continuity 
method}.

For each $t\in[0,1]$, define $f_t=tf+c_t$, where $c_t$ is the unique real 
constant such that ${\rm e}^{c_t}\int_M{\rm e}^{tf}\om^m=\int_M\om^m$.
Then $f_t$ depends smoothly on $t$, with $f_0\equiv 0$ and $f_1\equiv f$.
Define $S$ to be the set of $t\in[0,1]$ such that there exists a smooth
real function $\phi$ on $M$ satisfying parts (i) and (ii) above, and also
\begin{itemize}
\item[(iii$)'$] $(\om+i\pd\db\phi)^m={\rm e}^{f_t}\om^m$ on $M$.
\end{itemize}

The idea of the continuity method is to show that $S$ is both {\it open} 
and {\it closed}\/ in $[0,1]$. Thus, $S$ is a connected subset of $[0,1]$,
so $S=\emptyset$ or $S=[0,1]$. But $0\in S$, since as $f_0\equiv 0$ parts 
(i), (ii) and (iii$)'$ are satisfied by $\phi\equiv 0$. Thus $S=[0,1]$.
In particular, (i), (ii) and (iii$)'$ admit a solution $\phi$ when $t=1$.
As $f_1\equiv f$, this $\phi$ satisfies (iii), and the Calabi Conjecture 
is proved.

Showing that $S$ is open is fairly easy, and was done by Calabi. It
depends on the fact that (iii) is an {\it elliptic} p.d.e.\ --- basically,
the operator $\phi\mapsto(\om+i\pd\db\phi)^m$ is rather like a nonlinear
Laplacian --- and uses only standard facts about elliptic operators.

However, showing that $S$ is closed is much more difficult. One must 
prove that $S$ contains its limit points. That is, if $(t_n)_{n=1}^\iy$
is a sequence in $S$ converging to $t\in[0,1]$ then there exists
a sequence $(\phi_n)_{n=1}^\iy$ satisfying (i), (ii) and 
$(\om+i\pd\db\phi_n)^m={\rm e}^{f_{t_n}}\om^m$ for $n=1,2,\ldots$,
and we need to show that $\phi_n\ra\phi$ as $n\ra\iy$ for some
smooth real function $\phi$ satisfying (i), (ii) and (iii$)'$,
so that~$t\in S$.

The thing you have to worry about is that the sequence $(\phi_n)_{n=1}^\iy$
might converge to some horrible non-smooth function, or might not converge 
at all. To prove this doesn't happen you need {\it a priori estimates} on
the $\phi_n$ and all their derivatives. In effect, you need upper bounds 
on $\md{\nabla^k\phi_n}$ for all $n$ and $k$, bounds which are allowed to 
depend on $M,J,g,k$ and $f_{t_n}$, but not on $n$ or $\phi_n$. These a 
priori estimates were difficult to find, because the nonlinearities in 
$\phi$ of $(\om+i\pd\db\phi)^m={\rm e}^f\om^m$ are of a particularly 
nasty kind, and this is why it took so long to prove the Calabi Conjecture.

\subsection{Calabi--Yau 2-folds and $K3$ surfaces}
\label{l53}

Recall from \S\ref{l35} that the {\it K\"ahler holonomy groups} are
$\U(m)$, $\SU(m)$ and $\Sp(k)$. Calabi--Yau manifolds of complex
dimension $m$ have holonomy $\SU(m)$ for $m\ge 2$, and hyperk\"ahler
manifolds of complex dimension $2k$ have holonomy $\Sp(k)$ for
$k\ge 1$. In complex dimension 2 these coincide, as $\SU(2)=\Sp(1)$.
Because of this, Calabi--Yau 2-folds have special features which
are not present in Calabi--Yau $m$-folds for~$m\ge 3$.

Calabi--Yau 2-folds are very well understood, through the
classification of compact complex surfaces. A $K3$ {\it surface} 
is defined to be a compact, complex surface $(X,J)$ with 
$h^{1,0}(X)=0$ and trivial canonical bundle. All Calabi--Yau 
2-folds are $K3$ surfaces, and conversely, every $K3$ surface
$(X,J)$ admits a family of K\"ahler metrics $g$ making it into a
Calabi--Yau 2-fold. All $K3$ surfaces $(X,J)$ are diffeomorphic, 
sharing the same smooth 4-manifold $X$, which is simply-connected, 
with Betti numbers $b^2=22$, $b^2_+=3$, and~$b^2_-=19$. 

The moduli space ${\cal M}_{K3}$ of $K3$ surfaces is a connected 
20-dimensional singular complex manifold, which can be described 
very precisely via the `Torelli Theorems'. Some $K3$ surfaces are 
{\it algebraic}, that is, they can be embedded as complex submanifolds 
in $\CP^N$ for some $N$, and some are not. The set of algebraic $K3$ 
surfaces is a countable, dense union of 19-dimensional subvarieties
in ${\cal M}_{K3}$. Each $K3$ surface $(X,J)$ admits a real 
20-dimensional family of Calabi--Yau metrics $g$, so the family of 
Calabi--Yau 2-folds $(X,J,g)$ is a nonsingular 60-dimensional real 
manifold.

\subsection{General properties of Calabi--Yau $m$-folds for $m\ge 3$}
\label{l54}

Using general facts about Ricci-flat manifolds (the {\it Cheeger--Gromoll
Theorem}) one can show that every Calabi--Yau $m$-fold $(M,J,g,\Om)$ has 
finite fundamental group. Also, using the `Bochner argument' one
can show that any closed $(p,0)$-form $\xi$ on $M$ is constant under
the Levi-Civita connection $\nabla$ of~$g$.

However, the set of constant tensors on $M$ is determined by the
holonomy group $\Hol(g)$ of $g$, which is $\SU(m)$ by definition.
It is easy to show that the vector space of closed $(p,0)$-forms
on $M$ is $\C$ if $p=0,m$ and 0 otherwise. But the vector space of 
closed $(p,0)$ forms is the {\it Dolbeault cohomology group} 
$H^{p,0}(M)$, whose dimension is the {\it Hodge number} $h^{p,0}$ 
of $M$. Thus we prove:

\begin{prop} Let\/ $(M,J,g,\Om)$ be a Calabi--Yau $m$-fold with 
Hodge numbers $h^{p,q}$. Then $M$ has finite fundamental group, 
$h^{0,0}=h^{m,0}=1$ and\/ $h^{p,0}=0$ for~$p\ne 0,m$.
\label{l5prop2}
\end{prop}

For $m\ge 3$ this gives $h^{2,0}(M)=0$, and this has important 
consequences for the complex manifold $(M,J)$. It can be shown 
that a complex line bundle $L$ over a compact K\"ahler manifold $(M,J,g)$ 
admits a holomorphic structure if and only if $c_1(L)$ lies in 
$H^{1,1}(M)\subseteq H^2(M,\C)$. But $H^2(M,\C)=H^{2,0}(M)\op 
H^{1,1}(M)\op H^{0,2}(M)$, and $H^{2,0}(M)=H^{0,2}(M)=0$ as 
$h^{2,0}(M)=0$. Thus $H^{1,1}(M)=H^2(M,\C)$, and so {\it every}
complex line bundle $L$ over $M$ admits a holomorphic structure.

Thus, Calabi--Yau $m$-folds for $m\ge 3$ are richly endowed with
holomorphic line bundles. Using the {\it Kodaira Embedding Theorem}
one can show that some of these holomorphic line bundles admit
many holomorphic sections. By taking a line bundle with enough 
holomorphic sections (a {\it very ample} line bundle) we can
construct an embedding of $M$ in $\CP^N$ as a complex submanifold.
So we prove:

\begin{thm} Let\/ $(M,J,g,\Om)$ be a Calabi--Yau manifold of 
dimension $m\ge 3$. Then $M$ is projective. That is, $(M,J)$ is 
isomorphic as a complex manifold to a complex submanifold of\/ 
$\CP^N$, and is an algebraic variety.
\label{l5thm1}
\end{thm}

This shows that Calabi--Yau manifolds (or at least, the complex
manifolds underlying them) can be studied using {\it complex 
algebraic geometry}. 

\subsection{Constructions of Calabi--Yau $m$-folds}
\label{l55}

The easiest way to find examples of Calabi--Yau $m$-folds for 
$m\ge 3$ is to choose a method of generating a large number of 
complex algebraic varieties, and then check the topological 
conditions to see which of them are Calabi--Yau. Here are
some ways of doing this.

\begin{itemize} 
\item {\bf Hypersurfaces in $\CP^{m+1}$}. Suppose that
$X$ is a smooth degree $d$ hypersurface in $\CP^{m+1}$. When is
$X$ a Calabi--Yau manifold? Well, using the {\it adjunction formula}
one can show that the canonical bundle of $X$ is given by
$K_X=L^{d-m-2}\vert_X$, where $L\ra\CP^{m+1}$ is the hyperplane
line bundle on $\CP^{m+1}$. 

Therefore $K_X$ is trivial if and only if $d=m+2$. It is not
difficult to show that any smooth hypersurface of degree $m+2$
in $\CP^{m+1}$ is a Calabi--Yau $m$-fold. All such hypersurfaces
are diffeomorphic, for fixed $m$. For instance, the {\it Fermat 
quintic}
\begin{equation*}
\bigl\{[z_0,\ldots,z_4]\in\CP^4:z_0^5+\cdots+z_4^5=0\bigr\}
\end{equation*}
is a Calabi--Yau 3-fold, with Betti numbers $b^0=1$, $b^1=0$,
$b^2=1$ and~$b^3=204$.

\item {\bf Complete intersections in $\CP^{m+k}$}. In the
same way, suppose $X$ is a {\it complete intersection} of
transverse hypersurfaces $H_1,\ldots,H_k$ in $\CP^{m+k}$ of 
degrees $d_1,\ldots,d_k$, with each $d_j\ge 2$. It can be 
shown that $X$ is Calabi--Yau $m$-fold if and only if 
$d_1+\cdots+d_k=m+k+1$. This yields a finite number of 
topological types in each dimension~$m$.

\item {\bf Hypersurfaces in toric varieties}. A {\it toric variety} 
is a complex $m$-manifold $X$ with a holomorphic action of $(\C^*)^m$ 
which is transitive and free upon a dense open set in $X$. Toric
varieties can be constructed and studied using only a finite amount
of combinatorial data.

The conditions for a smooth hypersurface in a compact toric variety
to be a Calabi--Yau $m$-fold can be calculated using this combinatorial
data. Using a computer, one can generate a large (but finite) number
of Calabi--Yau $m$-folds, at least when $m=3$, and calculate their
topological invariants such as Hodge numbers. This has been
done by Candelas, and other authors.

\item {\bf Resolution of singularities}. Suppose you have some
way of producing examples of {\it singular} Calabi--Yau $m$-folds
$Y$. Often it is possible to find a {\it resolution} $X$ of $Y$ 
with holomorphic map $\pi:X\ra Y$, such that $X$ is a nonsingular
Calabi--Yau $m$-fold. Basically, each singular point in $Y$ is
replaced by a finite union of complex submanifolds in~$X$. 

Resolutions which preserve the Calabi--Yau property are called
{\it crepant resolutions}, and are well understood when $m=3$.
For certain classes of singularities, such as singularities
of Calabi--Yau 3-orbifolds, a crepant resolution always exists.

This technique can be applied in a number of ways. For instance,
you can start with a nonsingular Calabi--Yau $m$-fold $X$,
deform it till you get a singular Calabi--Yau $m$-fold $Y$,
and then resolve the singularities of $Y$ to get a second 
nonsingular Calabi--Yau $m$-fold $X'$ with different topology to~$X$.

Another method is to start with a nonsingular Calabi--Yau $m$-fold
$X$, divide by the action of a finite group $G$ preserving the
Calabi--Yau structure to get a singular Calabi--Yau manifold
(orbifold) $Y=X/G$, and then resolve the singularities of $Y$ 
to get a second nonsingular Calabi--Yau $m$-fold $X'$ with 
different topology to~$X$.
\end{itemize}

\subsection{Exercises}
\label{l56}

\begin{question}The most well-known examples of Calabi--Yau 
3-folds are {\it quintics} $X$ in $\CP^4$, defined by 
\begin{equation*}
X=\bigl\{[z_0,\ldots,z_4]\in\CP^4:p(z_0,\ldots,z_4)=0\bigr\},
\end{equation*}
where $p(z_0,\ldots,z_4)$ is a homogeneous quintic polynomial
in its arguments. Every nonsingular quintic has Hodge numbers
$h^{1,1}=h^{2,2}=1$ and~$h^{2,1}=h^{1,2}=101$.
\label{l5q1}

\inext Calculate the dimension of the vector space of homogeneous
quintic polynomials $p(z_0,\ldots,z_4)$. Hence find the dimension
of the moduli space of nonsingular quintics in $\CP^4$. (A generic
quintic is nonsingular).

\inext Identify the group of complex automorphisms of $\CP^4$ and
calculate its dimension.

\inext Hence calculate the dimension of the moduli space of
quintics in $\CP^4$ up to automorphisms of $\CP^4$.
\end{question}

It is a general fact that if $(X,J,g)$ is a Calabi--Yau
3-fold, then the moduli space of complex deformations of $(X,J)$
has dimension $h^{2,1}(X)$, and each nearby deformation is a
Calabi--Yau 3-fold. In this case, $h^{2,1}(X)=101$, and this
should be your answer to (iii). That is, deformations of quintics 
in $\CP^4$ are also quintics in~$\CP^4$.

\begin{question}One can also construct Calabi--Yau 3-folds as the 
complete intersection of two cubics in $\CP^5$,
\begin{equation*}
X=\bigl\{[z_0,\ldots,z_5]\in\CP^5:p(z_0,\ldots,z_5)=
q(z_0,\ldots,z_5)=0\bigr\},
\end{equation*}
where $p,q$ are linearly independent homogeneous cubic polynomials.
Using the method of Question \ref{l5}.\ref{l5q1}, calculate the 
dimension of the moduli space of such complete intersections up to 
automorphisms of $\CP^5$, and hence predict~$h^{2,1}(X)$.
\end{question}

\section{Introduction to calibrated geometry}
\label{l6}

The theory of {\it calibrated geometry} was invented by Harvey and 
Lawson \cite{HaLa}. It concerns {\it calibrated submanifolds}, a 
special kind of {\it minimal submanifold} of a Riemannian manifold 
$M$, which are defined using a closed form on $M$ called a 
{\it calibration}. It is closely connected with the theory of 
Riemannian holonomy groups because Riemannian manifolds with special 
holonomy usually come equipped with one or more natural calibrations. 

Some references for this section are Harvey and Lawson 
\cite[\S I, \S II]{HaLa}, Harvey \cite{Harv} and the author
\cite[\S 3.7]{Joyc2}. Some background reading on minimal
submanifolds and Geometric Measure Theory is Lawson \cite{Laws}
and Morgan~\cite{Morg}.

\subsection{Minimal submanifolds}
\label{l61}

Let $(M,g)$ be an $n$-dimensional Riemannian manifold, and $N$ a 
compact $k$-dimensional submanifold of $M$. Regard $N$ as an immersed
submanifold $(N,\iota)$, with immersion $\iota:N\ra M$. Using 
the metric $g$ we can define the {\it volume} $\Vol(N)$ of $N$, 
by integration over $N$. We call $N$ a {\it minimal submanifold}
if its volume is stationary under small variations of the
immersion $\iota:N\ra M$. When $k=1$, a curve in $M$ is minimal 
if and only if it is a {\it geodesic}.

Let $\nu\ra N$ be the normal bundle of $N$ in $M$, so that 
$TM\vert_N=TN\op\nu$ is an orthogonal direct sum. The {\it second 
fundamental form} is a section $B$ of $S^2T^*N\ot\nu$ such that 
whenever $v,w$ are vector fields on $M$ with $v\vert_N,w\vert_N$ 
sections of $TN$ over $N$, then $B\cdot \bigl(v\vert_N\ot w\vert_N
\bigr)=\pi_\nu\bigl(\nabla_vw\vert_N\bigr)$, where `$\cdot$' contracts 
$S^2T^*N$ with $TN\ot TN$, $\nabla$ is the Levi-Civita connection
of $g$, and $\pi_\nu$ is the projection to $\nu$ in the 
splitting~$TM\vert_N=TN\op\nu$.

The {\it mean curvature vector} $\ka$ of $N$ is the trace of the 
second fundamental form $B$ taken using the metric $g$ on $N$. It 
is a section of the normal bundle $\nu$. It can be shown by the
Euler--Lagrange method that a submanifold $N$ is minimal if and
only if its mean curvature vector $\ka$ is zero. Note that this 
is a local condition. Therefore we can also define noncompact 
submanifolds $N$ in $M$ to be minimal if they have zero mean 
curvature. This makes sense even when $N$ has infinite volume.

If $\iota:N\ra M$ is a immersed submanifold, then the mean 
curvature $\ka$ of $N$ depends on $\iota$ and its first and
second derivatives, so the condition that $N$ be minimal
is a {\it second-order} equation on $\iota$. Note that minimal 
submanifolds may not have minimal area, even amongst nearby 
homologous submanifolds. For instance, the equator in ${\cal S}^2$ 
is minimal, but does not minimize length amongst lines of latitude.

The following argument is important in the study of minimal
submanifolds. Let $(M,g)$ be a compact Riemannian manifold, and $\al$
a nonzero homology class in $H_k(M,\Z)$. We would like to find a
compact, minimal immersed, $k$-dimensional submanifold $N$ in $M$ 
with homology class $[N]=\al$. To do this, we choose a minimizing 
sequence $(N_i)_{i=1}^\iy$ of compact submanifolds $N_i$ with 
$[N_i]=\al$, such that $\Vol(N_i)$ approaches the infimum of 
volumes of submanifolds with homology class $\al$ as~$i\ra\iy$. 

Pretend for the moment that the set of all closed $k$-dimensional 
submanifolds $N$ with $\Vol(N)\le C$ is a {\it compact} topological 
space. Then there exists a subsequence $(N_{i_j})_{j=1}^\iy$ which 
converges to some submanifold $N$, which is the minimal submanifold 
we want. In fact this does not work, because the set of submanifolds 
$N$ does not have the compactness properties we need. 

However, if we work instead with {\it rectifiable currents}, which are a 
measure-theoretic generalization of submanifolds, one can show that every 
integral homology class $\al$ in $H_k(M,\Z)$ is represented by a minimal 
rectifiable current. One should think of rectifiable currents as a class 
of singular submanifolds, obtained by completing the set of nonsingular 
submanifolds with respect to some norm. They are studied in the subject 
of {\it Geometric Measure Theory}.

The question remains: how close are these minimal rectifiable currents
to being submanifolds? For example, it is known that a $k$-dimensional 
minimal rectifiable current in a Riemannian $n$-manifold is an embedded 
submanifold except on a singular set of Hausdorff dimension at most 
$k-2$. When $k=2$ or $k=n-1$ one can go further. In general, it is 
important to understand the possible singularities of such singular 
minimal submanifolds.

\subsection{Calibrations and calibrated submanifolds}
\label{l62}

Let $(M,g)$ be a Riemannian manifold. An {\it oriented
tangent $k$-plane} $V$ on $M$ is a vector subspace $V$ of
some tangent space $T_xM$ to $M$ with $\dim V=k$, equipped
with an orientation. If $V$ is an oriented tangent $k$-plane
on $M$ then $g\vert_V$ is a Euclidean metric on $V$, so 
combining $g\vert_V$ with the orientation on $V$ gives a 
natural {\it volume form} $\vol_V$ on $V$, which is a 
$k$-form on~$V$.

Now let $\vp$ be a closed $k$-form on $M$. We say that
$\vp$ is a {\it calibration} on $M$ if for every oriented
$k$-plane $V$ on $M$ we have $\vp\vert_V\le \vol_V$. Here
$\vp\vert_V=\al\cdot\vol_V$ for some $\al\in\R$, and 
$\vp\vert_V\le\vol_V$ if $\al\le 1$. Let $N$ be an 
oriented submanifold of $M$ with dimension $k$. Then 
each tangent space $T_xN$ for $x\in N$ is an oriented
tangent $k$-plane. We say that $N$ is a {\it calibrated 
submanifold} or $\vp$-{\it submanifold} if 
$\vp\vert_{T_xN}=\vol_{T_xN}$ for all~$x\in N$.

All calibrated submanifolds are automatically {\it minimal 
submanifolds}. We prove this in the compact case, but it is
true for noncompact submanifolds as well.

\begin{prop} Let\/ $(M,g)$ be a Riemannian manifold, $\vp$ a 
calibration on $M$, and\/ $N$ a compact $\vp$-submanifold 
in $M$. Then $N$ is volume-minimizing in its homology class.
\label{l6prop1}
\end{prop}

\begin{proof} Let $\dim N=k$, and let $[N]\in H_k(M,\R)$ 
and $[\vp]\in H^k(M,\R)$ be the homology and cohomology 
classes of $N$ and $\vp$. Then
\begin{equation*}
[\vp]\cdot[N]=\int_{x\in N}\vp\big\vert_{T_xN}=
\int_{x\in N}{\ts\vol_{T_xN}}=\Vol(N),
\end{equation*}
since $\vp\vert_{T_xN}=\vol_{T_xN}$ for each $x\in N$, as
$N$ is a calibrated submanifold. If $N'$ is any other compact 
$k$-submanifold of $M$ with $[N']=[N]$ in $H_k(M,\R)$, then
\begin{equation*}
[\vp]\cdot[N]=[\vp]\cdot[N']=\int_{x\in N'}\vp\big\vert_{T_xN'}
\le\int_{x\in N'}{\ts\vol_{T_xN'}}=\Vol(N'),
\end{equation*}
since $\vp\vert_{T_xN'}\le\vol_{T_xN'}$ because $\vp$ is a 
calibration. The last two equations give $\Vol(N)\le\Vol(N')$. 
Thus $N$ is volume-minimizing in its homology class.
\end{proof}

Now let $(M,g)$ be a Riemannian manifold with a calibration $\vp$,
and let $\iota:N\ra M$ be an immersed submanifold. Whether
$N$ is a $\vp$-submanifold depends upon the tangent spaces of 
$N$. That is, it depends on $\iota$ and its first derivative. 
So, to be calibrated with respect to $\vp$ is a {\it first-order}
equation on $\iota$. But if $N$ is calibrated then $N$ is minimal, 
and we saw in \S\ref{l61} that to be minimal is a {\it second-order} 
equation on~$\iota$. 

One moral is that the calibrated equations, being first-order,
are often easier to solve than the minimal submanifold equations,
which are second-order. So calibrated geometry is a fertile source
of examples of minimal submanifolds. 

\subsection{Calibrated submanifolds of $\R^n$}
\label{l63}

One simple class of calibrations is to take $(M,g)$ to be
$\R^n$ with the Euclidean metric, and $\vp$ to be a constant
$k$-form on $\R^n$, such that $\vp\vert_V\le \vol_V$ for
every oriented $k$-dimensional vector subspace $V\subseteq\R^n$.
Each such $\vp$ defines a class of minimal $k$-submanifolds
in $\R^n$. However, this class may be very small, or even empty.
For instance, $\vp=0$ is a calibration on $\R^n$, but has no 
calibrated submanifolds. 

For each constant calibration $k$-form $\vp$ on $\R^n$, define
${\cal F}_\vp$ to be the set of oriented $k$-dimensional vector 
subspaces $V$ of $\R^n$ such that $\ts\vp\vert_V=\vol_V$. Then 
an oriented submanifold $N$ of $\R^n$ is a $\vp$-submanifold 
if and only if each tangent space $T_xN$ lies in ${\cal F}_\vp$.
To be interesting, a calibration $\vp$ should define a fairly
abundant class of calibrated submanifolds, and this will only 
happen if ${\cal F}_\vp$ is reasonably large.

Define a {\it partial order} $\preceq$ on the set of constant 
calibration $k$-forms $\vp$ on $\R^n$ by $\vp\preceq\vp'$ if
${\cal F}_\vp\subseteq{\cal F}_{\vp'}$. A calibration $\vp$ is
{\it maximal} if it is maximal with respect to this partial 
order. A maximal calibration $\vp$ is one in which 
${\cal F}_\vp$ is as large as possible.

It is an interesting problem to determine the maximal calibrations
$\vp$ on $\R^n$. The symmetry group $G\subset\O(n)$ of a maximal
calibration is usually quite large. This is because if $V\in{\cal F}_\vp$
and $\ga\in G$ then $\ga\cdot V\in{\cal F}_\vp$, that is, $G$ acts 
on ${\cal F}_\vp$. So if $G$ is big we expect ${\cal F}_\vp$ to be 
big too. Symmetry groups of maximal calibrations are often possible 
holonomy groups of Riemannian metrics, and the classification problem 
for maximal calibrations can be seen as in some ways parallel to the 
classification problem for Riemannian holonomy groups.

\subsection{Calibrated submanifolds and special holonomy}
\label{l64}

Next we explain the connection with Riemannian holonomy. Let 
$G\subset{\rm O}(n)$ be a possible holonomy group of a Riemannian 
metric. In particular, we can take $G$ to be one of the
holonomy groups $\U(m)$, $\SU(m)$, $\Sp(m)$, $G_2$ or 
Spin(7) from Berger's classification. Then $G$ acts on 
the $k$-forms $\La^k(\R^n)^*$ on $\R^n$, so we can look
for $G$-invariant $k$-forms on~$\R^n$.

Suppose $\vp_0$ is a nonzero, $G$-invariant $k$-form on 
$\R^n$. By rescaling $\vp_0$ we can arrange that for each 
oriented $k$-plane $U\subset\R^n$ we have $\vp_0\vert_U\le\vol_U$, 
and that $\vp_0\vert_U=\vol_U$ for at least one such $U$. Thus
${\cal F}_{\vp_0}$ is nonempty. Since $\vp_0$ is $G$-invariant, if 
$U\in{\cal F}_{\vp_0}$ then $\ga\cdot U\in{\cal F}_{\vp_0}$ for all 
$\ga\in G$. Generally this means that ${\cal F}_{\vp_0}$ is
`reasonably large'.

Let $M$ be a manifold of dimension $n$, and $g$ a metric
on $M$ with Levi-Civita connection $\nabla$ and holonomy 
group $G$. Then by Theorem \ref{l2thm1} there is a $k$-form 
$\vp$ on $M$ with $\nabla\vp=0$, corresponding to $\vp_0$. 
Hence $\d\vp=0$, and $\vp$ is closed. Also, the condition 
$\vp_0\vert_U\le\vol_U$ for all oriented $k$-planes $U$ in $\R^n$ 
implies that $\vp\vert_V\le\vol_V$ for all oriented tangent 
$k$-planes $V$ in $M$. Thus $\vp$ is a {\it calibration} on~$M$. 

At each point $x\in M$ the family of oriented tangent $k$-planes 
$V$ with $\vp\vert_V=\vol_V$ is isomorphic to ${\cal F}_{\vp_0}$,
which is `reasonably large'. This suggests that locally there 
should exist many $\vp$-submanifolds $N$ in $M$, so the calibrated 
geometry of $\vp$ on $(M,g)$ is nontrivial.

This gives us a general method for finding interesting 
calibrations on manifolds with reduced holonomy. Here are 
the most important examples of this.

\begin{itemize}
\item Let $G=\U(m)\subset{\rm O}(2m)$. Then $G$ preserves a 2-form
$\om_0$ on $\R^{2m}$. If $g$ is a metric on $M$ with holonomy
$\U(m)$ then $g$ is {\it K\"ahler} with complex structure $J$, and 
the 2-form $\om$ on $M$ associated to $\om_0$ is the {\it K\"ahler 
form} of $g$. 

One can show that $\om$ is a calibration on $(M,g)$, and the 
calibrated submanifolds are exactly the {\it holomorphic curves} 
in $(M,J)$. More generally $\om^k/k!$ is a calibration on $M$ for 
$1\le k\le m$, and the corresponding calibrated submanifolds are 
the complex $k$-dimensional submanifolds of~$(M,J)$.

\item Let $G=\SU(m)\subset{\rm O}(2m)$. Compact manifolds $(M,g)$ 
with holonomy $\SU(m)$ extend to {\it Calabi--Yau $m$-folds}
$(M,J,g,\Om)$, as in \S\ref{l45}. The real part $\Re\Om$ is a 
calibration on $M$, and the corresponding calibrated submanifolds 
are called {\it special Lagrangian submanifolds}.

\item The group $G_2\subset{\rm O}(7)$ preserves a 3-form $\vp_0$ and 
a 4-form $*\vp_0$ on $\R^7$. Thus a Riemannian 7-manifold $(M,g)$ with 
holonomy $G_2$ comes with a 3-form $\vp$ and 4-form $*\vp$, which are
both calibrations. The corresponding calibrated submanifolds are
called {\it associative $3$-folds} and {\it coassociative $4$-folds}.

\item The group $\Spin(7)\subset{\rm O}(8)$ preserves a 4-form $\Om_0$
on $\R^8$. Thus a Riemannian 8-manifold $(M,g)$ with holonomy Spin(7) 
has a 4-form $\Om$, which is a calibration. We call $\Om$-submanifolds 
{\it Cayley $4$-folds}.
\end{itemize}

It is an important general principle that to each calibration
$\vp$ on an $n$-manifold $(M,g)$ with special holonomy we
construct in this way, there corresponds a constant calibration
$\vp_0$ on $\R^n$. Locally, $\vp$-submanifolds in $M$ will look 
very like $\vp_0$-submanifolds in $\R^n$, and have many of the
same properties. Thus, to understand the calibrated submanifolds 
in a manifold with special holonomy, it is often a good idea to 
start by studying the corresponding calibrated submanifolds of~$\R^n$. 

In particular, singularities of $\vp$-submanifolds in $M$ will be 
locally modelled on singularities of $\vp_0$-submanifolds in $\R^n$. 
(Formally, the {\it tangent cone} at a singular point of a 
$\vp$-submanifold in $M$ is a conical $\vp_0$-submanifold in $\R^n$.)
So by studying singular $\vp_0$-submanifolds in $\R^n$, we may 
understand the singular behaviour of $\vp$-submanifolds in~$M$.

\subsection{Exercises}
\label{l65}

\begin{question}The metric $g$ and K\"ahler form $\om$ on $\C^m$ are 
given by
\begin{equation*}
g=\ms{\d z_1}+\cdots+\ms{\d z_m}\quad\text{and}\quad
\om=\frac{i}{2}(\d z_1\w\d\bar z_1+\cdots+\d z_m\w\d\bar z_m).
\end{equation*}
Show that a tangent 2-plane in $\C^m$ is calibrated w.r.t.\ $\om$
if and only if it is a complex line in $\C^m$. (Harder) 
generalize to tangent $2k$-planes and~$\frac{1}{k!}\,\om^k$.
\end{question}

\section{Special Lagrangian submanifolds in $\C^m$}
\label{l7}

We now discuss special Lagrangian submanifolds in $\C^m$. A reference 
for this section is Harvey and Lawson~\cite[\S III.1--\S III.2]{HaLa}.

\begin{dfn} Let $\C^m\cong\R^{2m}$ have complex coordinates 
$(z_1,\dots,z_m)$ and complex structure $I$, and define a metric 
$g$, K\"ahler form $\om$ and complex volume form $\Om$ on $\C^m$ by
\e
\begin{split}
g=\ms{\d z_1}+\cdots+\ms{\d z_m},\quad
\om&=\frac{i}{2}(\d z_1\w\d\bar z_1+\cdots+\d z_m\w\d\bar z_m),\\
\text{and}\quad\Om&=\d z_1\w\cdots\w\d z_m.
\end{split}
\label{l7eq1}
\e
Then $\Re\Om$ and $\Im\Om$ are real $m$-forms on $\C^m$. Let
$L$ be an oriented real submanifold of $\C^m$ of real dimension 
$m$. We call $L$ a {\it special Lagrangian submanifold\/} of
$\C^m$, or {\it SL $m$-fold\/} for short, if $L$ is calibrated 
with respect to $\Re\Om$, in the sense of~\S\ref{l62}.
\end{dfn}

In fact there is a more general definition involving a {\it phase}
${\rm e}^{i\th}$: if $\th\in[0,2\pi)$, we say that $L$ is {\it special
Lagrangian with phase} ${\rm e}^{i\th}$ if it is calibrated with
respect to $\cos\th\,\Re\Om+\sin\th\,\Im\Om$. But we will not use
this.

We shall identify the family $\cal F$ of tangent $m$-planes in $\C^m$
calibrated with respect to $\Re\Om$. The subgroup of $\GL(2m,\R)$ 
preserving $g,\om$ and $\Om$ is the Lie group $\SU(m)$ of complex 
unitary matrices with determinant 1. Define a real vector subspace 
$U$ in $\C^m$ to be
\e
U=\bigl\{(x_1,\ldots,x_m):x_j\in\R\bigr\}\subset\C^m,
\label{l7eq2}
\e
and let $U$ have the usual orientation. Then $U$ is calibrated
w.r.t.~$\Re\Om$. 

Furthermore, any oriented real vector subspace $V$ in $\C^m$ calibrated 
w.r.t.\ $\Re\Om$ is of the form $V=\ga\cdot U$ for some $\ga\in\SU(m)$.
Therefore $\SU(m)$ acts transitively on $\cal F$. The stabilizer subgroup
of $U$ in $\SU(m)$ is the subset of matrices in $\SU(m)$ with real entries,
which is $\SO(m)$. Thus ${\cal F}\cong\SU(m)/\SO(m)$, and we prove:

\begin{prop} The family $\cal F$ of oriented real\/ $m$-dimensional
vector subspaces $V$ in $\C^m$ with\/ $\Re\Om\vert_V=\vol_V$ is
isomorphic to $\SU(m)/\SO(m)$, and has dimension~$\ha(m^2+m-2)$.
\label{l7prop1}
\end{prop}

The dimension follows because $\dim\SU(m)=m^2-1$ and $\dim\SO(m)=\ha m(m-1)$.
It is easy to see that $\om\vert_U=\Im\Om\vert_U=0$. As $\SU(m)$
preserves $\om$ and $\Im\Om$ and acts transitively on $\cal F$, it
follows that $\om\vert_V=\Im\Om\vert_V=0$ for any $V\in{\cal F}$.
Conversely, if $V$ is a real $m$-dimensional vector subspace of $\C^m$ 
and $\om\vert_V=\Im\Om\vert_V=0$, then $V$ lies in $\cal F$, with some 
orientation. This implies an alternative characterization of special 
Lagrangian submanifolds, \cite[Cor.~III.1.11]{HaLa}:

\begin{prop} Let\/ $L$ be a real\/ $m$-dimensional submanifold 
of\/ $\C^m$. Then $L$ admits an orientation making it into a
special Lagrangian submanifold of\/ $\C^m$ if and only if\/ 
$\om\vert_L\equiv 0$ and\/~$\Im\Om\vert_L\equiv 0$.
\label{l7prop2}
\end{prop}

Note that an $m$-dimensional submanifold $L$ in $\C^m$ is 
called {\it Lagrangian} if $\om\vert_L\equiv 0$. (This is a term
from symplectic geometry, and $\om$ is a symplectic structure.) 
Thus special Lagrangian submanifolds are Lagrangian submanifolds 
satisfying the extra condition that $\Im\Om\vert_L\equiv 0$, which 
is how they get their name.

\subsection{Special Lagrangian 2-folds in $\C^2$ and the quaternions}
\label{l71}

The smallest interesting dimension, $m=2$, is a special case.
Let $\C^2$ have complex coordinates $(z_1,z_2)$, complex
structure $I$, and metric $g$, K\"ahler form $\om$ and holomorphic 
2-form $\Om$ defined in \eq{l7eq1}. Define real coordinates 
$(x_0,x_1,x_2,x_3)$ on $\C^2\cong\R^4$ by $z_0=x_0+ix_1$, 
$z_1=x_2+ix_3$. Then
\begin{alignat*}{2}
g&=\d x_0^2+\cdots+\d x_3^2,&\qquad
\om&=\d x_0\w\d x_1+\d x_2\w\d x_3,\\
\Re\Om&=\d x_0\w\d x_2-\d x_1\w\d x_3&\quad\text{and}\quad
\Im\Om&=\d x_0\w\d x_3+\d x_1\w\d x_2.
\end{alignat*}
Now define a {\it different} set of complex coordinates $(w_1,w_2)$ 
on $\C^2=\R^4$ by $w_1=x_0+ix_2$ and $w_2=x_1-ix_3$. 
Then~$\om-i\Im\Om=\d w_1\w\d w_2$. 

But by Proposition \ref{l7prop2}, a real 2-submanifold $L\subset\R^4$ 
is special Lagrangian if and only if $\om\vert_L\equiv\Im\Om\vert_L
\equiv 0$. Thus, $L$ is special Lagrangian if and only if 
$(\d w_1\w\d w_2)\vert_L\equiv 0$. But this holds if and only if
$L$ is a {\it holomorphic curve} with respect to the complex 
coordinates~$(w_1,w_2)$. 

Here is another way to say this. There are {\it two different}
complex structures $I$ and $J$ involved in this problem,
associated to the two different complex coordinate systems 
$(z_1,z_2)$ and $(w_1,w_2)$ on $\R^4$. In the coordinates 
$(x_0,\ldots,x_3)$, $I$ and $J$ are given by
\begin{alignat*}{4}
I\bigl({\ts\frac{\pd}{\pd x_0}}\bigr)&={\ts\frac{\pd}{\pd x_1}},\quad &
I\bigl({\ts\frac{\pd}{\pd x_1}}\bigr)&=-{\ts\frac{\pd}{\pd x_0}},\quad &
I\bigl({\ts\frac{\pd}{\pd x_2}}\bigr)&={\ts\frac{\pd}{\pd x_3}},\quad &
I\bigl({\ts\frac{\pd}{\pd x_3}}\bigr)&=-{\ts\frac{\pd}{\pd x_2}},\\
J\bigl({\ts\frac{\pd}{\pd x_0}}\bigr)&={\ts\frac{\pd}{\pd x_2}},\quad &
J\bigl({\ts\frac{\pd}{\pd x_1}}\bigr)&=-{\ts\frac{\pd}{\pd x_3}},\quad &
J\bigl({\ts\frac{\pd}{\pd x_2}}\bigr)&=-{\ts\frac{\pd}{\pd x_0}},\quad &
J\bigl({\ts\frac{\pd}{\pd x_3}}\bigr)&={\ts\frac{\pd}{\pd x_1}}.
\end{alignat*}
The usual complex structure on $\C^2$ is $I$, but a 2-fold $L$ in
$\C^2$ is special Lagrangian if and only if it is holomorphic
w.r.t.\ the alternative complex structure $J$. This means that 
special Lagrangian 2-folds are already very well understood,
so we generally focus our attention on dimensions~$m\ge 3$.

We can express all this in terms of the {\it quaternions} $\H$\,.
The complex structures $I,J$ anticommute, so that $IJ=-JI$, and 
$K=IJ$ is also a complex structure on $\R^4$, and $\an{1,I,J,K}$ 
is an algebra of automorphisms of $\R^4$ isomorphic to~$\H$\,.

\subsection{Special Lagrangian submanifolds in $\C^m$ as graphs}
\label{l72}

In symplectic geometry, there is a well-known way of manufacturing 
{\it Lagrangian} submanifolds of $\R^{2m}\cong\C^m$, which works as 
follows. Let $f:\R^m\ra\R$ be a smooth function, and define
\begin{equation*}
\Ga_f\!=\!\bigl\{\bigl(x_1\!+\!i{\ts\frac{\pd f}{\pd x_1}}(x_1,\ldots,x_m),
\ldots,x_m\!+\!i{\ts\frac{\pd f}{\pd x_m}}(x_1,\ldots,x_m)\bigr):
x_1,\ldots,x_m\!\in\!\R\bigr\}.
\end{equation*}
Then $\Ga_f$ is a smooth real $m$-dimensional submanifold of $\C^m$,
with $\om\vert_{\Ga_f}\equiv 0$. Identifying $\C^m\cong\R^{2m}\cong
\R^m\t(\R^m)^*$, we may regard $\Ga_f$ as the graph of the 1-form
$\d f$ on $\R^m$, so that $\Ga_f$ is the {\it graph of a closed\/ 
$1$-form}. Locally, but not globally, every Lagrangian submanifold 
arises from this construction.

Now by Proposition \ref{l7prop2}, a special Lagrangian $m$-fold
in $\C^m$ is a Lagrangian $m$-fold $L$ satisfying the additional
condition that $\Im\Om\vert_L\equiv 0$. We shall find the
condition for $\Ga_f$ to be a special Lagrangian $m$-fold. Define
the {\it Hessian} $\Hess f$ of $f$ to be the $m\t m$ matrix 
$\bigl(\frac{\pd^2f}{\pd x_i\pd x_j}\bigr)_{i,j=1}^m$ of real 
functions on $\R^m$. Then it is easy to show that $\Im\Om
\vert_{\Ga_f}\equiv 0$ if and only if
\e
\ts\Im\det_{\sst\mathbb C}\bigl(I+i\Hess f\bigr)\equiv 0
\quad\text{on $\C^m$.}
\label{l7eq3}
\e
This is a {\it nonlinear second-order elliptic partial differential
equation} upon the function~$f:\R^m\ra\R$\,. 

\subsection{Local discussion of special Lagrangian deformations}
\label{l73}

Suppose $L_0$ is a special Lagrangian submanifold in $\C^m$ (or,
more generally, in some Calabi--Yau $m$-fold). What can we say
about the family of {\it special Lagrangian deformations} of $L_0$,
that is, the set of special Lagrangian $m$-folds $L$ that are
`close to $L_0$' in a suitable sense? Essentially, deformation
theory is one way of thinking about the question `how many special 
Lagrangian submanifolds are there in~$\C^m$'?

Locally (that is, in small enough open sets), every special Lagrangian
$m$-fold looks quite like $\R^m$ in $\C^m$. Therefore deformations of 
special Lagrangian $m$-folds should look like special Lagrangian 
deformations of $\R^m$ in $\C^m$. So, we would like to know what
special Lagrangian $m$-folds $L$ in $\C^m$ close to $\R^m$ look like.

Now $\R^m$ is the graph $\Ga_f$ associated to the function $f\equiv 0$.
Thus, a graph $\Ga_f$ will be close to $\R^m$ if the function $f$ and 
its derivatives are small. But then $\Hess f$ is small, so we can 
approximate equation \eq{l7eq3} by its {\it linearization}. For
\begin{equation*}
\ts\Im\det_{\sst\mathbb C}\bigl(I+i\Hess f\bigr)=
\Tr\Hess f+\text{higher order terms}.
\end{equation*}
Thus, when the second derivatives of $f$ are small, equation
\eq{l7eq3} reduces approximately to $\Tr\Hess f\equiv 0$. But 
$\Tr\Hess f=\frac{\pd^2f}{(\pd x_1)^2}+\cdots+\frac{\pd^2f}{(\pd x_m)^2}
=\De f$, where $\De$ is the {\it Laplacian} on~$\R^m$. 

Hence, the small special Lagrangian deformations of $\R^m$ in
$\C^m$ are approximately parametrized by small {\it harmonic 
functions} on $\R^m$. Actually, because adding a constant to
$f$ has no effect on $\Ga_f$, this parametrization is degenerate.
We can get round this by parametrizing instead by $\d f$, which
is a closed and coclosed 1-form. This justifies the following:
\medskip

\noindent{\bf Principle.} {\it Small special Lagrangian deformations 
of a special Lagrangian $m$-fold\/ $L$ are approximately parametrized 
by closed and coclosed\/ $1$-forms $\al$ on~$L$.}
\medskip

\noindent This is the idea behind McLean's Theorem, Theorem 
\ref{l9thm1} below.

We have seen using \eq{l7eq3} that the deformation problem for special 
Lagrangian $m$-folds can be written as an {\it elliptic equation}.
In particular, there are the same number of equations as functions, so
the problem is neither overdetermined nor underdetermined. Therefore we
do not expect special Lagrangian $m$-folds to be very few and very rigid
(as would be the case if \eq{l7eq3} were overdetermined), nor to 
be very abundant and very flabby (as would be the case if 
\eq{l7eq3} were underdetermined).

If we think about Proposition \ref{l7prop1} for a while, this may
seem surprising. For the set $\cal F$ of special Lagrangian $m$-planes
in $\C^m$ has dimension $\ha(m^2+m-2)$, but the set of all real $m$-planes
in $\C^m$ has dimension $m^2$. So the special Lagrangian $m$-planes
have codimension $\ha(m^2-m+2)$ in the set of all $m$-planes.

This means that the condition for a real $m$-submanifold $L$ in
$\C^m$ to be special Lagrangian is $\ha(m^2-m+2)$ real equations on 
each tangent space of $L$. However, the freedom to vary $L$ is the
sections of its normal bundle in $\C^m$, which is $m$ real functions. 
When $m\ge 3$, there are more equations than functions, so we would 
expect the deformation problem to be {\it overdetermined}.

The explanation is that because $\om$ is a {\it closed}\/ 2-form,
submanifolds $L$ with $\om\vert_L\equiv 0$ are much more abundant
than would otherwise be the case. So the closure of $\om$ is a
kind of integrability condition necessary for the existence of 
many special Lagrangian submanifolds, just as the integrability 
of an almost complex structure is a necessary condition for the
existence of many complex submanifolds of dimension greater than 
1 in a complex manifold.

\subsection{Exercises}
\label{l74}

\begin{question} Find your own proofs of Propositions 
\ref{l7prop1} and~\ref{l7prop2}.
\end{question}

\section{Constructions of SL $m$-folds in $\C^m$}
\label{l8}

We now describe five methods of constructing special Lagrangian
$m$-folds in $\C^m$, drawn from papers by the author
\cite{Joyc3,Joyc4,Joyc5,Joyc6,Joyc7,Joyc9,Joyc10,Joyc11},
Bryant \cite{Brya}, Castro and Urbano \cite{CaUr}, Goldstein
\cite{Gold1,Gold2}, Harvey \cite[p.~139--143]{Harv}, Harvey and
Lawson \cite[\S III]{HaLa}, Haskins \cite{Hask}, Lawlor
\cite{Lawl}, Ma and Ma \cite{MaMa}, McIntosh \cite{McIn} and
Sharipov \cite{Shar}. These yield many examples of singular SL
$m$-folds, and so hopefully will help in understanding what
general singularities of SL $m$-folds in Calabi--Yau $m$-folds
are like.

\subsection{SL $m$-folds with large symmetry groups}
\label{l81}

Here is a method used in \cite{Joyc3} (and also by Harvey and
Lawson \cite[\S III.3]{HaLa}, Haskins \cite{Hask} and Goldstein
\cite{Gold1,Gold2}) to construct examples of SL $m$-folds in 
$\C^m$. The group $\SU(m)\lt\C^m$ acts on $\C^m$ preserving 
all the structure $g,\om,\Om$, so that it takes SL $m$-folds 
to SL $m$-folds in $\C^m$. Let $G$ be a Lie subgroup of 
$\SU(m)\lt\C^m$ with Lie algebra $\g$, and $N$ a connected 
$G$-invariant SL $m$-fold in~$\C^m$.

Since $G$ preserves the symplectic form $\om$ on $\C^m$, one 
can show that it has a {\it moment map} $\mu:\C^m\ra\g^*$.
As $N$ is Lagrangian, one can show that $\mu$ is constant
on $N$, that is, $\mu\equiv c$ on $N$ for some $c\in Z(\g^*)$, 
the {\it centre} of~$\g^*$.

If the orbits of $G$ in $N$ are of codimension 1 (that is,
dimension $m-1$), then $N$ is a 1-parameter family of 
$G$-orbits ${\cal O}_t$ for $t\in\R$\,. After reparametrizing 
the variable $t$, it can be shown that the special Lagrangian 
condition is equivalent to an o.d.e.\ in $t$ upon the 
orbits~${\cal O}_t$.

Thus, we can construct examples of cohomogeneity one SL $m$-folds
in $\C^m$ by solving an o.d.e.\ in the family of $(m-1)$-dimensional 
$G$-orbits $\cal O$ in $\C^m$ with $\mu\vert_{\cal O}\equiv c$,
for fixed $c\in Z(\g^*)$. This o.d.e.\ usually turns out to be 
{\it integrable}.

Now suppose $N$ is a {\it special Lagrangian cone} in $\C^m$, 
invariant under a subgroup $G\subset\SU(m)$ which has orbits of 
dimension $m-2$ in $N$. In effect the symmetry group of $N$ is 
$G\t\R_+$, where $\R_+$ acts by {\it dilations}, as $N$ is a cone. 
Thus, in this situation too the symmetry group of $N$ acts with 
cohomogeneity one, and we again expect the problem to reduce to 
an o.d.e.

One can show that $N\cap{\cal S}^{2m-1}$ is a 1-parameter family of 
$G$-orbits ${\cal O}_t$ in ${\cal S}^{2m-1}\cap\mu^{-1}(0)$ satisfying 
an o.d.e. By solving this o.d.e.\ we construct SL cones in $\C^m$.
When $G=\U(1)^{m-2}$, the o.d.e. has many {\it periodic solutions} 
which give large families of distinct SL cones on $T^{m-1}$. In 
particular, we can find many examples of SL $T^2$-cones in~$\C^3$.

\subsection{Evolution equations for SL $m$-folds}
\label{l82}

The following method was used in \cite{Joyc4} and \cite{Joyc5} to 
construct many examples of SL $m$-folds in $\C^m$. A related but 
less general method was used by Lawlor \cite{Lawl}, and completed 
by Harvey~\cite[p.~139--143]{Harv}.

Let $P$ be a real analytic $(m-1)$-dimensional manifold, and 
$\chi$ a nonvanishing real analytic section of $\La^{m-1}TP$.
Let $\{\phi_t:t\in\R\}$ be a 1-parameter family of real analytic 
maps $\phi_t:P\ra\C^m$. Consider the o.d.e.
\e
\left(\frac{\d\phi_t}{\d t}\right)^b=(\phi_t)_*(\chi)^{a_1\ldots 
a_{m-1}}(\Re\Om)_{a_1\ldots a_{m-1}a_m}g^{a_mb},
\label{l8eq1}
\e
using the index notation for (real) tensors on $\C^m$, where 
$g^{ab}$ is the inverse of the Euclidean metric $g_{ab}$ on~$\C^m$.

It is shown in \cite[\S 3]{Joyc4} that if the $\phi_t$ satisfy 
\eq{l8eq1} and $\phi_0^*(\om)\equiv 0$, then $\phi_t^*(\om)\equiv 0$ 
for all $t$, and $N=\bigl\{\phi_t(p):p\in P$, $t\in\R\bigr\}$ is 
an SL $m$-fold in $\C^m$ wherever it is nonsingular. We think of 
\eq{l8eq1} as an {\it evolution equation}, and $N$ as the result 
of evolving a 1-parameter family of $(m\!-\!1)$-submanifolds
$\phi_t(P)$ in~$\C^m$.

Here is one way to understand this result. Suppose we are given 
$\phi_t:P\ra\C^m$ for some $t$, and we want to find an SL $m$-fold
$N$ in $\C^m$ containing the $(m\!-\!1)$-submanifold $\phi_t(P)$. 
As $N$ is Lagrangian, a necessary condition for this is that $\om
\vert_{\phi_t(P)}\equiv 0$, and hence $\phi_t^*(\om)\equiv 0$ on~$P$. 

The effect of equation \eq{l8eq1} is to flow $\phi_t(P)$ in the 
direction in which $\Re\Om$ is `largest'. The result is that $\Re\Om$ 
is `maximized' on $N$, given the initial conditions. But $\Re\Om$ is
maximal on $N$ exactly when $N$ is calibrated w.r.t.\ $\Re\Om$, that
is, when $N$ is special Lagrangian. The same technique also works for 
other calibrations, such as the associative and coassociative 
calibrations on $\R^7$, and the Cayley calibration on~$\R^8$.

Now \eq{l8eq1} evolves amongst the infinite-dimensional family of
real analytic maps $\phi:P\ra\C^m$ with $\phi^*(\om)\equiv 0$, so
it is an {\it infinite-dimensional} problem, and thus difficult
to solve explicitly. However, there are {\it finite-dimensional} 
families $\cal C$ of maps $\phi:P\ra\C^m$ such that evolution 
stays in $\cal C$. This gives a {\it finite-dimensional} o.d.e., 
which can hopefully be solved fairly explicitly. For example, 
if we take $G$ to be a Lie subgroup of $\SU(m)\lt\C^m$, $P$ to 
be an $(m\!-\!1)$-dimensional homogeneous space $G/H$, and 
$\phi:P\ra\C^m$ to be $G$-equivariant, we recover the 
construction of~\S\ref{l81}. 

But there are also other possibilities for $\cal C$ which do not 
involve a symmetry assumption. Suppose $P$ is a submanifold of 
$\R^n$, and $\chi$ the restriction to $P$ of a linear or affine 
map $\R^n\ra\La^{m-1}\R^n$. (This is a strong condition on $P$
and $\chi$.) Then we can take $\cal C$ to be the set of 
restrictions to $P$ of linear or affine maps~$\R^n\ra\C^m$.

For instance, set $m=n$ and let $P$ be a quadric in $\R^m$. Then
one can construct SL $m$-folds in $\C^m$ with few symmetries by
evolving quadrics in Lagrangian planes $\R^m$ in $\C^m$. When $P$ 
is a quadric cone in $\R^m$ this gives many SL cones on products
of spheres~${\cal S}^a\t{\cal S}^b\t{\cal S}^1$.

\subsection{Ruled special Lagrangian 3-folds}
\label{l83}

A 3-submanifold $N$ in $\C^3$ is called {\it ruled} if it is fibred by 
a 2-dimensional family $\cal F$ of real lines in $\C^3$. A {\it cone} 
$N_0$ in $\C^3$ is called {\it two-sided} if $N_0=-N_0$. Two-sided 
cones are automatically ruled. If $N$ is a ruled 3-fold in $\C^3$, 
we define the {\it asymptotic cone} $N_0$ of $N$ to be the two-sided 
cone fibred by the lines passing through 0 and parallel to those 
in~$\cal F$. 

Ruled SL 3-folds are studied in \cite{Joyc6}, and also by
Harvey and Lawson \cite[\S III.3.C, \S III.4.B]{HaLa} and
Bryant \cite[\S 3]{Brya}. Each (oriented) real line in $\C^3$ 
is determined by its {\it direction} in ${\cal S}^5$ together 
with an orthogonal {\it translation} from the origin. Thus a 
ruled 3-fold $N$ is determined by a 2-dimensional family of 
directions and translations. 

The condition for $N$ to be special Lagrangian turns out 
\cite[\S 5]{Joyc6} to reduce to two equations, the first 
involving only the direction components, and the second 
{\it linear} in the translation components. Hence, if a 
ruled 3-fold $N$ in $\C^3$ is special Lagrangian, then so is 
its asymptotic cone $N_0$. Conversely, the ruled SL 3-folds 
$N$ asymptotic to a given two-sided SL cone $N_0$ come from 
solutions of a linear equation, and so form a {\it vector space}.

Let $N_0$ be a two-sided SL cone, and let $\Si=N_0\cap{\cal S}^5$. 
Then $\Si$ is a {\it Riemann surface}. Holomorphic vector fields on 
$\Si$ give solutions to the linear equation (though not all solutions)
\cite[\S 6]{Joyc6}, and so yield new ruled SL 3-folds. In particular, 
each SL $T^2$-cone gives a 2-dimensional family of ruled SL 3-folds, 
which are generically diffeomorphic to $T^2\t\R$ as immersed 
3-submanifolds. 

\subsection{Integrable systems}
\label{l84}

Let $N_0$ be a special Lagrangian cone in $\C^3$, and set
$\Si=N_0\cap{\cal S}^5$. As $N_0$ is calibrated, it is minimal
in $\C^3$, and so $\Si$ is minimal in ${\cal S}^5$. That is, 
$\Si$ is a {\it minimal Legendrian surface} in ${\cal S}^5$.
Let $\pi:{\cal S}^5\ra\CP^2$ be the Hopf projection. One can 
also show that $\pi(\Si)$ is a {\it minimal Lagrangian surface}
in~$\CP^2$.

Regard $\Si$ as a {\it Riemann surface}. Then the inclusions 
$\iota:\Si\ra{\cal S}^5$ and $\pi\circ\iota:\Si\ra\CP^2$ 
are {\it conformal harmonic maps}. Now harmonic maps from 
Riemann surfaces into ${\cal S}^n$ and $\CP^m$ are an 
{\it integrable system}. There is a complicated theory 
for classifying them in terms of algebro-geometric `spectral 
data', and finding `explicit' solutions. In principle, this
gives all harmonic maps from $T^2$ into ${\cal S}^n$ and
$\CP^m$. So, the field of integrable systems offers the hope
of a {\it classification} of all SL $T^2$-cones in~$\C^3$.

For a good general introduction to this field, see Fordy and
Wood \cite{FoWo}. Sharipov \cite{Shar} and Ma and Ma \cite{MaMa}
apply this integrable systems machinery to describe minimal
Legendrian tori in ${\cal S}^5$, and minimal Lagrangian tori
in $\CP^2$, respectively, giving explicit formulae in terms
of Prym theta functions. McIntosh \cite{McIn} provides a more
recent, readable, and complete discussion of special Lagrangian
cones in $\C^3$ from the integrable systems perspective.

The families of SL $T^2$-cones constructed by $\U(1)$-invariance 
in \S\ref{l81}, and by evolving quadrics in \S\ref{l82}, turn out
to come from a more general, very explicit, `integrable systems' family 
of conformal harmonic maps $\R^2\ra{\cal S}^5$ with Legendrian image, 
involving two commuting, integrable o.d.e.s., described in \cite{Joyc7}.
So, we can fit some of our examples into the integrable systems framework.

However, we know a good number of other constructions of SL $m$-folds 
in $\C^m$ which have the classic hallmarks of integrable systems
--- elliptic functions, commuting o.d.e.s, and so on --- but which are 
not yet understood from the point of view of integrable systems. I would 
like to ask the integrable systems community: do SL $m$-folds in $\C^m$ 
for $m\ge 3$, or at least some classes of such submanifolds, constitute 
some kind of higher-dimensional integrable system?

\subsection{Analysis and $\U(1)$-invariant SL $3$-folds in $\C^3$}
\label{l85}

Next we summarize the author's three papers \cite{Joyc9,Joyc10,Joyc11},
which study SL 3-folds $N$ in $\C^3$ invariant under the $\U(1)$-action
\e
{\rm e}^{i\th}:(z_1,z_2,z_3)\mapsto
({\rm e}^{i\th}z_1,{\rm e}^{-i\th}z_2,z_3)
\quad\text{for ${\rm e}^{i\th}\in\U(1)$.}
\label{l8eq2}
\e
These papers are briefly surveyed in \cite{Joyc12}. Locally we can
write $N$ in the form
\e
\begin{split}
N=\bigl\{(z_1,z_2,z_3)\in\C^3:\,& z_1z_2=v(x,y)+iy,\quad z_3=x+iu(x,y),\\
&\ms{z_1}-\ms{z_2}=2a,\quad (x,y)\in S\bigr\},
\end{split}
\label{l8eq3}
\e
where $S$ is a domain in $\R^2$, $a\in\R$ and $u,v:S\ra\R$ are
continuous.

Here we may take $\ms{z_1}-\ms{z_2}=2a$ to be one of the equations
defining $N$ as $\ms{z_1}-\ms{z_2}$ is the {\it moment map} of the
$\U(1)$-action \eq{l8eq2}, and so $\ms{z_1}-\ms{z_2}$ is constant
on any $\U(1)$-invariant Lagrangian 3-fold in $\C^3$. Effectively
\eq{l8eq3} just means that we are choosing $x=\Re(z_3)$ and
$y=\Im(z_1z_2)$ as local coordinates on the 2-manifold $N/\U(1)$.
Then we find~\cite[Prop.~4.1]{Joyc9}:

\begin{prop} Let\/ $S,a,u,v$ and\/ $N$ be as above. Then
\begin{itemize}
\item[{\rm(a)}] If\/ $a=0$, then $N$ is a (possibly singular) SL\/
$3$-fold in $\C^3$ if\/ $u,v$ are differentiable and satisfy
\e
\frac{\pd u}{\pd x}=\frac{\pd v}{\pd y}
\quad\text{and}\quad
\frac{\pd v}{\pd x}=-2\bigl(v^2+y^2\bigr)^{1/2}\frac{\pd u}{\pd y},
\label{l8eq4}
\e
except at points $(x,0)$ in $S$ with\/ $v(x,0)=0$, where $u,v$ 
need not be differentiable. The singular points of\/ $N$ are those
of the form $(0,0,z_3)$, where $z_3=x+iu(x,0)$ for $(x,0)\in S$ 
with\/~$v(x,0)=0$.
\item[{\rm(b)}] If\/ $a\ne 0$, then $N$ is a nonsingular SL\/ $3$-fold
in $\C^3$ if and only if\/ $u,v$ are differentiable in $S$ and
satisfy
\e
\frac{\pd u}{\pd x}=\frac{\pd v}{\pd y}\quad\text{and}\quad
\frac{\pd v}{\pd x}=-2\bigl(v^2+y^2+a^2\bigr)^{1/2}\frac{\pd u}{\pd y}.
\label{l8eq5}
\e
\end{itemize}
\label{l8prop1}
\end{prop}

Now \eq{l8eq4} and \eq{l8eq5} are {\it nonlinear Cauchy--Riemann
equations}. Thus, we may treat $u+iv$ as like a holomorphic function
of $x+iy$. Many of the results in \cite{Joyc9,Joyc10,Joyc11} are
analogues of well-known results in elementary complex analysis.

In \cite[Prop.~7.1]{Joyc9} we show that solutions $u,v\in C^1(S)$ 
of \eq{l8eq5} come from a potential $f\in C^2(S)$ satisfying a
second-order quasilinear elliptic equation.

\begin{prop} Let\/ $S$ be a domain in $\R^2$ and\/ $u,v\in C^1(S)$
satisfy \eq{l8eq5} for $a\ne 0$. Then there exists $f\in C^2(S)$
with\/ $\frac{\pd f}{\pd y}=u$, $\frac{\pd f}{\pd x}=v$ and
\e
P(f)=\Bigl(\Bigl(\frac{\pd f}{\pd x}\Bigr)^2+y^2+a^2
\Bigr)^{-1/2}\frac{\pd^2f}{\pd x^2}+2\,\frac{\pd^2f}{\pd y^2}=0.
\label{l8eq6}
\e
This $f$ is unique up to addition of a constant, $f\mapsto f+c$.
Conversely, all solutions of\/ \eq{l8eq6} yield solutions 
of\/~\eq{l8eq5}. 
\label{l8prop2}
\end{prop}

In the following result, a condensation of \cite[Th.~7.6]{Joyc9}
and \cite[Th.s 9.20 \& 9.21]{Joyc10}, we prove existence and
uniqueness for the {\it Dirichlet problem} for~\eq{l8eq6}.

\begin{thm} Suppose $S$ is a strictly convex domain in $\R^2$ invariant
under $(x,y)\mapsto(x,-y)$, and\/ $\al\in(0,1)$. Let\/ $a\in\R$ and\/
$\phi\in C^{3,\al}(\pd S)$. Then if\/ $a\ne 0$ there exists a unique
solution $f$ of\/ \eq{l8eq6} in $C^{3,\al}(S)$ with\/ $f\vert_{\pd S}
=\phi$. If\/ $a=0$ there exists a unique $f\in C^1(S)$ with\/
$f\vert_{\pd S}=\phi$, which is twice weakly differentiable and
satisfies \eq{l8eq6} with weak derivatives. Furthermore, the map
$C^{3,\al}(\pd S)\t\R\ra C^1(S)$ taking $(\phi,a)\mapsto f$ is
continuous.
\label{l8thm1}
\end{thm}

Here a domain $S$ in $\R^2$ is {\it strictly convex} if it is
convex and the curvature of $\pd S$ is nonzero at each point.
Also domains are by definition compact, with smooth boundary,
and $C^{3,\al}(\pd S)$ and $C^{3,\al}(S)$ are {\it H\"older
spaces} of functions on $\pd S$ and $S$. For more details
see~\cite{Joyc9,Joyc10}.

Combining Propositions \ref{l8prop1} and \ref{l8prop2} and
Theorem \ref{l8thm1} gives existence and uniqueness for a large
class of $\U(1)$-invariant SL 3-folds in $\C^3$, with boundary
conditions, and including {\it singular} SL 3-folds. It is
interesting that this existence and uniqueness is {\it entirely
unaffected} by singularities appearing in~$S^\circ$. 

Here are some other areas covered in \cite{Joyc9,Joyc10,Joyc11}.
Examples of solutions $u,v$ of \eq{l8eq4} and \eq{l8eq5} are
given in \cite[\S 5]{Joyc9}. In \cite{Joyc10} we give more
precise statements on the regularity of singular solutions of
\eq{l8eq4} and \eq{l8eq6}. In \cite[\S 6]{Joyc9} and \cite[\S
7]{Joyc11} we consider the zeroes of $(u_1,v_1)-(u_2,v_2)$,
where $(u_j,v_j)$ are (possibly singular) solutions of \eq{l8eq4}
and~\eq{l8eq5}.

We show that if $(u_1,v_1)\not\equiv(u_2,v_2)$ then the zeroes
of $(u_1,v_1)-(u_2,v_2)$ in $S^\circ$ are {\it isolated}, with
a positive integer {\it multiplicity}, and that the zeroes of
$(u_1,v_1)-(u_2,v_2)$ in $S^\circ$ can be counted with
multiplicity in terms of boundary data on $\pd S$. In particular,
under some boundary conditions we can show $(u_1,v_1)-(u_2,v_2)$
has no zeroes in $S^\circ$, so that the corresponding SL 3-folds
do not intersect. This will be important in constructing
$\U(1)$-invariant SL fibrations in~\S\ref{l115}.

In \cite[\S 9--\S 10]{Joyc11} we study singularities of solutions
$u,v$ of \eq{l8eq4}. We show that either $u(x,-y)\equiv u(x,y)$
and $v(x,-y)\equiv -v(x,y)$, so that $u,v$ are singular all
along the $x$-axis, or else the singular points of $u,v$ in
$S^\circ$ are all {\it isolated}, with a positive integer
{\it multiplicity}, and one of two {\it types}. We also show
that singularities exist with every multiplicity and type, and
multiplicity $n$ singularities occur in codimension $n$ in the
family of all $\U(1)$-invariant SL 3-folds.

\subsection{Examples of singular special Lagrangian 3-folds in $\C^3$}
\label{l86}

We shall now describe four families of SL 3-folds in $\C^3$, as examples 
of the material of \S\ref{l81}--\S\ref{l84}. They have been chosen to 
illustrate different kinds of singular behaviour of SL 3-folds, and 
also to show how nonsingular SL 3-folds can converge to a singular SL 
3-fold, to serve as a preparation for our discussion of singularities 
of SL $m$-folds in~\S\ref{l10}.

Our first example derives from Harvey and Lawson \cite[\S III.3.A]{HaLa},
and is discussed in detail in \cite[\S 3]{Joyc1} and~\cite[\S 4]{Joyc13}.

\begin{ex} Define a subset $L_0$ in $\C^3$ by
\begin{equation*}
L_0=\bigl\{(r{\rm e}^{i\th_1},r{\rm e}^{i\th_2},r{\rm e}^{i\th_3}):
r\ge 0,\quad \th_1,\th_2,\th_3\in\R,\quad \th_1+\th_2+\th_3=0\bigr\}.
\end{equation*}
Then $L_0$ is a {\it special Lagrangian cone} on $T^2$. An 
alternative definition is
\begin{equation*}
L_0=\bigl\{(z_1,z_2,z_3)\in\C^3:\md{z_1}=\md{z_2}=\md{z_3},\;
\Im(z_1z_2z_3)=0,\; \Re(z_1z_2z_3)\ge 0\bigr\}.
\end{equation*}

Let $t>0$, write ${\cal S}^1=\bigl\{{\rm e}^{i\th}:\th\in\R\bigr\}$, 
and define a map $\phi_t:{\cal S}^1\t\C\ra\C^3$ by
\begin{equation*}
\phi_t:(e^{i\th},z)\mapsto
\bigl((\md{z}^2+t^2)^{1/2}{\rm e}^{i\th},z,e^{-i\th}\bar z\bigr).
\end{equation*}
Then $\phi_t$ is an {\it embedding}. Define $L_t=\Image\phi_t$.
Then $L_t$ is a nonsingular special Lagrangian 3-fold in $\C^3$
diffeomorphic to ${\cal S}^1\t\R^2$. An equivalent definition~is
\begin{align*}
L_t=\bigl\{(z_1,z_2,z_3)\in\C^3:\,&\ms{z_1}-t^2=\ms{z_2}=\ms{z_3},\\
&\Im(z_1z_2z_3)=0,\quad \Re(z_1z_2z_3)\ge 0\bigr\}.
\end{align*}

As $t\ra 0_+$, the nonsingular SL 3-fold $L_t$ converges to the
singular SL cone $L_0$. Note that $L_t$ is {\it asymptotic}
to $L_0$ at infinity, and that $L_t=t\,L_1$ for $t>0$, so that
the $L_t$ for $t>0$ are all homothetic to each other. Also,
each $L_t$ for $t\ge 0$ is invariant under the $T^2$ subgroup 
of $\SU(3)$ acting by
\begin{equation*}
(z_1,z_2,z_3)\mapsto({\rm e}^{i\th_1}z_1,{\rm e}^{i\th_2}z_2,
{\rm e}^{i\th_3}z_3) \;\>
\text{for $\th_1,\th_2,\th_3\in\R$ with $\th_1+\th_2+\th_3=0$,}
\end{equation*}
and so fits into the framework of \S\ref{l81}. By
\cite[Th.~5.1]{Joyc9} the $L_a$ may also be written in the form
\eq{l8eq3} for continuous $u,v:\R^2\ra\R$, as in~\S\ref{l85}.
\label{l8ex1}
\end{ex}

Our second example is adapted from Harvey and 
Lawson~\cite[\S III.3.B]{HaLa}. 

\begin{ex} For each $t>0$, define
\begin{align*}
L_t=\bigl\{({\rm e}^{i\th}x_1,{\rm e}^{i\th}x_2,{\rm e}^{i\th}x_3):&\,
x_j\in\R,\quad \th\in(0,\pi/3),\\
&x_1^2+x_2^2+x_3^2=t^2(\sin 3\th)^{-2/3}\bigr\}.
\end{align*}
Then $L_t$ is a nonsingular embedded SL 3-fold in $\C^3$ diffeomorphic
to ${\cal S}^2\t\R$\,. As $t\ra 0_+$ it converges to the singular union
$L_0$ of the two SL 3-planes
\begin{equation*}
\Pi_1=\bigl\{(x_1,x_2,x_3):x_j\in\R\bigr\}\;\>\text{and}\;\>
\Pi_2=\bigl\{({\rm e}^{i\pi/3}x_1,{\rm e}^{i\pi/3}x_2,
{\rm e}^{i\pi/3}x_3):x_j\in\R\bigr\},
\end{equation*}
which intersect at 0. Note that $L_t$ is invariant under the action 
of the Lie subgroup $\SO(3)$ of $\SU(3)$, acting on $\C^3$ in the 
obvious way, so again this comes from the method of \S\ref{l81}. 
Also $L_t$ is asymptotic to $L_0$ at infinity.
\label{l8ex2}
\end{ex}

Our third example is taken from~\cite[Ex.~9.4 \& Ex.~9.5]{Joyc3}.

\begin{ex} Let $a_1,a_2$ be positive, coprime integers, and set
$a_3=-a_1-a_2$. Let $c\in\R$, and define
\begin{equation*}
L^{a_1,a_2}_c=\bigl\{({\rm e}^{ia_1\th}x_1,{\rm e}^{ia_2\th}x_2,
i{\rm e}^{ia_3\th}x_3):\th\in\R,\; x_j\in\R,\; 
a_1x_1^2+a_2x_2^2+a_3x_3^2=c\bigr\}.
\end{equation*}
Then $L^{a_1,a_2}_c$ is a special Lagrangian 3-fold, which comes
from the `evolving quadrics' construction of \S\ref{l82}. It is
also symmetric under the $\U(1)$-action
\begin{equation*}
(z_1,z_2,z_3)\mapsto({\rm e}^{ia_1\th}z_1,{\rm e}^{ia_2\th}z_2,
i{\rm e}^{ia_3\th}z_3)\quad\text{for $\th\in\R$,}
\end{equation*}
but this is not a necessary feature of the construction; these
are just the easiest examples to write down.

When $c=0$ and $a_3$ is odd, $L^{a_1,a_2}_0$ is an embedded 
special Lagrangian cone on $T^2$, with one singular point at 0. 
When $c=0$ and $a_3$ is even, $L^{a_1,a_2}_0$ is two opposite 
embedded SL $T^2$-cones with one singular point at~0.

When $c>0$ and $a_3$ is odd, $L^{a_1,a_2}_c$ is an embedded 3-fold 
diffeomorphic to a nontrivial real line bundle over the Klein bottle.
When $c>0$ and $a_3$ is even, $L^{a_1,a_2}_c$ is an embedded 3-fold 
diffeomorphic to $T^2\t\R$\,. In both cases, $L^{a_1,a_2}_c$ is
a {\it ruled}\/ SL 3-fold, as in \S\ref{l83}, since it is fibred
by hyperboloids of one sheet in $\R^3$, which are ruled in two
different ways.

When $c<0$ and $a_3$ is odd, $L^{a_1,a_2}_c$ an immersed
copy of ${\cal S}^1\t\R^2$. When $c<0$ and $a_3$ is even, 
$L^{a_1,a_2}_c$ two immersed copies of~${\cal S}^1\t\R^2$. 
\label{l8ex3}
\end{ex}

All the singular SL 3-folds we have seen so far have been {\it 
cones} in $\C^3$. Our final example, taken from \cite{Joyc5},
has more complicated singularities which are not cones. They
are difficult to describe in a simple way, so we will not say
much about them. For more details, see~\cite{Joyc5}.
                                                   
\begin{ex} In \cite[\S 5]{Joyc5} the author constructed a
family of maps $\Phi:\R^3\ra\C^3$ with special Lagrangian image
$N=\Image\Phi$. It is shown in \cite[\S 6]{Joyc5} that generic 
$\Phi$ in this family are immersions, so that $N$ is nonsingular
as an immersed SL 3-fold, but in codimension 1 in the family they 
develop isolated singularities.

Here is a rough description of these singularities, taken 
from \cite[\S 6]{Joyc5}. Taking the singular point to be at 
$\Phi(0,0,0)=0$, one can write $\Phi$ as
\e
\begin{split}
\Phi(x,y,t)=
&\bigl(x+{\ts\frac{1}{4}}g({\bf u},{\bf v})t^2\bigr)\,{\bf u}
+\bigl(y^2-{\ts\frac{1}{4}}\ms{{\bf u}}t^2\bigr)\,{\bf v}\\
&+2yt\,{\bf u}\t{\bf v}+O\bigl(x^2+\md{xy}+\md{xt}+\md{y}^3+\md{t}^3\bigr),
\end{split}
\label{l8eq7}
\e
where ${\bf u},{\bf v}$ are linearly independent vectors in $\C^3$ with 
$\om({\bf u},{\bf v})=0$, and $\t:\C^3\t\C^3\ra\C^3$ is defined by
\begin{equation*}
(r_1,r_2,r_3)\t(s_1,s_2,s_3)={\ts\frac{1}{2}}(\bar r_2\bar s_3-\bar r_3\bar 
s_2,\bar r_3\bar s_1-\bar r_1\bar s_3,\bar r_1\bar s_2-\bar r_2\bar s_1).
\end{equation*}
The next few terms in the expansion \eq{l8eq7} can also be given very 
explicitly, but we will not write them down as they are rather complex,
and involve further choices of vectors~${\bf w},{\bf x},\ldots$. 

What is going on here is that the lowest order terms in $\Phi$ are a
{\it double cover} of the special Lagrangian plane $\an{{\bf u},{\bf v},
{\bf u}\t{\bf v}}_{\sst\mathbb R}$ in $\C^3$, {\it branched} along the 
real line $\an{\bf u}_{\sst\mathbb R}$. The branching occurs when $y=t=0$. 
Higher order terms deviate from the 3-plane $\an{{\bf u},{\bf v},{\bf u}
\t{\bf v}}_{\sst\mathbb R}$, and make the singularity isolated.
\label{l8ex4}
\end{ex}

\subsection{Exercises}
\label{l87}

\begin{question}The group of automorphisms of $\C^m$ preserving 
$g,\om$ and $\Om$ is $\SU(m)\lt\C^m$, where $\C^m$ acts by 
translations. Let $G$ be a Lie subgroup of $\SU(m)\lt\C^m$, let 
$\g$ be its Lie algebra, and let $\phi:\g\ra\Vect(\C^m)$ be the 
natural map associating an element of $\g$ to the corresponding 
vector field on~$\C^m$.

A {\it moment map} for the action of $G$ on $\C^m$ is a smooth
map $\mu:\C^m\ra\g^*$, such that $\phi(x)\cdot\om=x\cdot\d\mu$ 
for all $x\in\g$, and $\mu:\C^m\ra\g^*$ is equivariant with 
respect to the $G$-action on $\C^m$ and the coadjoint $G$-action 
on $\g^*$. Moment maps always exist if $G$ is compact or 
semisimple, and are unique up to the addition of a constant 
in the centre $Z(\g^*)$ of $\g^*$, that is, the $G$-invariant
subspace of~$\g^*$.

Suppose $L$ is a (special) Lagrangian $m$-fold in $\C^m$ 
invariant under a Lie subgroup $G$ in $\SU(m)\lt\C^m$, 
with moment map $\mu$. Show that $\mu\equiv c$ on $L$ for 
some~$c\in Z(\g^*)$.
\label{l8q1}
\end{question}

\begin{question}Define a smooth map $f:\C^3\ra\R^3$ by
\begin{equation*}
f(z_1,z_2,z_3)=\bigl(\ms{z_1}-\ms{z_3},\ms{z_2}-\ms{z_3},
\Im(z_1z_2z_3)\bigr).
\end{equation*}
For each $a,b,c\in\R^3$, define $N_{a,b,c}=f^{-1}(a,b,c)$.
Then $N_{a,b,c}$ is a real 3-dimensional submanifold
of $\C^3$, which may be singular.
\inext At ${\bf z}=(z_1,z_2,z_3)\in\C^3$, determine
$\d f\vert_{\bf z}:\C^3\ra\R^3$. Find the conditions
on $\bf z$ for $\d f\vert_{\bf z}$ to be surjective.

Now $N_{a,b,c}$ is nonsingular at ${\bf z}\in N_{a,b,c}$ if 
and only if $\d f\vert_{\bf z}$ is surjective. Hence determine
which of the $N_{a,b,c}$ are singular, and find their singular
points.
\inext If $\bf z$ is a nonsingular point of $N_{a,b,c}$, then
$T_{\bf z}N_{a,b,c}=\Ker\d f\vert_{\bf z}$. Determine 
$\Ker\d f\vert_{\bf z}$ in this case, and show that it is a 
special Lagrangian 3-plane in~$\C^3$. 

Hence prove that $N_{a,b,c}$ is a special Lagrangian 3-fold 
wherever it is nonsingular, and that $f:\C^3\ra\R^3$ is
a {\it special Lagrangian fibration}.
\inext Observe that $N_{a,b,c}$ is invariant under the
Lie group $G=\U(1)^2$, acting by
\begin{equation*}
({\rm e}^{i\th_1},{\rm e}^{i\th_2}):
(z_1,z_2,z_3)\mapsto({\rm e}^{i\th_1}z_1,{\rm e}^{i\th_2}z_2,
{\rm e}^{-i\th_1-i\th_2}z_3).
\end{equation*}
How is the form of $f$ related to the ideas of question 
\ref{l8}.\ref{l8q1}? How might $G$-invariance have been 
used to construct the fibration~$f$?
\inext Describe the topology of $N_{a,b,c}$, distinguishing
different cases according to the singularities.
\end{question}

\section{Compact SL $m$-folds in Calabi--Yau $m$-folds}
\label{l9}

In this section we shall discuss {\it compact} special Lagrangian 
submanifolds in Calabi--Yau manifolds. Here are three important 
questions which motivate work in this area.

\begin{itemize}
\item[\bf 1.] Let $N$ be a compact special
Lagrangian $m$-fold in a fixed Calabi--Yau $m$-fold
$(M,J,g,\Om)$. Let ${\cal M}_N$ be the moduli space of 
{\it special Lagrangian deformations} of $N$, that is, 
the connected component of the set of special Lagrangian 
$m$-folds containing $N$. What can we say about ${\cal M}_N$? 
For instance, is it a smooth manifold, and of what dimension?
\smallskip

\item[\bf 2.] Let $\bigl\{(M,J_t,g_t,\Om_t):t\in(-\ep,\ep)\bigr\}$ 
be a smooth 1-parameter family of Calabi--Yau $m$-folds. Suppose $N_0$ 
is an SL $m$-fold in $(M,J_0,g_0,\Om_0)$. Under what conditions can we 
extend $N_0$ to a smooth family of special Lagrangian $m$-folds 
$N_t$ in $(M,J_t,g_t,\Om_t)$ for~$t\in(-\ep,\ep)$?
\smallskip

\item[\bf 3.] In general the moduli space ${\cal M}_N$ in
Question 1 will be noncompact. Can we enlarge ${\cal M}_N$ to
a compact space $\,\,\ov{\!\!\cal M}_N$ by adding a `boundary'
consisting of {\it singular} special Lagrangian $m$-folds?
If so, what is the nature of the singularities that develop?
\end{itemize}

Briefly, these questions concern the {\it deformations} of
special Lagrangian $m$-folds, {\it obstructions} to their
existence, and their {\it singularities} respectively. The
local answers to Questions 1 and 2 are well understood, and 
we shall discuss them in this section. Question 3 is the
subject of~\S\ref{l10}--\S\ref{l11}.

\subsection{SL $m$-folds in Calabi--Yau $m$-folds}
\label{l91}

Here is the definition.

\begin{dfn} Let $(M,J,g,\Om)$ be a Calabi--Yau $m$-fold. Then 
$\Re\Om$ is a {\it calibration} on the Riemannian manifold $(M,g)$. 
An oriented real $m$-dimensional submanifold $N$ in $M$ is called 
a {\it special Lagrangian submanifold (SL\/ $m$-fold)} if it is 
calibrated with respect to~$\Re\Om$.
\label{l9def1}
\end{dfn}

From Proposition \ref{l7prop2} we deduce an {\it alternative
definition} of SL $m$-folds. It is often more useful than
Definition~\ref{l9def1}.

\begin{prop} Let\/ $(M,J,g,\Om)$ be a Calabi--Yau $m$-fold, with
K\"ahler form $\om$, and\/ $L$ a real\/ $m$-dimensional submanifold 
in $M$. Then $N$ admits an orientation making it into an SL\/ 
$m$-fold in $M$ if and only if\/ $\om\vert_N\equiv 0$ 
and\/~$\Im\Om\vert_N\equiv 0$.
\label{l9prop}
\end{prop}

Regard $N$ as an immersed submanifold, with immersion $\iota:N\ra M$. 
Then $[\om\vert_N]$ and $[\Im\Om\vert_N]$ are unchanged under
continuous variations of the immersion $\iota$. Thus, $[\om\vert_N]
=[\Im\Om\vert_N]=0$ is a necessary condition not just for $N$ to
be special Lagrangian, but also for any isotopic submanifold $N'$ in
$M$ to be special Lagrangian. This proves:

\begin{cor} Let\/ $(M,J,g,\Om)$ be a Calabi--Yau $m$-fold, and\/
$N$ a compact real\/ $m$-submanifold in $M$. Then a necessary
condition for $N$ to be isotopic to a special Lagrangian 
submanifold\/ $N'$ in $M$ is that\/ $[\om\vert_N]=0$ in 
$H^2(N,\R)$ and\/ $[\Im\Om\vert_N]=0$ in~$H^m(N,\R)$.
\label{l9cor}
\end{cor}

This gives a simple, necessary topological condition for
an isotopy class of $m$-submanifolds in a Calabi--Yau 
$m$-fold to contain a special Lagrangian submanifold.

\subsection{Deformations of compact special Lagrangian $m$-folds}
\label{l92}

The deformation theory of compact special Lagrangian manifolds
was studied by McLean \cite{McLe}, who proved the following result.

\begin{thm} Let\/ $(M,J,g,\Om)$ be a Calabi--Yau $m$-fold, and\/ 
$N$ a compact special Lagrangian $m$-fold in $M$. Then the moduli 
space ${\cal M}_N$ of special Lagrangian deformations of\/ $N$ is 
a smooth manifold of dimension $b^1(N)$, the first Betti number of\/~$N$.
\label{l9thm1}
\end{thm}
\medskip

\noindent{\it Sketch proof.} Suppose for simplicity that $N$ is an
embedded submanifold. There is a natural orthogonal decomposition 
$TM\vert_N=TN\op\nu$, where $\nu\ra N$ is the {\it normal bundle} 
of $N$ in $M$. As $N$ is Lagrangian, the complex structure 
$J:TM\ra TM$ gives an isomorphism $J:\nu\ra TN$. But the metric $g$ 
gives an isomorphism $TN\cong T^*N$. Composing these two gives an
isomorphism~$\nu\cong T^*N$.

Let $T$ be a small {\it tubular neighbourhood} of $N$ in $M$. Then 
we can identify $T$ with a neighbourhood of the zero section in $\nu$.
Using the isomorphism $\nu\cong T^*N$, we have an identification
between $T$ and a neighbourhood of the zero section in $T^*N$. This
can be chosen to identify the K\"ahler form $\om$ on $T$ with the natural
symplectic structure on $T^*N$. Let $\pi:T\ra N$ be the obvious projection.

Under this identification, submanifolds $N'$ in $T\subset M$ which 
are $C^1$ close to $N$ are identified with the graphs of small smooth 
sections $\al$ of $T^*N$. That is, submanifolds $N'$ of $M$ close to
$N$ are identified with 1-{\it forms} $\al$ on $N$. We need to know: 
which 1-forms $\al$ are identified with {\it special Lagrangian}
submanifolds~$N'$?

Well, $N'$ is special Lagrangian if $\om\vert_{N'}\equiv
\Im\Om\vert_{N'}\equiv 0$. Now $\pi\vert_{N'}:N'\ra N$ is a
diffeomorphism, so we can push $\om\vert_{N'}$ and
$\Im\Om\vert_{N'}$ down to $N$, and regard them as functions 
of $\al$. Calculation~shows~that
\begin{equation*}
\pi_*\bigl(\om\vert_{N'}\bigr)=\d\al
\quad\text{and}\quad
\pi_*\bigl(\Im\Om\vert_{N'}\bigr)=F(\al,\nabla\al),
\end{equation*}
where $F$ is a nonlinear function of its arguments. Thus, the moduli 
space ${\cal M}_N$ is locally isomorphic to the set of small 1-forms 
$\al$ on $N$ such that $\d\al\equiv 0$ and~$F(\al,\nabla\al)\equiv 0$.

Now it turns out that $F$ satisfies $F(\al,\nabla\al)\approx \d(*\al)$ 
when $\al$ is small. Therefore ${\cal M}_N$ is locally approximately 
isomorphic to the vector space of 1-forms $\al$ with $\d\al=\d(*\al)=0$.
But by Hodge theory, this is isomorphic to the de Rham cohomology
group $H^1(N,\R)$, and is a manifold with dimension~$b^1(N)$.

To carry out this last step rigorously requires some technical
machinery: one must work with certain {\it Banach spaces} of 
sections of $T^*N$, $\La^2T^*N$ and $\La^mT^*N$, use {\it elliptic 
regularity results} to prove that the map $\al\mapsto\bigl(\d\al,
F(\al,\nabla\al)\bigr)$ has {\it closed image} in these Banach spaces,
and then use the {\it Implicit Function Theorem for Banach spaces}
to show that the kernel of the map is what we expect.
\relax\unskip\nobreak ~\hfill$\square$\medskip

\subsection{Obstructions to the existence of compact SL $m$-folds}
\label{l93}

Next we address Question 2 above. Let $\bigl\{(M,J_t,g_t,\Om_t):
t\in(-\ep,\ep)\bigr\}$ be a smooth 1-parameter family of Calabi--Yau 
$m$-folds. Suppose $N_0$ is a special Lagrangian $m$-fold of 
$(M,J_0,g_0,\Om_0)$. When can we extend $N_0$ to a smooth family 
of special Lagrangian $m$-folds $N_t$ in $(M,J_t,g_t,\Om_t)$ 
for~$t\in(-\ep,\ep)$?

By Corollary \ref{l9cor}, a necessary condition is that 
$[\om_t\vert_{N_0}]=[\Im\Om_t\vert_{N_0}]=0$ for all $t$.
Our next result shows that locally, this is also a {\it sufficient}\/
condition.

\begin{thm} Let\/ $\bigl\{(M,J_t,g_t,\Om_t):t\in(-\ep,\ep)\bigr\}$
be a smooth\/ $1$-parameter family of Calabi--Yau $m$-folds, with
K\"ahler forms $\om_t$. Let\/ $N_0$ be a compact SL\/ $m$-fold in
$(M,J_0,g_0,\Om_0)$, and suppose that\/ $[\om_t\vert_{N_0}]=0$ in
$H^2(N_0,\R)$ and\/ $[\Im\Om_t\vert_{N_0}]=0$ in $H^m(N_0,\R)$ for
all\/ $t\in(-\ep,\ep)$. Then $N_0$ extends to a smooth\/ $1$-parameter
family $\bigl\{N_t:t\in(-\de,\de)\bigr\}$, where $0<\de\le\ep$ and\/
$N_t$ is a compact SL\/ $m$-fold in~$(M,J_t,g_t,\Om_t)$.
\label{l9thm2}
\end{thm}

This can be proved using similar techniques to Theorem 
\ref{l9thm1}, though McLean did not prove it. Note that the 
condition $[\Im\Om_t\vert_{N_0}]=0$ for all $t$ can be satisfied
by choosing the phases of the $\Om_t$ appropriately, and if the 
image of $H_2(N,\Z)$ in $H_2(M,\R)$ is zero, then the condition 
$[\om\vert_N]=0$ holds automatically. 

Thus, the obstructions $[\om_t\vert_{N_0}]=[\Im\Om_t\vert_{N_0}]=0$ 
in Theorem \ref{l9thm2} are actually fairly mild restrictions, and 
special Lagrangian $m$-folds should be thought of as pretty stable 
under small deformations of the Calabi--Yau structure.
\medskip

\noindent{\bf Remark.} The deformation and obstruction theory of 
compact special Lagrangian $m$-folds are {\it extremely well-behaved}
compared to many other moduli space problems in differential geometry.
In other geometric problems (such as the deformations of complex 
structures on a complex manifold, or pseudo-holomorphic curves in 
an almost complex manifold, or instantons on a Riemannian 4-manifold, 
and so on), the deformation theory often has the following general 
structure.

There are vector bundles $E,F$ over a compact manifold $M$, and an
elliptic operator $P:C^\iy(E)\ra C^\iy(F)$, usually first-order. The 
kernel $\Ker P$ is the set of {\it infinitesimal deformations}, and 
the cokernel $\Coker P$ the set of {\it obstructions}. The actual 
moduli space $\cal M$ is locally the zeros of a nonlinear 
map~$\Psi:\Ker P\ra\Coker P$.

In a {\it generic} case, $\Coker P=0$, and then the moduli space
$\cal M$ is locally isomorphic to $\Ker P$, and so is locally a 
manifold with dimension $\ind(P)$. However, in nongeneric situations 
$\Coker P$ may be nonzero, and then the moduli space $\cal M$ may be 
nonsingular, or have an unexpected dimension.

However, special Lagrangian submanifolds do not follow this pattern. 
Instead, the obstructions are {\it topologically determined}, and the
moduli space is {\it always} smooth, with dimension given by a
topological formula. This should be regarded as a minor mathematical
miracle.

\subsection{Natural coordinates on the moduli space ${\cal M}_N$}
\label{l94}

Let $N$ be a compact SL $m$-fold in a Calabi--Yau $m$-fold 
$(M,J,g,\Om)$. Theorem \ref{l9thm1} shows that the moduli space 
${\cal M}_N$ has dimension $b^1(N)$. By Poincar\'e duality 
$b^1(N)=b^{m-1}(N)$. Thus ${\cal M}_N$ has the same dimension as 
the de Rham cohomology groups $H^1(M,\R)$ and~$H^{m-1}(M,\R)$. 

We shall construct natural local diffeomorphisms $\Phi$ from
${\cal M}_N$ to $H^1(N,\R)$, and $\Psi$ from ${\cal M}_N$ to
$H^{m-1}(N,\R)$. These induce two natural {\it affine structures} 
on ${\cal M}_N$, and can be thought of as two {\it natural 
coordinate systems} on ${\cal M}_N$. The material of this
section can be found in Hitchin~\cite[\S 4]{Hitc}.

Here is how to define $\Phi$ and $\Psi$. Let $U$ be a connected
and simply-connected open neighbourhood of $N$ in ${\cal M}_N$. 
We will construct smooth maps $\Phi:U\ra H^1(N,\R)$ and
$\Psi:U\ra H^{m-1}(N,\R)$ with $\Phi(N)=\Psi(N)=0$, which are 
local diffeomorphisms.

Let $N'\in U$. Then as $U$ is connected, there exists a smooth
path $\ga:[0,1]\ra U$ with $\ga(0)=N$ and $\ga(1)=N'$, and as
$U$ is simply-connected, $\ga$ is unique up to isotopy. Now 
$\ga$ parametrizes a family of submanifolds of $M$ diffeomorphic
to $N$, which we can lift to a smooth map $\Ga:N\t[0,1]\ra M$
with~$\Ga(N\t\{t\})=\ga(t)$.

Consider the 2-form $\Ga^*(\om)$ on $N\t[0,1]$. As each fibre
$\ga(t)$ is Lagrangian, we have $\Ga^*(\om)\vert_{N\t\{t\}}\equiv 0$
for each $t\in[0,1]$. Therefore we may write $\Ga^*(\om)=\al_t\w\d t$,
where $\al_t$ is a closed 1-form on $N$ for $t\in[0,1]$. Define
$\Phi(N')=\bigl[\int_0^1\al_t\,\d t\bigr]\in H^1(N,\R)$. That is, we 
integrate the 1-forms $\al_t$ with respect to $t$ to get a closed 
1-form $\int_0^1\al_t\,\d t$, and then take its cohomology class. 

Similarly, write $\Ga^*(\Im\Om)=\be_t\w\d t$, where $\be_t$ is a 
closed $(m\!-\!1)$-form on $N$ for $t\in[0,1]$, and define
$\Psi(N')=\bigl[\int_0^1\be_t\,\d t\bigr]\in H^{m-1}(N,\R)$.
Then $\Phi$ and $\Psi$ are independent of choices made in the 
construction (exercise). We need to restrict to a simply-connected 
subset $U$ of ${\cal M}_N$ so that $\ga$ is unique up to isotopy. 
Alternatively, one can define $\Phi$ and $\Psi$ on the universal 
cover $\,\,\widetilde{\!\!\cal M}_N$ of~${\cal M}_N$.

\subsection{SL $m$-folds in almost Calabi--Yau manifolds}
\label{l95}

Next we explain a generalization of special Lagrangian geometry
to the class of {\it almost Calabi--Yau manifolds}.

\begin{dfn} Let $m\ge 2$. An {\it almost Calabi--Yau $m$-fold}, or
{\it ACY\/ $m$-fold}\/ for short, is a quadruple $(M,J,g,\Om)$ 
such that $(M,J,g)$ is a compact $m$-dimensional K\"ahler manifold, 
and $\Om$ is a non-vanishing holomorphic $(m,0)$-form on~$M$.
\label{l9def2}
\end{dfn}

The difference between this and Definition \ref{l4def} is that
we do not require $\Om$ and the K\"ahler form $\om$ of $g$ and $\Om$ 
to satisfy equation \eq{omOmeq}, and hence $g$ need not be
Ricci-flat, nor have holonomy $\SU(m)$. Here is the appropriate 
definition of special Lagrangian $m$-folds in ACY $m$-folds.

\begin{dfn} Let $(M,J,g,\Om)$ be an almost Calabi--Yau $m$-fold with
K\"ahler form $\om$, and $N$ a real $m$-dimensional submanifold of $M$.
We call $N$ a {\it special Lagrangian submanifold}, or {\it SL $m$-fold} 
for short, if $\om\vert_N\equiv\Im\Om\vert_N\equiv 0$. It easily
follows that $\Re\Om\vert_N$ is a nonvanishing $m$-form on $N$.
Thus $N$ is orientable, with a unique orientation in which
$\Re\Om\vert_N$ is positive.
\label{l9def3}
\end{dfn}

By Proposition \ref{l9prop}, if $(M,J,g,\Om)$ is Calabi--Yau
rather than almost Calabi--Yau, then $N$ is special Lagrangian 
in the sense of Definition \ref{l4def}. Thus, this is a genuine
extension of the idea of special Lagrangian submanifold. Many 
of the good properties of special Lagrangian submanifolds in 
Calabi--Yau manifolds also apply in almost Calabi--Yau manifolds.
In particular:

\begin{thm} Corollary \ref{l9cor} and Theorems \ref{l9thm1} and 
\ref{l9thm2} also hold in almost Calabi--Yau manifolds rather 
than Calabi--Yau manifolds.
\label{l9thm3}
\end{thm}

This is because the proofs of these results only really depend on
the conditions $\om\vert_N\equiv\Im\Om\vert_N\equiv 0$, and the
pointwise connection \eq{omOmeq} between $\om$ and $\Om$ is not
important. 

Let $(M,J,g,\Om)$ be an ACY $m$-fold, with metric $g$. In general, 
SL $m$-folds in $M$ are neither calibrated nor minimal with respect 
to $g$. However, let $f:M\ra(0,\iy)$ be the unique smooth function 
such that $f^{2m}\om^m/m!=(-1)^{m(m-1)/2}(i/2)^m\Om\w\bar\Om$, and 
define $\ti g$ to be the conformally equivalent metric $f^2g$ on $M$. 
Then $\Re\Om$ is a calibration on the Riemannian manifold $(M,\ti g)$, 
and that SL $m$-folds $N$ in $(M,J,g,\Om)$ are calibrated with respect 
to it, so that they are minimal with respect to~$\ti g$.

The idea of extending special Lagrangian geometry to almost Calabi--Yau 
manifolds appears in the work of Goldstein \cite[\S 3.1]{Gold1}, Bryant 
\cite[\S 1]{Brya}, who uses the term `special K\"ahler' instead of 
`almost Calabi--Yau', and the author~\cite{Joyc13}. 

One important reason for considering SL $m$-folds in almost Calabi--Yau
rather than Calabi--Yau $m$-folds is that they have much stronger
{\it genericness properties}. There are many situations in geometry
in which one uses a genericity assumption to control singular behaviour.

For instance, pseudo-holomorphic curves in an arbitrary almost complex 
manifold may have bad singularities, but the possible singularities
in a generic almost complex manifold are much simpler. In the same 
way, it is reasonable to hope that in a {\it generic} Calabi--Yau 
$m$-fold, compact SL $m$-folds may have better singular behaviour 
than in an arbitrary Calabi--Yau $m$-fold. 

But because Calabi--Yau manifolds come in only finite-dimensional 
families, choosing a generic Calabi--Yau structure is a fairly weak
assumption, and probably will not help very much. However, almost
Calabi--Yau manifolds come in {\it infinite-dimensional}\/ families,
so choosing a generic almost Calabi--Yau structure is a much more 
powerful thing to do, and will probably simplify the singular
behaviour of compact SL $m$-folds considerably. We will return
to this idea in~\S\ref{l10}.

\subsection{Exercises}
\label{l96}

\begin{question}
\label{l96q1}
Show that the maps $\Phi,\Psi$ between special 
Lagrangian moduli space ${\cal M}_N$ and $H^1(N,\R)$, $H^{m-1}(N,\R)$ 
defined in \S\ref{l94} are well-defined and independent of choices. 

Prove also that $\Phi$ and $\Psi$ are {\it local diffeomorphisms}, that 
is, that $\d\Phi\vert_{N'}$ and $\d\Psi\vert_{N'}$ are isomorphisms
between $T_{N'}{\cal M}_N$ and $H^1(N,\R)$, $H^{m-1}(N,\R)$ for
each~$N'\in U$.
\end{question}

\begin{question}Putting together the maps $\Phi,\Psi$ of Question
\ref{l9}.\ref{l96q1} gives a map $\Phi\t\Psi:U\ra H^1(N,\R)\t 
H^{m-1}(N,\R)$. Now $H^1(N,\R)$ and $H^{m-1}(N,\R)$ are dual by 
Poincar\'e duality, so $H^1(N,\R)\t H^{m-1}(N,\R)$ has a natural 
{\it symplectic structure}. Show that the image of $U$ is a 
{\it Lagrangian submanifold} in $H^1(N,\R)\t H^{m-1}(N,\R)$. 
\smallskip

{\bf Hint:} From the proof of McLean's theorem in \S\ref{l92}, the 
tangent space $T_N{\cal M}_N$ is isomorphic to the vector space of 
1-forms $\al$ with $\d\al=\d(*\al)=0$. Then $\d\Phi\vert_N:T_N{\cal M}
\ra H^1(M,\R)$ takes $\al\mapsto[\al]$, and $\d\Psi\vert_N:T_N{\cal M}
\ra H^{m-1}(M,\R)$ takes $\al\mapsto[*\al]$. Use the fact that for 
1-forms $\al,\be$ on an oriented Riemannian manifold we 
have~$\al\w(*\be)=\be\w(*\al)$.
\end{question}

\section{Singularities of special Lagrangian $m$-folds}
\label{l10}

Now we move on to Question 3 of \S\ref{l9}, and discuss the
{\it singularities} of special Lagrangian $m$-folds. We can
divide it into two sub-questions:
\begin{itemize}
\item[\bf 3(a)] What kinds of singularities are possible in
singular special Lagrangian $m$-folds, and what do
they look like?
\item[\bf 3(b)] How can singular SL $m$-folds arise as limits
of nonsingular SL $m$-folds, and what does the limiting
behaviour look like near the singularities?
\end{itemize}

The basic premise of the author's approach to special Lagrangian
singularities is that singularities of SL $m$-folds in Calabi--Yau
$m$-folds should look locally like singularities of SL $m$-folds
in $\C^m$, to the first few orders of approximation. That is, if
$M$ is a Calabi--Yau $m$-fold and $N$ an SL $m$-fold in $M$ with
a singular point at $x\in M$, then near $x$, $M$ resembles
$\C^m=T_xM$, and $N$ resembles an SL $m$-fold $L$ in $\C^m$ with a 
singular point at 0. We call $L$ a {\it local model} for $N$ near~$x$.

Therefore, to understand singularities of SL $m$-folds in Calabi--Yau
manifolds, we begin by studying singularities of SL $m$-folds in $\C^m$,
first by constructing as many examples as we can, and then by aiming for
some kind of rough classification of the most common kinds of special 
Lagrangian singularities, at least in low dimensions such as~$m=3$.

\subsection{Cones, and asymptotically conical SL $m$-folds}
\label{l101}

In Examples \ref{l8ex1}--\ref{l8ex3}, the singular SL 3-folds we 
constructed were cones. Here a closed SL $m$-fold $C$ in $\C^m$ is 
called a {\it cone} if $C=tC$ for all $t>0$, where $tC=\{tx:x\in C\}$. 
Note that 0 is always a singular point of $C$, unless $C$ is a
special Lagrangian plane $\R^m$ in $\C^m$. The simplest kind of 
SL cones (from the point of view of singular behaviour) are those in 
which 0 is the only singular point. Then $\Si=C\cap{\cal S}^{2m-1}$
is a nonsingular, compact, minimal, Legendrian $(m-1)$-submanifold
in the unit sphere ${\cal S}^{2m-1}$ in~$\C^m$.

In one sense, {\it all}\/ singularities of SL $m$-folds are 
modelled on special Lagrangian cones, to highest order. It 
follows from \cite[\S II.5]{HaLa} that if $M$ is a Calabi--Yau
$m$-fold and $N$ an SL $m$-fold in $M$ with a singular point
at $x$, then $N$ has a {\it tangent cone} at $x$ in the sense
of Geometric Measure Theory, which is a special Lagrangian 
cone $C$ in~$\C^m=T_xM$.

If $C$ has multiplicity one and 0 is its only singular point,
then $N$ really does look like $C$ near $x$. However, things
become more complicated if the singularities of $C$ are not
isolated (for instance, if $C$ is the union of two SL 3-planes
$\R^3$ in $\C^3$ intersecting in a line) or if the multiplicity
of $C$ is greater than 1 (this happens in Example \ref{l8ex4},
as the tangent cone is a double $\R^3$\,). In these cases, the
tangent cone captures only the simplest part of the singular
behaviour, and we have to work harder to find out what is going on.

Now suppose for simplicity that we are interested in SL $m$-folds
with singularities modelled on a multiplicity 1 SL cone $C$
in $\C^m$ with an isolated singular point. To answer question
3(b) above, and understand how such singularities arise as limits 
of nonsingular SL $m$-folds, we need to study {\it asymptotically
conical (AC)}\/ SL $m$-folds $L$ in $\C^m$ asymptotic to~$C$. 

We shall be interested in two classes of asymptotically conical 
SL $m$-folds $L$, {\it weakly AC}, which converge to $C$ like
$o(r)$, and {\it strongly AC}, which converge to $C$ like
$O(r^{-1})$. Here is a more precise definition.

\begin{dfn} Let $C$ be a special Lagrangian cone in $\C^m$ with 
isolated singularity at 0, and let $\Si=C\cap{\cal S}^{2m-1}$,
so that $\Si$ is a compact, nonsingular $(m-1)$-manifold. Define
the {\it number of ends of\/ $C$ at infinity} to be the number of
connected components of $\Si$. Let $h$ be the metric on $\Si$
induced by the metric $g$ on $\C^m$, and $r$ the radius function
on $\C^m$. Define $\iota:\Si\t(0,\iy)\ra\C^m$ by $\iota(\si,r)=r\si$.
Then the image of $\iota$ is $C\sm\{0\}$, and $\iota^*(g)=r^2h+\d r^2$
is the cone metric on~$C\sm\{0\}$.

Let $L$ be a closed, nonsingular SL $m$-fold in $\C^m$. We call $L$
{\it weakly asymptotically conical (weakly AC)} with cone $C$ if
there exists a compact subset $K\subset L$ and a diffeomorphism
$\phi:\Si\t(R,\iy)\ra L\sm K$ for some $R>0$, such that
$\md{\phi-\iota}=o(r)$ and $\bmd{\nabla^k(\phi-\iota)}=o(r^{1-k})$
as $r\ra\iy$ for $k=1,2,\ldots$, where $\nabla$ is the Levi-Civita
connection of the cone metric $\iota^*(g)$, and $\md{\,.\,}$ is
computed using~$\iota^*(g)$.

Similarly, we call $L$ {\it strongly asymptotically conical 
(strongly AC)} with cone $C$ if $\md{\phi-\iota}=O(r^{-1})$ 
and $\bmd{\nabla^k(\phi-\iota)}=O(r^{-1-k})$ as $r\ra\iy$ for
$k=1,2,\ldots$, using the same notation.
\label{l10def1}
\end{dfn}

These two asymptotic conditions are useful for different purposes.
If $L$ is a weakly AC SL $m$-fold then $tL$ converges to $C$ as 
$t\ra 0_+$. Thus, weakly AC SL $m$-folds provide models for how 
singularities modelled on cones $C$ can arise as limits of 
nonsingular SL $m$-folds. The weakly AC condition is in practice 
the weakest asymptotic condition which ensures this; if the $o(r)$ 
condition were made any weaker then the asymptotic cone $C$ at 
infinity might not be unique or well-defined.

On the other hand, explicit constructions tend to produce strongly
AC SL $m$-folds, and they appear to be the easiest class to prove
results about. For example, one can show:

\begin{prop} Suppose $C$ is an SL cone in $\C^m$ with an isolated
singular point at\/ $0$, invariant under a connected Lie subgroup
$G$ of\/ $\SU(m)$. Then any strongly AC SL $m$-fold $L$ in $\C^m$
with cone $C$ is also $G$-invariant.
\label{l10prop1}
\end{prop}

This should be a help in classifying strongly AC SL $m$-folds with
cones with a lot of symmetry. For instance, using $\U(1)^2$ symmetry 
one can show that the only strongly AC SL 3-folds in $\C^3$ with cone 
$L_0$ from Example \ref{l8ex1} are the SL 3-folds $L_t$ from Example 
\ref{l8ex1} for $t>0$, and two other families obtained from the $L_t$
by cyclic permutations of coordinates~$(z_1,z_2,z_3)$.

\subsection{Moduli spaces of AC SL $m$-folds}
\label{l102}

Next we discuss moduli space problems for AC SL $m$-folds. I shall state
our problems as conjectures, because the proofs are not yet complete.
One should be able to prove analogues of Theorem \ref{l9thm1} for AC SL
$m$-folds. Here is the appropriate result for strongly AC SL $m$-folds.

\begin{conj} Let\/ $L$ be a strongly AC SL $m$-fold in $\C^m$, 
with cone $C$, and let\/ $k$ be the number of ends of\/ $C$
at infinity. Then the moduli space ${\cal M}^s_{\sst L}$ of 
strongly AC SL $m$-folds in $\C^m$ with cone $C$ is near $L$ 
a smooth manifold of dimension~$b^1(L)+k-1$.
\label{l10conj1}
\end{conj}

Before generalizing this to the weak case, here is a definition.

\begin{dfn} Let $C$ be a special Lagrangian cone in $\C^m$ with 
an isolated singularity at 0, and let $\Si=C\cap{\cal S}^{2m-1}$.
Regard $\Si$ as a compact Riemannian manifold, with metric induced 
from the round metric on ${\cal S}^{2m-1}$. Let $\De=\d^*\d$ be the
Laplacian on functions on $\Si$. Define the {\it Legendrian index}
$\lind(C)$ to be the number of eigenvalues of $\De$ in $(0,2m)$,
counted with multiplicity.
\label{l10def2}
\end{dfn}

We call this the {\it Legendrian index} since it is the
index of the area functional at $\Si$ under Legendrian
variations of the submanifold $\Si$ in ${\cal S}^{2m-1}$.
It is not difficult to show that the restriction of any real
linear function on $\C^m$ to $\Si$ is an eigenfunction of $\De$
with eigenvalue $m-1$. These contribute $m$ to $\lind(C)$ for
each connected component of $\Si$ which is a round unit sphere
${\cal S}^{m-1}$, and $2m$ to $\lind(C)$ for each other
connected component of $\Si$. This gives a useful lower bound 
for $\lind(C)$. In particular,~$\lind(C)\ge 2m$. 

\begin{conj} Let\/ $L$ be a weakly AC SL\/ $m$-fold in $\C^m$, 
with cone $C$, and let\/ $k$ be the number of ends of\/ $C$
at infinity. Then the moduli space ${\cal M}^w_{\sst L}$ of 
weakly AC SL\/ $m$-folds in $\C^m$ with cone $C$ is near $L$ 
a smooth manifold of dimension~$b^1(L)+k-1+\lind(C)$.
\label{l10conj2}
\end{conj}

Here are some remarks on these conjectures:

\begin{itemize}
\item The author's student, Stephen Marshall, is working on proofs
of Conjectures \ref{l10conj1} and \ref{l10conj2}, which we hope to
be able to publish soon. Some related results have recently been
proved by Tommaso Pacini,~\cite[\S 3]{Paci}.

\item If $L$ is a weakly AC SL $m$-fold in $\C^m$, then any translation 
of $L$ is also weakly AC, with the same cone. Since $C$ has an 
isolated singularity by assumption, it cannot have translation 
symmetries. Hence $L$ also has no translation symmetries, so the 
translations of $L$ are all distinct, and ${\cal M}^w_{\sst L}$ 
has dimension at least $2m$. The inequality $\lind(C)\ge 2m$ 
above ensures this.
\item The dimension of ${\cal M}^s_{\sst L}$ in Conjecture
\ref{l10conj1} is purely topological, as in Theorem \ref{l9thm1},
which is another indication that strongly AC is in many ways
the nicest asymptotic condition to work with. But the dimension 
of ${\cal M}^w_{\sst L}$ in Conjecture \ref{l10conj2} has an
analytic component, the eigenvalue count in~$\lind(C)$.
\item It is an interesting question whether moduli spaces of 
weakly AC SL $m$-folds always contain a strongly AC SL $m$-fold.
\end{itemize}

\subsection{SL singularities in generic almost Calabi--Yau $m$-folds}
\label{l103}

We move on to discuss the singular behaviour of compact SL $m$-folds
in Calabi--Yau $m$-folds. For simplicity we shall restrict our 
attention to a class of SL cones with no nontrivial deformations.

\begin{dfn} Let $C$ be a special Lagrangian cone in $\C^m$ with 
an isolated singularity at 0 and $k$ ends at infinity, and set
$\Si=C\cap{\cal S}^{2m-1}$. Let $\Si_1,\ldots,\Si_k$ be the
connected components of $\Si$. Regard each $\Si_j$ as a compact
Riemannian manifold, and let $\De_j=\d^*\d$ be the Laplacian on
functions on $\Si_j$. Let $G_j$ be the Lie subgroup of $\SU(m)$
preserving $\Si_j$, and $V_j$ the eigenspace of $\De_j$ with
eigenvalue $2m$. We call the SL cone $C$ {\it rigid} if
$\dim V_j=\dim\SU(m)-\dim G_j$ for each~$j=1,\ldots,k$.
\label{l10def3}
\end{dfn}

Here is how to understand this definition. We can regard $C$ as
the union of one-ended SL cones $C_1,\ldots,C_k$ intersecting at 0,
where $C_j\sm\{0\}$ is naturally identified with $\Si_j\t(0,\iy)$.
The cone metric on $C_j$ is $g_j=r^2h_j+\d r^2$, where $h_j$ is the
metric on $\Si_j$. Suppose $f_j$ is an eigenfunction of $\De_j$ on
$\Si_j$ with eigenvalue $2m$. Then $r^2f_j$ is {\it harmonic} on~$C_j$.

Hence, $\d(r^2f_j)$ is a closed, coclosed 1-form on $C_j$ which is
linear in $r$. By the Principle in \S\ref{l73}, the basis of the proof
of Theorem \ref{l9thm1}, such 1-forms correspond to small deformations 
of $C_j$ as an SL cone in $\C^m$. Therefore, we can interpret $V_j$ as
the space of {\it infinitesimal deformations} of $C_j$ as a special
Lagrangian cone. 

Clearly, one way to deform $C_j$ as a special Lagrangian cone is 
to apply elements of $\SU(m)$. This gives a family $\SU(m)/G_j$ of
deformations of $C_j$, with dimension $\dim\SU(m)-\dim G_j$, so that
$\dim V_j\ge\dim\SU(m)-\dim G_j$. (The corresponding functions in $V_j$
are moment maps of $\su(m)$ vector fields.) We call $C$ {\it rigid} if
equality holds for all $j$, that is, if all infinitesimal deformations
of $C$ come from applying motions in $\su(m)$ to the component
cones~$C_1,\ldots,C_k$.

Not all SL cones in $\C^m$ are rigid. One can show using integrable
systems that there exist families of SL $T^2$-cones $C$ in $\C^3$
up to $\SU(3)$ equivalence, of arbitrarily high dimension. If the
dimension of the family is greater than $\dim\SU(3)-\dim G_1$,
where $G_1$ is the Lie subgroup of $\SU(3)$ preserving $C$, then
$C$ is not rigid.

Now we can give a first approximation to the kinds of results 
the author expects to hold for singular SL $m$-folds in (almost) 
Calabi--Yau $m$-folds.

\begin{conj} Let\/ $C$ be a rigid special Lagrangian cone in $\C^m$
with an isolated singularity at\/ $0$ and\/ $k$ ends at infinity,
and\/ $L$ a weakly AC SL\/ $m$-fold in $\C^m$ with cone $C$. Let\/
$(M,J,g,\Om)$ be a generic almost Calabi--Yau $m$-fold, and\/ $\cal M$
a connected moduli space of compact nonsingular SL\/ $m$-folds $N$ in~$M$. 

Suppose that at the boundary of\/ $\cal M$ there is a moduli space 
${\cal M}_C$ of compact, singular SL\/ $m$-folds with one isolated 
singular point modelled on the cone $C$, which arise as limits of
SL $m$-folds in $\cal M$ by collapsing weakly AC SL\/ $m$-folds
with the topology of\/ $L$. Then
\e
\dim{\cal M}=\dim{\cal M}_C+b^1(L)+k-1+\lind(C)-2m.
\label{l10eq2}
\e
\label{l10conj3}
\end{conj}

Here are some remarks on the conjecture:
\begin{itemize}
\item I have an outline proof of this conjecture which works
when $m<6$. The analytic difficulties increase with dimension; 
I am not sure whether the conjecture holds in high dimensions.
\item Similar results should hold for non-rigid singularities,
but the dimension formulae will be more complicated.
\item Closely associated to this result is an analogue of 
Theorem \ref{l9thm1} for SL $m$-folds with isolated conical 
singularities of a given kind, under an appropriate genericity
assumption.
\item Here is one way to arrive at equation \eq{l10eq2}. Assuming
Conjecture \ref{l10conj2}, the moduli space ${\cal M}^w_{\sst L}$ of 
weakly AC SL $m$-folds containing $L$ has dimension $b^1(L)+k-1+\lind(C)$.
Now translations in $\C^m$ act freely on ${\cal M}^w_{\sst L}$, so the
family of weakly AC SL $m$-folds up to translations has
dimension~$b^1(L)+k-1+\lind(C)-2m$.

The idea is that each singular SL $m$-fold $N_0$ in ${\cal M}_C$ 
can be `resolved' to give a nonsingular SL $m$-fold $N$ by gluing
in any `sufficiently small' weakly AC SL $m$-fold $L'$, up to
translation. Thus, desingularizing should add $b^1(L)+k-1+\lind(C)-2m$
degrees of freedom, which is how we get equation~\eq{l10eq2}.
\item However, \eq{l10eq2} may not give the right answer in every 
case. One can imagine situations in which there are {\it cohomological 
obstructions} to gluing in AC SL $m$-folds $L'$. For instance, it might 
be necessary that the symplectic area of a disc in $\C^m$ with boundary 
in $L'$ be zero to make $N=L'\# N_0$ Lagrangian. These could reduce
the number of degrees of freedom in desingularizing $N_0$, and then
\eq{l10eq2} would require correction. An example of this is considered
in~\cite[\S 4.4]{Joyc1}.
\end{itemize}

Now we can introduce the final important idea in this section.
Suppose we have a suitably generic (almost) Calabi--Yau $m$-fold 
$M$ and a compact, singular SL $m$-fold $N_0$ in $M$, which is the 
limit of a family of compact nonsingular SL $m$-folds $N$ in~$M$. 

We (loosely) define the {\it index} of the singularities of $N_0$ 
to be the codimension of the family of singular SL $m$-folds with 
singularities like those of $N_0$ in the family of nonsingular SL 
$m$-folds $N$. Thus, in the situation of Conjecture \ref{l10conj3}, 
the index of the singularities is~$b^1(L)+k-1+\lind(C)-2m$.

More generally, one can work not just with a fixed generic almost
Calabi--Yau $m$-fold, but with a {\it generic family} of almost 
Calabi--Yau $m$-folds. So, for instance, if we have a generic 
$k$-dimensional family of almost Calabi--Yau $m$-folds $M$, and in 
each $M$ we have an $l$-dimensional family of SL $m$-folds, then 
in the total $(k\!+\!l)$-dimensional family of SL $m$-folds we are 
guaranteed to meet singularities of index at most~$k\!+\!l$.

Now singularities with {\it small index} are the most commonly
occurring, and so arguably the most interesting kinds of singularity.
Also, for various problems it will only be necessary to know about
singularities with index up to a certain value.

For example, in \cite{Joyc1} the author proposed to define an
invariant of almost Calabi--Yau 3-folds by counting special 
Lagrangian homology 3-spheres (which occur in 0-dimensional 
moduli spaces) in a given homology class, with a certain 
topological weight. This invariant will only be interesting 
if it is essentially conserved under deformations of the 
underlying almost Calabi--Yau 3-fold. During such a deformation, 
nonsingular SL 3-folds can develop singularities and disappear, 
or new ones appear, which might change the invariant.

To prove the invariant is conserved, we need to show that it is 
unchanged along generic 1-parameter families of Calabi--Yau 3-folds.
The only kinds of singularities of SL homology 3-spheres that arise 
in such families will have index 1. Thus, to resolve the conjectures 
in \cite{Joyc1}, we only have to know about index 1 singularities
of SL 3-folds in almost Calabi--Yau 3-folds.

Another problem in which the index of singularities will be important is 
the {\it SYZ Conjecture}, to be discussed in \S\ref{l11}. This has to 
do with dual 3-dimensional families ${\cal M}_X$, ${\cal M}_{\hat X}$ 
of SL 3-tori in (almost) Calabi--Yau 3-folds $X,\hat X$. If $X,\hat X$ 
are generic then the only kinds of singularities that can occur at the 
boundaries of ${\cal M}_X,{\cal M}_{\hat X}$ are of index 1, 2 or 3. So, 
to study the SYZ Conjecture in the generic case, we only have to know 
about singularities of SL 3-folds with index 1, 2, 3 (and possibly~4).

It would be an interesting and useful project to find examples
of, and eventually to classify, special Lagrangian singularities
with small index, at least in dimension 3. For instance, consider
rigid SL cones $C$ in $\C^3$ as in Conjecture \ref{l10conj3}, of
index 1. Then $b^1(L)+k-1+\lind(C)-6=1$, and $\lind(C)\ge 6$, so
$b^1(L)+k\le 2$. But $k\ge 1$ and $b^1(L)\ge\ha b^1(\Si)$, so
$b^1(\Si)\le 2$. As $\Si$ is oriented, one can show that either
$k=1$, $\lind(C)=6$ and $\Si$ is a torus $T^2$, or $k=2$,
$\lind(C)=6$, $\Si$ is 2 copies of ${\cal S}^2$, and $C$ is the
union of two SL 3-planes in $\C^3$ intersecting only at~0.

The eigenvalue count $\lind(C)$ implies an upper bound for the area 
of $\Si$. Hopefully, one can then use integrable systems results as
in \S\ref{l84} to pin down the possibilities for $C$. The author 
guesses that the $T^2$-cone $L_0$ of Example \ref{l8ex1}, and perhaps
also the $T^2$-cone $L^{1,2}_0$ of Example \ref{l8ex3}, are the only 
examples of index 1 SL $T^2$-cones in $\C^3$ up to $\SU(3)$ isomorphisms.

\subsection{Exercises}
\label{l104}

\begin{question} Prove Proposition~\ref{l10prop1}.
\end{question}

\section{The SYZ Conjecture, and SL fibrations}
\label{l11}

{\it Mirror Symmetry} is a mysterious relationship between pairs 
of Calabi--Yau 3-folds $X,\hat X$, arising from a branch of 
physics known as {\it String Theory}, and leading to some very 
strange and exciting conjectures about Calabi--Yau 3-folds, 
many of which have been proved in special cases.

The {\it SYZ Conjecture} is an attempt to explain Mirror Symmetry
in terms of dual `fibrations' $f:X\ra B$ and $\hat f:\hat X\ra B$
of $X,\hat X$ by special Lagrangian 3-folds, including singular
fibres. We give brief introductions to String Theory, Mirror 
Symmetry, and the SYZ Conjecture, and then a short survey of
the state of mathematical research into the SYZ Conjecture, 
biased in favour of the author's own interests.

\subsection{String Theory and Mirror Symmetry}
\label{l111}

String Theory is a branch of high-energy theoretical physics in 
which particles are modelled not as points but as 1-dimensional 
objects -- `strings' -- propagating in some background space-time
$M$. String theorists aim to construct a {\it quantum theory} of 
the string's motion. The process of quantization is extremely 
complicated, and fraught with mathematical difficulties that 
are as yet still poorly understood.

The most popular version of String Theory requires the universe
to be 10-dimensional for this quantization process to work.
Therefore, String Theorists suppose that the space we live in 
looks locally like $M=\R^4\t X$, where $\R^4$ is Minkowski space, 
and $X$ is a compact Riemannian 6-manifold with radius of order 
$10^{-33}$cm, the Planck length. Since the Planck length is so 
small, space then appears to macroscopic observers to be 4-dimensional.

Because of supersymmetry, $X$ has to be a {\it Calabi--Yau $3$-fold}. 
Therefore String Theorists are very interested in Calabi--Yau 3-folds. 
They believe that each Calabi--Yau 3-fold $X$ has a quantization, which 
is a {\it Super Conformal Field Theory} (SCFT), a complicated
mathematical object. Invariants of $X$ such as the Dolbeault 
groups $H^{p,q}(X)$ and the number of holomorphic curves in $X$ 
translate to properties of the SCFT.

However, two entirely different Calabi--Yau 3-folds $X$ and $\hat X$ 
may have the {\it same} SCFT. In this case, there are powerful 
relationships between the invariants of $X$ and of $\hat X$ that
translate to properties of the SCFT. This is the idea behind
{\it Mirror Symmetry} of Calabi--Yau 3-folds.

It turns out that there is a very simple automorphism of the
structure of a SCFT --- changing the sign of a $\U(1)$-action 
--- which does {\it not\/} correspond to a classical automorphism
of Calabi--Yau 3-folds. We say that $X$ and $\hat X$ are {\it mirror}
Calabi--Yau 3-folds if their SCFT's are related by this automorphism.
Then one can argue using String Theory that
\begin{equation*}
H^{1,1}(X)\cong H^{2,1}(\hat X) \quad\text{and}\quad
H^{2,1}(X)\cong H^{1,1}(\hat X).
\end{equation*}
Effectively, the mirror transform exchanges even- and
odd-dimensional cohomology. This is a very surprising result!

More involved String Theory arguments show that, in effect,
the Mirror Transform exchanges things related to the complex
structure of $X$ with things related to the symplectic structure
of $\hat X$, and vice versa. Also, a generating function for the 
number of holomorphic rational curves in $X$ is exchanged
with a simple invariant to do with variation of complex 
structure on $\hat X$, and so on.

Because the quantization process is poorly understood and
not at all rigorous --- it involves non-convergent path-integrals 
over horrible infinite-dimensional spaces --- String Theory
generates only conjectures about Mirror Symmetry, not proofs.
However, many of these conjectures have been verified in
particular cases.

\subsection{Mathematical interpretations of Mirror Symmetry}
\label{l112}

In the beginning (the 1980's), Mirror Symmetry seemed mathematically 
completely mysterious. But there are now two complementary conjectural
theories, due to Kontsevich and Strominger--Yau--Zaslow, which explain 
Mirror Symmetry in a fairly mathematical way. Probably both are true, 
at some level.

The first proposal was due to Kontsevich \cite{Kont} in 1994. This 
says that for mirror Calabi--Yau 3-folds $X$ and $\hat X$, the derived 
category $D^b(X)$ of coherent sheaves on $X$ is equivalent to the derived 
category $D^b({\rm Fuk}(\hat X))$ of the Fukaya category of $\hat X$, 
and vice versa. Basically, $D^b(X)$ has to do with $X$ as a complex
manifold, and $D^b({\rm Fuk}(\hat X))$ with $\hat X$ as a symplectic 
manifold, and its Lagrangian submanifolds. We shall not discuss this here.

The second proposal, due to Strominger, Yau and Zaslow \cite{SYZ} in 1996, 
is known as the {\it SYZ Conjecture}. Here is an attempt to state it.
\medskip

\noindent{\bf The SYZ Conjecture} {\it Suppose $X$ and\/ $\hat X$ are 
mirror Calabi--Yau $3$-folds. Then (under some additional conditions) 
there should exist a compact topological\/ $3$-manifold\/ $B$ and 
surjective, continuous maps $f:X\ra B$ and\/ $\hat f:\hat X\ra B$, 
such that
\begin{itemize}
\item[{\rm(i)}] There exists a dense open set\/ $B_0\subset B$, such 
that for each\/ $b\in B_0$, the fibres $f^{-1}(b)$ and\/ $\hat f^{-1}(b)$
are nonsingular special Lagrangian $3$-tori $T^3$ in $X$ and\/ $\hat X$.
Furthermore, $f^{-1}(b)$ and\/ $\hat f^{-1}(b)$ are in some sense
dual to one another.
\item[{\rm(ii)}] For each\/ $b\in \De=B\sm B_0$, the fibres $f^{-1}(b)$ 
and\/ $\hat f^{-1}(b)$ are expected to be singular special Lagrangian
$3$-folds in $X$ and\/~$\hat X$.
\end{itemize}}
\medskip

We call $f$ and $\hat f$ {\it special Lagrangian fibrations}, and
the set of singular fibres $\De$ is called the {\it discriminant}.
In part (i), the nonsingular fibres of $f$ and $\hat f$ are supposed
to be {\it dual tori}. What does this mean?

On the topological level, we can define duality between two tori 
$T,\hat T$ to be a choice of isomorphism $H^1(T,\Z)\cong H_1(\hat T,\Z)$. 
We can also define duality between tori equipped with flat Riemannian
metrics. Write $T=V/\La$, where $V$ is a Euclidean vector space
and $\La$ a {\it lattice} in $V$. Then the dual torus $\hat T$ is 
defined to be $V^*/\La^*$, where $V^*$ is the dual vector space and 
$\La^*$ the dual lattice. However, there is no notion of duality
between non-flat metrics on dual tori.

Strominger, Yau and Zaslow argue only that their conjecture holds
when $X,\hat X$ are close to the `large complex structure limit'.
In this case, the diameters of the fibres $f^{-1}(b),\hat f^{-1}(b)$
are expected to be small compared to the diameter of the base space
$B$, and away from singularities of $f,\hat f$, the metrics on the
nonsingular fibres are expected to be approximately flat. 

So, part (i) of the SYZ Conjecture says that for $b\in B\sm B_0$, 
$f^{-1}(b)$ is approximately a flat Riemannian 3-torus, and $\hat f^{-1}(b)$ 
is approximately the dual flat Riemannian torus. Really, the SYZ Conjecture
makes most sense as a statement about the limiting behaviour of 
{\it families} of mirror Calabi--Yau 3-folds $X_t$, $\hat X_t$ 
which approach the `large complex structure limit' as~$t\ra 0$.

\subsection{The symplectic topological approach to SYZ}
\label{l113}

The most successful approach to the SYZ Conjecture so far could 
be described as {\it symplectic topological}. In this approach, 
we mostly forget about complex structures, and treat $X,\hat X$ 
just as {\it symplectic manifolds}. We mostly forget about the 
`special' condition, and treat $f,\hat f$ just as {\it Lagrangian 
fibrations}. We also impose the condition that $B$ is a {\it smooth\/} 
3-manifold and $f:X\ra B$ and $\hat f:\hat X\ra B$ are {\it smooth 
maps}. (It is not clear that $f,\hat f$ can in fact be smooth at 
every point, though).

Under these simplifying assumptions, Gross \cite{Gros1,Gros2,Gros3,Gros4}, 
Ruan \cite{Ruan1,Ruan2}, and others have built up a beautiful, detailed 
picture of how dual SYZ fibrations work at the global topological level, 
in particular for examples such as the quintic and its mirror, and for 
Calabi--Yau 3-folds constructed as hypersurfaces in toric 4-folds, using 
combinatorial data.

\subsection{Local geometric approach, and SL singularities}
\label{l114}

There is also another approach to the SYZ Conjecture, begun by the
author in \cite{Joyc11,Joyc13}, and making use of the ideas and
philosophy set out in \S\ref{l10}. We could describe it as a
{\it local geometric} approach.

In it we try to take the special Lagrangian condition seriously 
from the outset, and our focus is on the local behaviour of special 
Lagrangian submanifolds, and especially their singularities, rather 
than on global topological questions. Also, we are interested in what 
fibrations of {\it generic} (almost) Calabi--Yau 3-folds might look like.

One of the first-fruits of this approach has been the understanding that
for {\it generic} (almost) Calabi--Yau 3-folds $X$, special Lagrangian
fibrations $f:X\ra B$ will not be smooth maps, but only piecewise
smooth. Furthermore, their behaviour at the singular set is rather 
different to the smooth Lagrangian fibrations discussed in~\S\ref{l113}. 

For smooth special Lagrangian fibrations $f:X\ra B$, the discriminant
$\De$ is of codimension 2 in $B$, and the typical singular fibre is 
singular along an ${\cal S}^1$. But in a generic special Lagrangian 
fibration $f:X\ra B$ the discriminant $\De$ is of codimension 1 in $B$, 
and the typical singular fibre is singular at finitely many points.

One can also show that if $X,\hat X$ are a mirror pair of generic
(almost) Calabi--Yau 3-folds and $f:X\ra B$ and $\hat f:\hat X\ra B$ are
dual special Lagrangian fibrations, then in general the discriminants
$\De$ of $f$ and $\hat\De$ of $\hat f$ cannot coincide in $B$, because
they have different topological properties in the neighbourhood of a
certain kind of codimension 3 singular fibre. 

This contradicts part (ii) of the SYZ Conjecture, as we have stated it
in \S\ref{l112}. In the author's view, these calculations support the 
idea that the SYZ Conjecture in its present form should be viewed 
primarily as a limiting statement, about what happens at the `large 
complex structure limit', rather than as simply being about pairs of 
Calabi--Yau 3-folds. A similar conclusion is reached by Mark Gross
in~\cite[\S 5]{Gros4}.

\subsection{$\U(1)$-invariant SL fibrations in $\C^3$}
\label{l115}

We finish by describing work of the author in \cite[\S 8]{Joyc11}
and \cite{Joyc13}, which aims to describe what the singularities
of SL fibrations of {\it generic} (almost) Calabi--Yau 3-folds
look like, providing they exist. This proceeds by first studying
SL fibrations of subsets of $\C^3$ invariant under the
$\U(1)$-action \eq{l8eq2}, using the ideas of \S\ref{l85}.
For a brief survey of the main results, see~\cite{Joyc12}.

Then we argue (without a complete proof, as yet) that the kinds
of singularities we see in codimension 1 and 2 in generic
$\U(1)$-invariant SL fibrations in $\C^3$, also occur in
codimension 1 and 2 in SL fibrations of generic (almost)
Calabi--Yau 3-folds, without any assumption of $\U(1)$-invariance.

Following \cite[Def.~8.1]{Joyc11}, we use the results of
\S\ref{l85} to construct a family of SL 3-folds $N_{\bs\al}$ in
$\C^3$, depending on boundary data~$\Phi({\bs\al})$.

\begin{dfn} Let $S$ be a strictly convex domain in $\R^2$ invariant
under $(x,y)\mapsto(x,-y)$, let $U$ be an open set in $\R^3$, and
$\al\in(0,1)$. Suppose $\Phi:U\ra C^{3,\al}(\pd S)$ is a continuous
map such that if $(a,b,c)\ne(a,b',c')$ in $U$ then $\Phi(a,b,c)-
\Phi(a,b',c')$ has exactly one local maximum and one local minimum
in~$\pd S$.

For ${\bs\al}=(a,b,c)\in U$, let $f_{\bs\al}\in C^{3,\al}(S)$
or $C^1(S)$ be the unique (weak) solution of \eq{l8eq6} with
$f_{\bs\al}\vert_{\pd S}=\Phi({\bs\al})$, which exists by
Theorem \ref{l8thm1}. Define $u_{\bs\al}=\frac{\pd f_{\bs\al}}{\pd y}$
and $v_{\bs\al}=\frac{\pd f_{\bs\al}}{\pd x}$. Then $(u_{\bs\al},
v_{\bs\al})$ is a solution of \eq{l8eq5} in $C^{2,\al}(S)$ if
$a\ne 0$, and a weak solution of \eq{l8eq4} in $C^0(S)$ if $a=0$.
Also $u_{\bs\al},v_{\bs\al}$ depend continuously on ${\bs\al}\in U$
in $C^0(S)$, by Theorem~\ref{l8thm1}.

For each ${\bs\al}=(a,b,c)$ in $U$, define $N_{\bs\al}$ in $\C^3$ by
\e
\begin{split}
N_{\bs\al}=\bigl\{(z_1,z_2,z_3)\in\C^3:\,&
z_1z_2=v_{\bs\al}(x,y)+iy,\quad z_3=x+iu_{\bs\al}(x,y),\\
&\ms{z_1}-\ms{z_2}=2a,\quad (x,y)\in S^\circ\bigr\}.
\end{split}
\label{l11eq1}
\e
Then $N_{\bs\al}$ is a noncompact SL 3-fold without boundary in $\C^3$,
which is nonsingular if $a\ne 0$, by Proposition~\ref{l8prop1}.
\label{l11def}
\end{dfn}

In \cite[Th.~8.2]{Joyc11} we show that the $N_{\bs\al}$ are the
fibres of an {\it SL fibration}.

\begin{thm} In the situation of Definition \ref{l11def}, if\/
${\bs\al}\ne{\bs\al}'$ in $U$ then $N_{\bs\al}\cap N_{{\bs\al}'}
=\emptyset$. There exists an open set\/ $V\subset\C^3$ and a
continuous, surjective map $F:V\ra U$ such that\/ $F^{-1}({\bs\al})
=N_{\bs\al}$ for all\/ ${\bs\al}\in U$. Thus, $F$ is a special
Lagrangian fibration of\/ $V\subset\C^3$, which may include
singular fibres.
\label{l11thm}
\end{thm}

It is easy to produce families $\Phi$ satisfying Definition
\ref{l11def}. For example \cite[Ex.~8.3]{Joyc11}, given any
$\phi\in C^{3,\al}(\pd S)$ we may define $U=\R^3$ and $\Phi:
\R^3\ra C^{3,\al}(\pd S)$ by $\Phi(a,b,c)=\phi+bx+cy$. So this
construction produces very large families of $\U(1)$-invariant
SL fibrations, including singular fibres, which can have any
multiplicity and type.

Here is a simple, explicit example. Define $F:\C^3\ra\R\t\C$ by
\e
\begin{gathered}
F(z_1,z_2,z_3)=(a,b),\quad\text{where}\quad 2a=\ms{z_1}-\ms{z_2} \\
\text{and}\quad
b=\begin{cases}
z_3, & a=z_1=z_2=0, \\
z_3+\bar z_1\bar z_2/\md{z_1}, & a\ge 0,\;\> z_1\ne 0,\\
z_3+\bar z_1\bar z_2/\md{z_2}, & a<0.
\end{cases}
\end{gathered}
\label{l11eq2}
\e
This is a piecewise-smooth SL fibration of $\C^3$. It is not
smooth on~$\md{z_1}=\md{z_2}$.

The fibres $F^{-1}(a,b)$ are $T^2$-cones singular at $(0,0,b)$ 
when $a=0$, and nonsingular ${\cal S}^1\t\R^2$ when $a\ne 0$. 
They are isomorphic to the SL 3-folds of Example \ref{l8ex1}
under transformations of $\C^3$, but they are assembled to
make a fibration in a novel way.

As $a$ goes from positive to negative the fibres undergo a surgery, 
a Dehn twist on ${\cal S}^1$. The reason why the fibration is only 
piecewise-smooth, rather than smooth, is really this topological 
transition, rather than the singularities themselves. The fibration 
is not differentiable at every point of a singular fibre, rather
than just at singular points, and this is because we are jumping
from one moduli space of SL 3-folds to another at the singular fibres.

I conjecture that $F$ is the local model for codimension one 
singularities of SL fibrations of generic almost Calabi--Yau 3-folds. 
The justification for this is that the $T^2$-cone singularities
have `index one' in the sense of \S\ref{l103}, and so should
occur in codimension one in families of SL 3-folds in generic
almost Calabi--Yau 3-folds. Since they occur in codimension one in 
this family, the singular behaviour should be stable under small 
perturbations of the underlying almost Calabi--Yau structure.

I also have a $\U(1)$-invariant model for codimension two 
singularities, described in \cite{Joyc13}, in which two of the
codimension one $T^2$-cones come together and cancel out. I
conjecture that it too is a typical codimension two singular
behaviour in SL fibrations of generic almost Calabi--Yau 3-folds.
I do not expect codimension three singularities in generic SL
fibrations to be locally $\U(1)$-invariant, and so this
approach will not help.

\subsection{Exercises}
\label{l116}

\begin{question} Let $f:\C^3\ra\R\t\C$ be as in equation~\eq{l11eq2}.
\anext Show that $f$ is continuous, surjective, and piecewise smooth.
\anext Show that $f^{-1}(a,b)$ is a (possibly singular) special 
Lagrangian 3-fold for all~$(a,b)\in\R\t\C$.
\anext Identify the singular fibres and describe their
singularities. Describe the topology of the singular
and the nonsingular fibres.
\end{question}

The idea of a `special Lagrangian fibration' $f:X\ra B$ is in 
some ways a rather unnatural one. One of the problems is that
the map $f$ doesn't satisfy a particularly nice equation, locally; 
the level sets of $f$ do, but the `coordinates' on $B$ are determined
globally rather than locally. To understand the problems with special 
Lagrangian fibrations, try the following (rather difficult) exercise. 

\begin{question}Let $X$ be a Calabi--Yau 3-fold, $N$ a compact SL
3-fold in $X$ diffeomorphic to $T^3$, ${\cal M}_N$ the family of 
special Lagrangian deformations of $N$, and $\,\,\ov{\!\!\cal M}_N$ 
be ${\cal M}_N$ together with the singular SL 3-folds occurring as 
limits of elements of~${\cal M}_N$.

In good cases, SYZ hope that $\,\,\ov{\!\!\cal M}_N$ is the
family of level sets of an SL fibration $f:X\ra B$, where $B$ is 
homeomorphic to $\,\,\ov{\!\!\cal M}_N$. How many different 
ways can you think of for this not to happen? (There are at least 
two mechanisms not involving singular fibres, and others which do). 
\end{question}

\end{document}